\documentclass[11pt,reqno]{amsart}
\usepackage{graphicx,tikz}
\usepackage{amsfonts,amsmath,amssymb}
\usetikzlibrary{matrix,arrows}
\newcommand{\C}{\mathbb C}
\newcommand{\R}{\mathbb R}

\newtheorem{theorem}{Th\'eor\`eme}[section]

\newtheorem{prop}[theorem]{Proposition}
\newtheorem{cor}[theorem]{Corollaire}

\newtheorem{defn}[theorem]{Definition}
\newtheorem{rem}[theorem]{Remarque}
\newtheorem{ex}[theorem]{Exemple}

\begin{document}

\title[]{La formule de Cauchy dans les alg\`ebres}
\date\today

\author{Pierre Bonneau*  and Emmanuel Mazzilli**}

\maketitle

\thanks{\address{(*)Equipe Emile Picard, Institut de Mathématiques, 
 Universit\'{e} Paul Sabatier, 31062 Toulouse Cedex 9 \footnote{\email{\it pierre.bonneau@math.univ-toulouse.fr}}}}.

 \thanks{ \address{(**)Laboratoire Paul Painlev\'e, Universit\'e Nord de France, 59655 Villeneuve-d'ascq \footnote{\email{\it emmanuel.mazzilli@univ-lille.fr}}}}.

\pagestyle{myheadings}
 \markboth{\today}{\today}
 \setcounter{part}{-1}
\vspace{2cm}

\setcounter{section}{-1}
\section{Introduction.}\setcounter{section}{0}
\smallbreak

De nombreux essais de généralisation algébrique de 
$\C$, en dimensions supérieures, eurent lieu, au 19° siècle, mais se heurtèrent à de fortes
limitations, puisque d'après un théorème de 
Bott-Milnor  (\cite{BM}) et Kervaire (\cite{K}),  l'espace vectoriel $\R^{n}$ possède une opération produit $\R-$bilinéaire sans diviseur de $0$ seulement pour $n=1,\;2,\;4\;ou\; 8$ -- et
les seules algèbres de division qui sont associatives sont (à isomorphisme près),
d'après le théorème de Frobenius,   $\R$, $\C$ et l'algèbre
$\mathbb{H}$ des quaternions. Nous n'adopterons pas ce point de vue algébrique, ici, mais allons plutôt  tenter d'obtenir des généralisations analytiques de $\C.$ \\

En effet, le succès de $\mathbb{C}$
repose, pour beaucoup, sur les propriétés remarquables des
fonctions de $\mathbb{C}$ dans $\mathbb{C}$ qui sont
$\mathbb{C}$-différentiables (les fonctions holomorphes). C'est l'analyse complexe. Deux façons de l'aborder sont traditionnelles ; ce sont les points de vue de Weierstrass d'une part, de Cauchy et Riemann d'autre part. C'est ce dernier point de vue que nous essaierons ici de généraliser. L'outil principal en est la représentation intégrale des fonctions holomorphes. Cet outil permet une approche particulièrement féconde de l'analyse complexe, dans \cite{R}, par exemple. Cauchy et Riemann caractérisent les fonctions holomorphes $f$ par certaines conditions, les conditions de Cauchy-Riemann, qui se traduisent par l'équation $\overline{\partial} f\;=\;0$ où $\overline{\partial}$ est l'opérateur de Cauchy-Riemann dont les fonctions holomorphes constituent ainsi le noyau. Une solution fondamentale de l'opérateur $\overline{\partial}$ fournit alors une représentation intégrale de la fonction holomorphe $f$ à partir de ses valeurs sur le bord du domaine de définition $D,$ c'est la formule (de représentation intégrale) de Cauchy. On déduit de cette formule de représentation intégrale des propriétés de la fonction $f$ elle même. \\ 
Nous avions déjà, dans le cadre de la superanalyse, étudié ainsi, pour une superalgèbre $\Lambda$ les fonctions superdifférentiables (voir \cite{BC1} et \cite{BC2}). \\
Nous nous proposons, ici, de revenir sur cette question de la représentation intégrale pour des fonctions à valeurs dans une algèbre, fonctions auxquelles nous imposerons certaines conditions du premier ordre que nous appellerons de Cauchy et qui caractériseront le noyau d'un opérateur différentiel que nous dirons de Cauchy-Riemann. Nous chercherons dès lors une solution fondamentale de cet opérateur, donc une représentation intégrale des fonctions vérifiant ces conditions de Cauchy. Nous serons ainsi en mesure de démontrer certaines propriétés intéressantes pour ces fonctions.\\ 
Le premier paragraphe nous permettra de tester notre méthode pour en déterminer les contraintes, dans un cas trivial et sans intérêt. \\Dans le deuxième paragraphe, nous étudierons le cas de conditions de Cauchy-Riemann à coefficients constants. Après avoir considéré le cas général, nous nous intéresserons au cas d'une algèbre commutative et utiliserons les déterminants à valeurs dans $A.$
 Dans ce second paragraphe, nous donnerons des exemples précis, dans diverses algèbres de dimension 2 ou 4, de classes de fonctions admettant une formule de Cauchy, et démontrerons l'impossibilité d'une telle formule en dimension 3. \\
 Dans le troisième paragraphe, nous expliquerons le cas de l'alg\`ebre des quaternions muni des \'equations de Cauchy-Fueter qui avait d\'ej\`a \'et\'e trait\'e dans \cite{WW1} et \cite{WW2}. Enfin, nous d\'emontrerons les applications classiques dans les alg\`ebres munies d' \'equations de Cauchy-Riemann admettant des formules de Cauchy pour les fonctions v\'erifiant ces \'equations, le th\'eor\`eme de Hartogs ainsi que les solutions \`a support compact pour les donn\'ees \`a support compact.\\
 Dans le dernier paragraphe, nous étudierons, succinctement, le cas où les conditions de Cauchy-Riemann ont des coefficients qui sont des fonctions. Comme au paragraphe précédent, nous envisagerons les cas commutatif et non commutatif.\\
Comme dans \cite{BC1} et \cite{BC2}, nous obtiendrons des contraintes sur l'algèbre $A$ utilisée, et, comme nous l'avions déjà expliqué alors, cela est bien dû à la nature du problème étudié. Si le système des conditions de Cauchy n'est pas elliptique, nous ne pouvons obtenir de formules de représentation intégrale pour la classe des fonctions vérifiant ces conditions. \\

Pour terminer, nous remarquerons, d'abord, que nous définissons, dans cet article, la singularité du noyau de la formule de Cauchy à l'aide d'une norme euclidienne. Ceci, bien sûr, n'est pas obligatoire, et nous pourrions tout aussi bien utiliser une norme non euclidienne, par exemple la jauge induite par un voisinage de 0, à frontière lisse, mais éventuellement non analytique.\\ 
Par ailleurs, nous pouvons aussi remarquer que l'outil essentiel pour établir une formule de Cauchy est le théorème de Stokes. Or il s'agit là d'un outil d'analyse réelle. Cela nous fait présumer que le cadre naturel pour l'étude de classes de fonctions vérifiant une formule de Cauchy n'est pas l'étude de fonctions à valeurs dans une algèbre, mais de fonctions de $\R^{n}$ dans $\R^{p}.$ Nous n'abordons pas cela ici, mais ce sera l'objet d'un prochain article. 
Enfin, nous renvoyons \`a un article en pr\'eparation pour les formules de repr\'esentation des formes et non plus uniquement des fonctions, autrement dit l'\'equivalent de la formule de Koppelmann en analyse pluri-complexe.

\section{Le cas des fonctions réelles de variables réelles.}
Notre objectif est donc de généraliser le point de vue de Cauchy-Riemann. Pour étudier l'apport de la méthode proposée, mais aussi les limitations, nous allons d'abord considérer le cas le plus naif, celui où l'algèbre considérée est $\R$ et où la singularité du noyau intégral est fournie par une norme euclidienne. Comme une algèbre de dimension finie $n$ est isomorphe, en tant qu'espace vectoriel, à $\mathbb{R}^{n}$, les composantes d'une fonction définie sur une puissance d'une algèbre $A$ et à valeurs dans $A$ relèvent de ce cas. \\

Nous considérons une fonction définie sur un ouvert $D$ de $\R^{n}$ à bords lisses, $f\,:\,D\longrightarrow \R$ dont les variables seront notées $(x_{1},...,x_{n})$, fonction que nous supposons, a priori, de classe $C^{1}$ (nous verrons, plus loin, que cette hypothèse est inutile) et à laquelle nous imposons $q$ conditions de Cauchy-Riemann d'ordre 1:
$$P_{m}(\frac{\partial }{\partial x})(f)\;:=\;\sum^{n}_{j=1}\;a^{j}_{m}\;\frac{\partial f}{\partial x_{j}}\;=\;0;\quad m=1,...,q;\quad a^{j}_{m}\in \R.$$
Nous définissons alors l'opérateur de Cauchy-Riemann $d''$ par: 
$$d''\;=\;\sum^{q}_{m=1}\,dx_{m}\sum^{n}_{j=1}\,a^{j}_{m}\,\frac{\partial}{\partial x_{j}}$$ 
de sorte que $f$ est dans le noyau de $d''$ si et seulement si elle satisfait les conditions de Cauchy-Riemann.\\ 
Si $D$ est un domaine de $\R^{n}$ borné et à
frontière lisse, si $f$ est une fonction continue dans $\overline{D},$
 de classe $C^{1}$ dans $D$ et à valeurs réelles, nous cherchons une représentation
intégrale de $f$, c'est-à-dire nous essayons d'écrire
\begin{eqnarray}\label{RepInt}
f(x)=\int_{\partial D}\,f(y)K(y,x)\,-\,\int_{D}\,d''f(y)K(y,x),\;\;x\in D
\end{eqnarray}
Cette représentation intégrale peut s'écrire en termes de courants
$d_{y}''K(y,x)=[\Delta]_{n-1}$ où  $[\Delta]$ est le courant d'intégration
sur la diagonale $\Delta=\{(x,x):\,x\in D\}$ et $[\Delta]_{n-1}$ est la composante de degré $n-1$ en $y$ de $[\Delta]$ (voir \cite{HP}). Un noyau
$K(y,x)$ vérifiant cette dernière égalité est lié à une  solution
fondamentale de l'opérateur de Cauchy-Riemann $d''$. \\
Soit $\Psi:\,\overline{D}\times\overline{D}\rightarrow\R^{n}$
définie par $\Psi(y,x)=y-x$. Alors $[\Delta]=\Psi^{\ast}([0])$ où
$[0]$ est le courant d'évaluation en
$0\in\R^{n}$. Donc $d''K=[\Delta]$ s'écrit
aussi $d''\Omega=[0]$ si l'on pose $K=\Psi^{\ast}\Omega$ (le
premier $d''K$ porte sur les deux variables (x,y),
$d''=d_{y}''+d_{x}''$, tandis que le second $d''\Omega$ porte sur
la variable $x'=y-x$ que nous noterons x). \\
Nous devons donc chercher
 $\Omega$ vérifiant $d''\Omega\,=\,[0]$ i.e. une solution fondamentale de $d''.$ \\
  Notons $dV=dx_{1}\wedge dx_{2}\wedge ... \wedge dx_{n},$ et $\widehat{dx_{m}}=dx_{1}\wedge dx_{2}\wedge\cdots\wedge dx_{m-1}\wedge dx_{m+1}\wedge ... \wedge dx_{n}.$
 Nous cherchons cette solution fondamentale $\Omega$ sous la forme: 
\begin{equation}\label{Fond}
\Omega\;=\;\frac{1}{\left\|x\right\|^{n}}\sum^{q}_{m=1}\,(-1)^{m-1}\widehat{dx_{m}}\sum^{n}_{j=1}b^{j}_{m}x_{j}\;=\;\frac{1}{\left\|x\right\|^{n}}\sum^{q}_{m=1}\,(-1)^{m-1}\widehat{dx_{m}}Q_{m}(x)\;=\;\sum^{q}_{m=1}(-1)^{m-1}\widehat{dx_{m}}\Omega_{m}.
\end{equation}
Hors de 0, $\Omega$ est lisse, et, par un calcul direct: 
\begin{equation}
\begin{split}
d''\Omega\;&=\;\frac{1}{\left\|x\right\|^{n+2}}\left[\left\|x\right\|^{2}\sum^{q}_{m=1}\sum^{n}_{j=1}\,a^{j}_{m}b^{j}_{m}\;-\;n\sum^{q}_{m=1}\sum^{n}_{j,j'=1}\,a^{j}_{m}b^{j'}_{m}x_{j}x_{j'}\right]\\
&=\frac{1}{\left\|x\right\|^{n+2}}\left[\left\|x\right\|^{2}\sum_{m,j}a^{j}_{m}b^{j}_{m}\;-\;n\sum_{m,j}a^{j}_{m}b^{j}_{m}x^{2}_{j}\;-\;n\sum_{m,j<j'}\,(a^{j}_{m}b^{j'}_{m}\,+\,a^{j'}_{m}b^{j}_{m})x_{j}x_{j'}\right].
\end{split}
\end{equation}
La nullité, hors de $0$, de $d''\Omega$ s'écrit:
\begin{equation}
\begin{cases}
\sum^{q}_{m=1}\;a^{j}_{m}b^{j}_{m}\;=\;\frac{1}{n}\sum^{q}_{m=1}\sum^{n}_{j=1}a^{j}_{m}b^{j}_{m};\quad \forall j=1,...,n.\\
\sum^{q}_{m=1}\;(a^{j}_{m}b^{j'}_{m}\;+\;a^{j'}_{m}b^{j}_{m})\;=\;0;\quad j\neq j'.
\end{cases}
\end{equation}
La première ligne dit, tout simplement, qu'il existe une constante $C$ telle que, pour tout $j=1,...,n,$ on ait $\sum^{q}_{m=1}\;a^{j}_{m}b^{j}_{m}\;=\;C.$ Le système précédent s'écrit donc: 
\begin{equation}\label{Cond}
\sum^{q}_{m=1}\;(a^{j}_{m}b^{j'}_{m}\;+\;a^{j'}_{m}b^{j}_{m})\;=\;2C\delta^{j'}_{j}.
\end{equation}
Il reste maintenant à étudier la singularité en 0. De façon tout à fait analogue au lemme 1 de \cite{BC2}, nous démontrons que si $\varphi$ est une fonction test, alors $<d''\Omega ,\varphi>\;=\;nC\varphi (0).Vol(B(0,1)).$ Si nous choisissons $C\,=\,[nVol(B(0,1))]^{-1},$ nous avons pour $\Omega$ une solution fondamentale de $d''$ et donc une formule de représentation intégrale \ref{RepInt}.\\
La condition \ref{Cond} s'écrit aussi, en notant $Q_{m}(x)$ le polynôme $Q_{m}(x)\,=\,\sum^{n}_{j=1}b^{j}_{m}x_{j},$ 
\begin{equation}\label{C1a}
\sum^{q}_{m=1}\,P_{m}(x)Q_{m}(x)\;=\;C\sum^{n}_{j=1}\,x^{2}_{j},
\end{equation}
ce qui prouve que l'opérateur $d''$ est elliptique s'il admet une solution fondamentale du type \ref{Fond}. En effet, $P_{1}(x)\,=\,...\,=\,P_{q}(x)\,=\,0$ implique que $x=0.$ 
Comme nous l'avons annoncé en début de paragraphe, la condition de régularité $C^{1}$ faite sur $f$ est donc superflue.\\
 Par ailleurs, la condition \ref{C1a} impose que les formes linéaires $P_{1}(x),...,P_{q}(x)$ dont l'intersection des noyaux se réduit à $\{0\}$ sont au moins au nombre de $n,$ donc $q\geq n,$ et, parmi elles, $n$ sont indépendantes (par exemple, et sans particulariser $P_{1},...,P_{n}$ qui forment donc une base du dual de $\R^{n},$) tandis que les autres sont des combinaisons linéaires de ces $n$ premières. Notons $U_{1},...,U_{n}$ les vecteurs suivants de $\R^{n}$ écrits dans la base canonique $(e_{1},...,e_{n})$ de $\R^{n}$
$$U_{m}\,=\,\sum^{n}_{j=1}a^{j}_{m}e_{j}.$$
Ces $n$ vecteurs de $\R^{n}$ sont indépendants, puisque $P_{1},...,P_{n}$ le sont, et les conditions $P_{m}(\frac{\partial }{\partial x})(f)\,=0$ s'écrivent $\frac{\partial f}{\partial U_{m}}=0.$ Ainsi, les $n$ conditions imposées à $f$ traduisent le fait que $df\,=\,0,$ donc que $f$ est constante. \\
Ainsi, si $d''f\,=\,0,$ et si $d''$ admet une solution fondamentale, alors $f$ est une constante.\\
Nous déduisons, des calculs ci-dessus, trois enseignements.\\
Tout d'abord, plus nous imposons de conditions de Cauchy-Riemann, plus il est aisé d'obtenir une solution fondamentale; mais moins il existe de fonctions satisfaisant les conditions imposées. Il s'agit donc de maintenir une sorte d'équilibre entre la quantité de fonctions vérifiant les conditions de Cauchy-Riemann et l'existence d'une solution fondamentale pour l'opérateur de Cauchy-Riemann $d''$ associé.\\ 
Ensuite, si nous voulons obtenir une représentation intégrale de la forme \ref{RepInt} pour une fonction de $D$ dans $\R^{p},$ alors le noyau $K(y,x)$ ne peut être réel. Sinon, la représentation intégrale s'effectue composante par composante, et donc, pour chaque composante, nous sommes ramenés au cas précédent, avec la limitation $d''f=0$ impliquant que $f$ est constante.\\
Enfin, comme nous avons besoin, dans \ref{RepInt}, du produit de $f(y)$ par $K(y,x),$ il faut définir dans l'espace d'arrivée $\R^{p}$ un produit bilinéaire, donc une structure d'algèbre, et cette structure ne peut se résumer au produit composante par composante, sous peine de revenir aux limitations précédentes, et pour les mêmes raisons.

\section{Conditions de Cauchy à coefficients constants.}
Soit une fonction définie sur l'adhérence d'un domaine $D$ de $\R^{n}$ à bords lisses, $f\,:\,\overline{D}\longrightarrow A$ dont les variables seront notées $(x_{1},...,x_{n})$, fonction que nous supposons, a priori, de classe $C^{1}$, à valeurs dans l'algèbre unitaire
$A\,=\,Vect(e_{0},e_{1},...,e_{p})$, $e_{0}$ étant l'élément unité et le produit étant défini, dans la base $(e_{0},e_{1},...,e_{p}),$ par: 
$$e_{i}e_{j}\;=\;\sum^{p}_{k=0}\Gamma^{k}_{i,j}e_{k},\quad \quad \Gamma^{k}_{i,j}\in \R.$$

Nous allons, dans ce paragraphe, imposer à $f,$ q conditions de Cauchy-Riemann à coefficients constants dans $A$ et d'ordre 1: 
\begin{equation}\label{C2a}
P_{m}(\frac{\partial }{\partial x})(f)\;:=\;\sum^{n}_{j=1}\;\frac{\partial f}{\partial x_{j}}\;a^{j}_{m}\;=\;0;\quad m=1,...,q;\quad a^{j}_{m}\in A.
\end{equation}
Nous définissons, par conséquent, les opérateurs de Cauchy-Riemann par: 
$$d''\;=\;\sum^{q}_{m=1}dx_{i_{m}}\sum^{n}_{j=1}\;\frac{\partial}{\partial x_{j}}a^{j}_{m}, \quad 1\leq i_{m}<i_{m+1}\leq n,$$
$$\widehat{d''}\;=\;\sum^{q}_{m=1}dx_{i_{m}}\sum^{n}_{j=1}\;a^{j}_{m}\frac{\partial}{\partial x_{j}}, \quad 1\leq i_{m}<i_{m+1}\leq n,$$
de façon que $f$ soit dans le noyau de $d''$ si et seulement il satisfait les $q$ conditions de Cauchy-Riemann. \\
Ces deux opérateurs $d''$ et $\widehat{d''}$ se confondent dans le cas d'une algèbre $A$ commutative, mais nous les explicitons précisément pour être en mesure de traiter le cas non commutatif.\\
Nous cherchons une représentation intégrale de la forme \ref{RepInt} pour $f$, c'est-à-dire une solution fondamentale pour $\widehat{d''_{y}}$ de la forme: 
\begin{equation}
\begin{split}
K(y,x)\;&=\;\frac{1}{\left\|y-x\right\|^{n}}\sum^{q}_{m=1}(-1)^{(i_{m}-1)}\widehat{dy_{i_{m}}}\varphi_{m}(x,y)\\
&=\;\sum(-1)^{(i_{m}-1)}\widehat{dy_{i_{m}}}\Omega_{m}.
\end{split}
\end{equation}
où $\varphi_{m}$ est une fonction à valeurs dans $A,$ et au moins $C^{1}$ qui s'annule sur $\Delta.$\\
Hors de $\Delta,$ le calcul de $\widehat{d''_{y}}K(y,x)$ est tout à fait analogue au calcul correspondant du paragraphe 1, en gardant toutefois à l'esprit que $A$ peut n'être pas commutatif. On note $dV\,=\,dx_{1}\wedge dx_{2}\wedge ...\wedge dx_{n}.$
\begin{equation}
\begin{split}
\widehat{d''_{y}}K(y,x)\;&=\;\frac{1}{\left\|y-x\right\|^{n}}\sum^{q}_{m=1}\sum^{p}_{j=0}\,a^{j}_{m}\frac{\partial \varphi_{m}(x,y)}{\partial y_{j}}dV\,-\,\frac{n}{2\left\|y-x\right\|^{n+2}}\sum^{q}_{m=1}\,dy_{i_{m}}\sum_{j}\,a^{j}_{m}\frac{\partial\left\|y-x\right\|^{2}}{\partial y_{j}}\\
&\qquad \varphi_{m}(x,y)dV\\
&=\;\frac{1}{\left\|y-x\right\|^{n+2}}\big[\left\|y-x\right\|^{2}\sum^{q}_{m=1}\sum^{p}_{j=0}\,a^{j}_{m}\frac{\partial\varphi_{m}}{\partial y_{j}}\,-\,\frac{n}{2}\sum_{m,j}a^{j}_{m}\frac{\partial\left\|y-x\right\|^{2}}{\partial y_{j}}\varphi_{m}(x,y)\big]dV\,
\end{split}
\end{equation}
de sorte que la nullité de $\widehat{d''}K(y,x)$ hors de $\Delta$ s'écrit:
\begin{equation}\label{cond1}
\left\|y-x\right\|^{2}\sum_{m,j}a^{j}_{m}\frac{\partial\varphi_{m}}{\partial x_{j}}\,=\,n\sum_{m,j}a^{j}_{m}(y_{j}-x_{j})\varphi_{m}(x,y).
\end{equation}
Nous allons maintenant étudier la singularité de $\widehat{d''}K(y,x),$ et, comme $A$ n'est pas commutatif, nous allons détailler un peu ce calcul.  
\begin{equation}
\begin{split}
\int_{D}d''f(y)\wedge K(y,x)\;&=\;\int_{D}\frac{1}{\left\|y-x\right\|^{n}}\sum^{q}_{m=1}\sum^{n}_{j=1}\frac{\partial f}{\partial y_{j}}(y)a^{j}_{m}\varphi_{m}(x,y)dV\\
&=\;\lim_{\epsilon\rightarrow 0}\int_{D-B(x,\epsilon)}\frac{1}{\left\|y-x\right\|^{n}}\sum^{q}_{m=1}\sum^{n}_{j=1}\frac{\partial f}{\partial y_{j}}(y)a^{j}_{m}\varphi_{m}(x,y)dV\\
&=\;\lim_{\epsilon\rightarrow 0}\int_{D-B(x,\epsilon)}d\Big[\sum_{m,j}f(y)a^{j}_{m}\frac{\varphi_{m}(x,y)}{\left\|y-x\right\|^{n}}(-1)^{j-1}\widehat{dy_{j}}\Big]-f(y)\widehat{d''_{y}}K(y,x)\\
&=\;\int_{\partial D}f(y)\sum_{m,j}a^{j}_{m}\frac{\varphi_{m}(x,y)}{\left\|y-x\right\|^{n}}(-1)^{j-1}\widehat{dy_{j}}\\
&\qquad\qquad-\lim_{\epsilon\rightarrow 0}\int_{\partial B(x,\epsilon)}f(y)\sum_{m,j}a^{j}_{m}\frac{\varphi_{m}(x,y)}{\left\|y-x\right\|^{n}}(-1)^{j-1}\widehat{dy_{j}}\\
&=\;\;\int_{\partial D}f(y)\sum_{m,j}a^{j}_{m}\frac{\varphi_{m}(x,y)}{\left\|y-x\right\|^{n}}(-1)^{j-1}\widehat{dy_{j}}\\
&\qquad\qquad-\lim_{\epsilon\rightarrow 0}\int_{\partial B(x,\epsilon)}\sum_{m,j}\Big(f(x)+O(\epsilon)\Big)a^{j}_{m}\frac{O(\epsilon)}{\epsilon^{n}}(-1)^{j-1}\widehat{dy_{j}}\\
&=\;\int_{\partial D}f(y)\sum_{m,j}a^{j}_{m}\frac{\varphi_{m}(x,y)}{\left\|y-x\right\|^{n}}(-1)^{j-1}\widehat{dy_{j}}\\
&\qquad\qquad-f(x)\lim_{\epsilon\rightarrow 0}\int_{\partial B(x,\epsilon)}\sum_{m,j}\frac{a^{j}_{m}\varphi_{m}(x,y)(-1)^{j-1}}{\epsilon^{n}}\widehat{dy_{j}}\\
&=\;\int_{\partial D}f(y)\sum_{m,j}a^{j}_{m}\frac{\varphi_{m}(x,y)}{\left\|y-x\right\|^{n}}(-1)^{j-1}\widehat{dy_{j}}\\
&\qquad\qquad-f(x)\lim_{\epsilon\rightarrow 0}\frac{1}{\epsilon^{n}}\int_{B(x,\epsilon)}\sum_{m,j}a^{j}_{m}\frac{\partial \varphi_{m}(x,y)}{\partial y_{j}}dV\\
&=\;\int_{\partial D}f(y)\sum_{m,j}a^{j}_{m}\frac{\varphi_{m}(x,y)}{\left\|y-x\right\|^{n}}(-1)^{j-1}\widehat{dy_{j}}\\
&\qquad\qquad-f(x)\lim_{\epsilon\rightarrow 0}\frac{1}{\epsilon^{n}}\int_{B(x,\epsilon)}\sum_{m,j}a^{j}_{m}\Big(\frac{\partial \varphi_{m}(x,x)}{\partial y_{j}}+O(\epsilon)\Big)dV\\
&=\;\int_{\partial D}f(y)\sum^{q}_{m=1}\Big(\sum^{n}_{j=1}a^{j}_{m}\frac{\partial}{\partial y_{j}}\Big)\lrcorner \Big(\frac{\varphi_{m}(x,y)}{\left\|y-x\right\|^{n}}dV\Big)\\
&\qquad\qquad-f(x)Vol(B(0,1))\sum^{q}_{m=1}\sum^{n}_{j=1}a^{j}_{m}\frac{\partial \varphi_{m}(x,x)}{\partial y_{j}}. 
\end{split}
\end{equation}
On veut que le dernier terme, à la dernière ligne,  soit $f(x).$ Tenant compte de \ref{cond1} , pour que $K(y,x)$ soit une solution   fondamentale, nous avons donc :
\begin{equation}\label{cond2}
\begin{cases}
\sum^{q}_{m=1}\,\sum^{n}_{j=1}a^{j}_{m}\frac{\partial\varphi_{m}}{\partial y_{j}}(x,x)\;=\;\frac{e_{0}}{Vol(B(0,1))},\\
\left\|y-x\right\|^{2}\sum_{m,j}a^{j}_{m}\frac{\partial\varphi_{m}}{\partial y_{j}}\,=\,n\sum_{m,j}a^{j}_{m}(y_{j}-x_{j})\varphi_{m}(x,y), 
\end{cases}
\end{equation}
et nous obtenons la formule de représentation intégrale 
\begin{equation}\label{repint}
f(x)\;=\;\int_{\partial D}f(y)\sum^{q}_{m=1}\Big(\sum^{n}_{j=1}a^{j}_{m}\frac{\partial}{\partial y_{j}}\Big)\lrcorner \Big(\frac{\varphi_{m}(x,y)}{\left\|y-x\right\|^{n}}dV\Big)-\int_{D}d''f(y)\wedge K(y,x). 
\end{equation}
Si $d''f(y)=0,$ alors le dernier terme est nul. Nous obtenons, par conséquent, la formule de Cauchy 
\begin{equation}\label{CondCau}
f(x)\;=\;\int_{\partial D}f(y)\sum^{q}_{m=1}\Big(\sum^{n}_{j=1}a^{j}_{m}\frac{\partial}{\partial y_{j}}\Big)\lrcorner \Big(\frac{\varphi_{m}(x,y)}{\left\|y-x\right\|^{n}}dV\Big) 
\end{equation}
pour la classe des fonctions vérifiant $d''f=0,$ c'est-à-dire la classe des fonctions vérifiant les conditions de Cauchy \ref{C2a}.\\
Au vu des conditions \ref{cond2}, il semble naturel de chercher les fonctions $\varphi_{m}$ sous la forme de polynômes homogènes de degré 1. C'est pourquoi nous posons: 
\begin{equation}\label{phi}
\varphi_{m}(x,y)\;=\;\sum^{n}_{i=1}\,b^{i}_{m}(y_{i}-x_{i}), \quad b^{i}_{m}\in A. 
\end{equation}
Les conditions \ref{cond2} s'écrivent alors: 
\begin{equation}\label{cond11}
\begin{cases}
\sum^{q}_{m=1}\sum^{n}_{j=1}\,a^{j}_{m}b^{j}_{m}\;=\;\frac{e_{0}}{Vol(B(0,1))},\\
\left\|y-x\right\|^{2}\sum_{m,j}a^{j}_{m}b^{j}_{m}\;=\;n\sum_{m,j}a^{j}_{m}(y_{j}-x_{j})\sum^{n}_{i=1}b^{i}_{m}(y_{i}-x_{i}).
\end{cases}
\end{equation}
Le terme de droite dans la dernière équation s'écrit:\\ $n\sum^{n}_{j=1}(y_{j}-x_{j})^{2}\sum^{q}_{m=1}a^{j}_{m}b^{j}_{m}\;+\;n\sum_{j<i}(y_{j}-x_{j})(y_{i}-x_{i})\sum^{q}_{m=1}(a^{j}_{m}b^{i}_{m}\,+\,a^{i}_{m}b^{j}_{m}),$ \\
de sorte qu'en identifiant les membres de gauche et de droite, nous obtenons: 
\begin{equation}
\begin{cases}\label{cond12}
\sum^{q}_{m=1}\,a^{j}_{m}b^{j}_{m}\;=\;\frac{e_{0}}{nVol(B(0,1))},\quad pour \quad tout \quad j=1,...,n;\\
\sum^{q}_{m=1}a^{j}_{m}b^{i}_{m}\,+\,a^{i}_{m}b^{j}_{m}\;=\;0,\quad pour \quad i\neq j.
\end{cases}
\end{equation}
Si nous appelons $\widetilde{A}$ (respectivement B) les matrices à $n$ lignes et $q$ colonnes comportant comme entrée, à la ligne $j$ et à la colonne $m$, l'élément $a^{j}_{m}$ (respectivement $b^{j}_{m}$) de l'algèbre $A$, et $I$ la matrice carrée d'ordre $n$ comportant $e_{0}$ sur chaque entrée de la diagonale principale, et $0$ partout ailleurs, alors les conditions \ref{cond12} s'écrivent: 
\begin{equation}
\widetilde{A}\,^{t}B\;+\;^{t}(\widetilde{A}\,^{t}B)\;=\;\frac{2}{nVol(B(0,1))}I,
\end{equation}
ou encore: 
\begin{equation}\label{condmat}
\widetilde{A}\,^{t}B\;=\;\frac{1}{nVol(B(0,1))}I\;+\;\mathbb{A}
\end{equation}
où $\mathbb{A}$ est une matrice antisymétrique à coefficients dans l'algèbre $A.$\\
\\
Une autre façon de traduire les conditions \ref{cond12} est de définir les polynômes $P_{m}(X)$ et $Q_{m}(X)$ à coefficients dans l'algèbre $A$ et variable $X=(X_{1},...,X_{n})\in \R^{n}$ par 
$$P_{m}(X)=\sum^{n}_{j=1}a^{j}_{m}X_{j}$$
et $$Q_{m}(X)=\sum^{n}_{j=1}b^{j}_{m}X_{j}.$$
Les conditions \ref{cond12} deviennent alors
\begin{equation}\label{condell}
\sum^{q}_{m=1}P_{m}(X)Q_{m}(X)=\frac{e_{0}}{nVol(B(0,1))}\left\|X\right\|^{2}.
\end{equation}
Ainsi, $P_{1}(X)=P_{2}(X)=...=P_{q}(X)=0$ implique $X=0.$\\

\begin{rem}\label{rq3}
Ellipticité des conditions de Cauchy.
\end{rem}
Les conditions \ref{cond12} nécessaires et suffisantes pour l'obtention d'une formule de Cauchy impliquent donc l'ellipticité des ces conditions. \\
Bien sûr, nous pourrions aussi traduire la nullité de chaque polynôme $P_{i}(X)$ par la nullité de toutes ses composantes si nous voulions travailler avec des polynômes à coefficients réels.\\
\\
\begin{rem}
	Dans ces calculs, l'algèbre $A$ peut ne pas être associative.\\ 
\end{rem}

Nous allons encore traduire d'une autre façon les conditions \ref{cond12}.\\
Puisque les fonctions $\varphi_{m}$ sont des polynômes homogènes du premier degré, nous posons, pour tout  $j=1,...,n,$
$$\Psi^{j}(x,y)\;=\;Vol(B(0,1))\sum^{q}_{m=1}a^{j}_{m}\varphi_{m}(x,y)\;=\;\sum^{n}_{i=1}\,c^{j}_{i}(y_{i}-x_{i}),$$ 
de sorte que les conditions \ref{cond2} s'écrivent: 
\begin{equation}\label{cond3}
\begin{cases}
\sum^{n}_{j=1}\frac{\partial\Psi^{j}}{\partial y_{j}}(x,x)\;=\;e_{0},\\
\left\|y-x\right\|^{2}\sum_{j}\frac{\partial\Psi^{j}}{\partial y_{j}}(x,y)\,=\,n\sum_{j}(y_{j}-x_{j})\Psi^{j}(x,y).
\end{cases}
\end{equation}
Dans ce système, toute référence aux conditions de Cauchy a disparu. Il est aisé de le résoudre car la deuxième équation s'écrit: 
$$\left\|y-x\right\|^{2}\sum^{n}_{j=1}\,c^{j}_{j}\;=\;n\sum^{n}_{i=1}\,c^{i}_{i}(y_{i}-x_{i})^{2}\;+\;n\sum_{i<j}(y_{i}-x_{i})(y_{j}-x_{j})(c^{j}_{i}\,+\,c^{i}_{j})$$
donc, en identifiant les deux membres des deux équations, nous obtenons: 
\begin{equation}
\begin{cases}
\sum^{n}_{j=1}\,c^{j}_{j}\;=\;nc^{i}_{i}\;=\;e_{0},\\
c^{i}_{j}\;+\;c^{j}_{i}\;=\;0, \quad si\quad i\neq j,
\end{cases}
\end{equation}
soit aussi: 
\begin{equation}\label{cond4}
\begin{cases}
c^{i}_{i}\;=\;\frac{e_{0}}{n}\quad pour \quad tout \quad i=1,...,n,\\
c^{i}_{j}\;+\;c^{j}_{i}\;=\;0, \quad si\quad i\neq j.
\end{cases}
\end{equation}
Nous avons ainsi
$$\Psi^{j}(x,y)\;=\;(y_{j}-x_{j})\frac{e_{0}}{n}\;+\;\sum_{i\neq j}\,(y_{i}-x_{i})c^{j}_{i}\quad avec \quad  c^{j}_{i}\;=\;-c^{i}_{j}.$$
Il reste maintenant à résoudre le système de $n$ équations à $q$ inconnues $ \varphi_{m}(x,y)$
\begin{equation}\label{syst1}
\forall j=1,...,n,\quad\frac{\Psi^{j}(x,y)}{Vol(B(0,1))}\;=\;\sum^{q}_{m=1}\,a^{j}_{m}\varphi_{m}(x,y).
\end{equation}
soit, en explicitant,

\begin{equation}
\frac{\;(y_{j}-x_{j})\frac{e_{0}}{n}\;+\;\sum_{i\neq j}\,(y_{i}-x_{i})c^{j}_{i}}{Vol(B(0,1))}\;=\;\sum^{q}_{m=1}\,a^{j}_{m}\sum^{n}_{i=1}\,b^{i}_{m}(y_{i}-x_{i}),
\end{equation}
ou encore
\begin{equation}\label{syst4}
\begin{cases}
\sum_{m=1}^{q}a_{m}^{j}b_{m}^{j}=\frac{e_{0}}{nVol}\\
\sum_{m=1}^{q}a_{m}^{j}b_{m}^{i}=\frac{c^{j}_{i}}{Vol}\;\; si\;\; i\neq j,
\end{cases}
\end{equation}
qui, pour chaque $i=1,...,n$ et un système à $n$ équations (dans $A$) et $q$ inconnues (dans $A$).
Si nous arrivons à résoudre ce système, la connaissance des fonctions $\varphi_{m}$ permet d'expliciter le noyau $K(y,x),$ et d'écrire la formule de représentation intégrale \ref{repint}. Pour la classe des fonctions vérifiant \ref{C2a}, nous disposons, alors, de la formule de Cauchy \ref{CondCau}.\\

Nous recherchons donc à quelles conditions ce système a des solutions.\\
Il s'écrit, en notant $A^{j}_{m,l}=a^{j}_{m}e_{l}=\sum_{s=0}^{p}A^{j,s}_{m,l}e_{s},$
\begin{equation}
\begin{split}
\sum_{m=1}^{q}a^{j}_{m}&\sum_{i=1}^{n}b^{i}_{m}(y_{i}-x_{i})=\sum_{m=1}^{q}a^{j}_{m}\sum_{l=0}^{p}\sum_{i=1}^{n}b^{i,l}_{m}(y_{i}-x_{i})e_{l}\\
&\qquad =\sum_{m=1}^{q}\sum_{i=1}^{n}\sum_{l=0}^{p}b^{i,l}_{m}(y_{i}-x_{i})A^{j}_{m,l}=\sum_{m=1}^{q}\sum_{i=1}^{n}\sum_{l=0}^{p}b^{i,l}_{m}(y_{i}-x_{i})\sum_{s=0}^{p}A^{j,s}_{m,l}e_{s}\\
&\quad =(y_{j}-x_{j})\frac{e_{0}}{nVol(B(0,1))}+\sum_{i\neq j}(y_{i}-x_{i})c^{j}_{i}\quad avec \;\; c^{i}_{j}=-c^{j}_{i}\\
&\qquad = (y_{j}-x_{j})\frac{e_{0}}{nVol(B(0,1))}+\sum_{i\neq j}(y_{i}-x_{i})\sum_{s=0}c^{j,s}_{i}e_{s},
\end{split}
\end{equation}
donc, en identifiant les composantes $\forall i,j=1,,n\;\; avec\;\;i\neq j ,\; \forall s=0,...,p,$

\begin{equation}\label{S4}
\begin{cases}
\sum_{m=1}^{q}\sum_{l=0}^{p}b^{j,l}_{m}A^{j,s}_{m,l}=\frac{\delta^{s}_{0}}{nVol}:=c^{j,s}_{j}\\
\sum_{m=1}^{q}\sum_{l=0}^{p}b^{i,l}_{m}A^{j,s}_{m,l}=c^{j,s}_{i}.
\end{cases}
\end{equation}
Pour chaque $i=1,...,n,$ nous avons le système 
\begin{equation}\label{Rel4}
S_{i}\;\;:\;\;\begin{cases}
\sum_{m,l}b^{i,l}_{m}A^{j,s}_{m,l}=c^{j,s}_{i}.
\end{cases}
\end{equation}
Pour tous les $i=1,...,n,$ les coefficients du premier membre de $S_{i}$ sont  les mêmes, mais pas les inconnues $b^{i,l}_{m},$ ni les seconds membres. Chaque système $S_{i}$ comporte $n(p+1)$ équations et $q(p+1)$ inconnues.Chaque équation est repérée par son indice $(j,s).$ Nous notons $r$ le rang de chaque système $S_{i}.$ Nous choisissons $r$ équations principales. Ces $r$ équations principales sont repérées par $r$ indices $(j,s)$ que nous classons dans l'ordre alphabétique, ce qui donne les indices $(j_{1},s_{1}),\;\; (j_{2},s_{2}),\;\; ....(j_{r},s_{r}), $ pour les équations principales. Pour les équations non principales, nous aurons les indices\\ 
$(j_{r+1},s_{r+1}),\;\; (j_{r+2},s_{r+2}),\;\; ,....,(j_{n(p+1)},s_{n(p+1)}).$\\
Remarquons que, d'après le classement dans l'ordre alphabétique,
\begin{equation}\label{Rel}
\begin{split}
si \;\; k\leq l\;\; alors\;\;j_{k}\leq j_{l},\\
si\;\; j_{k}<j_{l}\:\: alors\:\: k<l.
\end{split}
\end{equation}
Nous notons
\begin{equation}
\begin{split}
I&=\{(j_{1},s_{1}),\;(j_{2},s_{2}),...., (j_{r},s_{r})\}\\
\overline{I}&=[1,n]\times [0,p]-I=\{(j_{r+1},s_{r+1}),....,(j_{n(p+1)},s_{n(p+1)})\}\\
\forall i=1,...,n,\quad I_{i}&=\{k=1,...,r\;:\;(i,s_{k})\in I\}\\
\forall i=1,...,n,\quad \overline{I}_{i}&=\{k=1,...,r\;:\; (i,s_{k})\in \overline{I}\}
\end{split}
\end{equation}
Le système des équations principales admet un déterminant d'ordre $r$ non nul. Nous le noterons $D_{0}.$\\
\begin{equation}
D_{0}=det\Big( A^{j_{k},s_{k}}_{m_{t},l_{t}}\Big)_{k,t=1,...,r}\neq 0.
\end{equation}
Il détermine, bien sûr, les inconnues principales.\\
Pour que le système complet $S_{i}$ admette comme solutions les solutions du système principal, il faut et il suffit que les $n(p+1)-r$ déterminants caractéristiques s'annulent, ce qui s'écrit, $\forall (j_{0},s_{0})\in \overline{I},$
\begin{equation}
\begin{vmatrix}
A^{j_{1},s_{1}}_{m_{1},l_{1}}&.....&A^{j_{1},s_{1}}_{m_{r},s_{r}}&c^{j_{1},s^{1}}_{i}\\
....&.....&....&....\\
A^{j_{r},s_{r}}_{m_{1},l_{1}}&.....&A^{j_{r},s_{r}}_{m_{r},s_{r}}&c^{j_{r},s_{r}}_{i}\\
A^{j_{0},s_{0}}_{m_{1},l_{1}}&.....&A^{j_{0},s_{0}}_{m_{r}l_{r}}&c^{j_{0},s_{0}}_{i}
\end{vmatrix}=0.
\end{equation}

En développant cela par rapport à la dernière colonne, et en notant $D^{j_{k},s_{k}}_{j_{0},s_{0},i}$ le déterminant obtenu en supprimant la ligne $(j_{k},s_{k})$ et le dernières colonne, 
\begin{equation}
\sum_{k=1}^{r}(-1)^{r+1-k}c^{j_{k}s_{k}}_{i}D^{j_{k}s_{k}}_{j_{0},s_{0},i}+c^{j_{0},s_{0},i}_{i}D_{0}=0,
\end{equation}
c'est à dire $\forall (j,s)\in \overline{I},\;\; \forall i=1,...,n,$
\begin{equation}\label{am2}
c^{j,s}_{i}=-\frac{1}{D_{0}}\sum_{k=1}^{r}(-1)^{r+1-k}c^{j_{k}s_{k}}_{i}D^{j_{k}s_{k}}_{j,s,i}.
\end{equation}
Si $j\neq i,$ alors \ref{am2} peut être considéré comme une définition de $c^{j,s}_{i}$ et s'écrit aussi 

\begin{equation}\label{am2b}
\begin{split}
c^{j,s}_{i}&=-\frac{1}{D_{0}}\sum_{k=1}^{r}(-1)^{r+1-k}c^{j_{k},s_{k}}_{i}D^{j_{k},s_{k}}_{j,s,i}\\
&\qquad =-\frac{1}{D_{0}}\sum_{k\in I_{i}}(-1)^{r+1-k}c^{j_{k},s_{k}}_{i}D^{j_{k},s_{k}}_{j,s,i}+\frac{1}{D_{0}}\sum_{k\in \overline{I}_{i}}(-1)^{r+1-k}c^{i,s_{k}}_{j_{k}}D^{j_{k},s_{k}}_{j,s,i}\\
&\qquad =-\frac{1}{D_{0}}\sum_{k\in I_{i}}(-1)^{r+1-k}c^{j_{k},s_{k}}_{i}D^{j_{k},s_{k}}_{j,s,i}\\
&\qquad \qquad -\frac{1}{D_{0}^{2}}\sum_{k\in \overline{I}_{i},\;\;l=1,...,r.\;\;}(-1)^{k+l}D^{j_{k},s_{k}}_{j,s,i}c^{j_{l},s_{l}}_{j_{k}}D^{j_{l},s_{l}}_{i,s_{k},j_{k}},
\end{split}
\end{equation}

Mais, \ref{am2} impose aussi, d'après \ref{cond4}, la condition nécessaire $\forall (i,s)\in \overline{I},$
\begin{equation}\label{cond6}
c^{i,s}_{i}=\frac{\delta_{0}^{s}}{n}=-\frac{1}{D_{0}}\sum_{k=1}^{r}(-1)^{r+1-k}c^{j_{k}s_{k}}_{i}D^{j_{k}s_{k}}_{i,s,i}.
\end{equation}

Ainsi, tous les $c^{j,s}_{i}$ s'expriment en fonctions des $c^{j_{k},s_{k}}_{j_{l}}$ avec $k,l=1,...,,r$. Nous allons préciser cela en envisageant quatre cas. \\

\underline{Premier cas}: $(j,s)\in I$ et $(i,s)\in I.$\\
Alors $(j,s)=(j_{k},s_{k})\in I\;\; et\;\; (i,s)=(j_{l},s_{l})\in I$ donc $k,l\leqslant r, $ donc $c^{j,s}_{i}$ est connu. Il faut toutefois garder à l'esprit que $c_{j}^{i,s}=-c^{j,s}_{i} \;\; si\;\;i\neq j,,$ et aussi $c^{i,s}_{i}=\frac{\delta_{0}^{s}}{nVol}.$

\underline{Second cas}: $(j,s)\in \overline{I}\;\; et \;\;(i,s)\in I.$\\
Alors, $(j,s)=(j_{k_{1}},s_{k_{1}})$ et $(i,s)=(j_{k_{2}},s_{k_{2}})$ avec $k_{2}\leqslant r<k_{1},$
donc
\begin{equation}
c^{j,s}_{i}=-c^{i,s}_{j}=-\frac{1}{D_{0}}\sum_{k=1}^{r}(-1)^{r+1-k}c^{j_{k},s_{k}}_{j_{k_{2}}}D^{j_{k},s_{k}}_{j,s,i}
\end{equation}

\underline{Troisième cas} : $(j,s)\in I\;\; et\;\; (i,s)\in\overline{I}.$\\
Alors $(j,s)=(j_{k_{1}},s_{k_{1}})\;\; et\;\; (i,s)=(j_{k_{2}},s_{k_{2}}),$ donc
\begin{equation}
c^{j,s}_{i}=-c^{i,s}_{j}=\frac{1}{D_{0}}\sum_{k=1}^{r}(-1)^{r+1-k}c^{j_{k},s_{k}}_{j_{k_{1}}}D^{j_{k},s_{k}}_{i,s,j}
\end{equation}
D'après ces second et troisième cas, $c^{j,s}_{i}$ est connu, en fonction des $c^{j_{k},s_{k}}_{j_{l}}$ avec $k,l=1,...,,r,$ si $(j,s)\;\; ou\;\; (i,s)\in \overline{I}.$\\
Nous remarquons que ce deuxième et ce troisième cas sont les mêmes : on passe de l'un à l'autre en intervertissant $i$ et $j,$ et en changeant de signe.\\
\underline{Quatrième cas} : $(j,s)\;\; et\;\;(i,s)\in \overline{I},\;\; i\neq j,.$\\
Alors, $(j,s)=(j_{k_{1}},s_{k_{1}})\;\;et\;\;(i,s)=(j_{k_{2}},s_{k_{2}})\;\; avec\;\;k_{1}\;\; et \;\; k_{2}>r,$ donc 

\begin{equation}
\begin{split}
c^{j,s}_{i}& =-\frac{1}{D_{0}}\sum_{k\in I_{j_{k_{2}}}}(-1)^{r+1-k}c^{j_{k},s_{k}}_{j_{k_{2}}}D^{j_{k},s_{k}}_{j,s,i}\\
&\qquad \qquad -\frac{1}{D_{0}^{2}}\sum_{k\in \overline{I}_{j_{k_{2}}},\;\;l=1,...,r.\;\;}(-1)^{k+l}D^{j_{k},s_{k}}_{j,s,i}c^{j_{l},s_{l}}_{j_{k}}D^{j_{l},s_{l}}_{i,s_{k},j_{k}},
\end{split}
\end{equation}
mais aussi
\begin{equation}
\begin{split}
c^{j,s}_{i}&=-c^{i,s}_{j}=\frac{1}{D_{0}}\sum_{k\in I_{j_{k_{1}}}}(-1)^{r+1-k}c^{j_{k},s_{k}}_{j_{k_{1}}}D^{j_{k},s_{k}}_{i,s,j}\\
&\qquad \qquad +\frac{1}{D_{0}^{2}}\sum_{k\in \overline{I}_{j_{k_{1}}},\;\;l=1,...,r.\;\;}(-1)^{k+l}D^{j_{k},s_{k}}_{i,s,j}c^{j_{l},s_{l}}_{j_{k}}D^{j_{l},s_{l}}_{j,s_{k},j_{k}}.
\end{split}
\end{equation}
Ainsi, $c^{j,s}_{i}$ est écrit de deux façons qui doivent être égales, et donc, si $(j,s)=(j_{k_{1},s_{k_{1}}})=\in\overline{I},\;\; (i,s)=(j_{k_{2}},s_{k_{2}})\in \overline{I},\;\; et\;\; i\neq j,$
\begin{equation}\label{Rel2}
\sum_{k=1}^{r}(-1)^{k}[c^{j_{k},s_{k}}_{j}D^{j_{k},s_{k}}_{i,s,j}+c^{j_{k},s_{k}}_{i}D^{j_{k},s_{k}}_{j,s,i}]=0,
\end{equation}
ou
\begin{equation}\label{Rel18}
\begin{split}
&\sum_{k\in I_{j_{k_{2}}}}(-1)^{r+1-k}c^{j_{k},s_{k}}_{j_{k_{2}}}D^{j_{k},s_{k}}_{j,s,i}+\sum_{k\in I_{j_{k_{1}}}}(-1)^{r+1-k}c^{j_{k},s_{k}}_{j_{k_{1}}}D^{j_{k},s_{k}}_{i,s,j}\\
&\quad +\frac{1}{D_{0}}\Big[\sum_{k\in \overline{I}_{j_{k_{2}}},\;\;l=1,...,r.\;\;}(-1)^{k+l}D^{j_{k},s_{k}}_{j,s,i}c^{j_{l},s_{l}}_{j_{k}}D^{j_{l},s_{l}}_{i,s_{k},j_{k}}\\
&\quad +\sum_{k\in \overline{I}_{j_{k_{1}}},\;\;l=1,...,r.\;\;}(-1)^{k+l}D^{j_{k},s_{k}}_{i,s,j}c^{j_{l},s_{l}}_{j_{k}}D^{j_{l},s_{l}}_{j,s_{k},j_{k}}\Big]=0.
\end{split}
\end{equation}
De plus, d'après \ref{cond6}, si $(j,s)\in \overline{I},$ alors
\begin{equation}
\begin{split}
c^{j,s}_{j}=\frac{\delta_{0}^{s}}{n}&=-\frac{1}{D_{0}}\sum_{l=1}^{r}(-1)^{r+1-l}c^{j_{l},s_{l}}_{j}D^{j_{l}s_{l}}_{j,s,j}\\
&=\frac{-1}{D_{0}}\sum_{l\in I_{j}}(-1)^{r+1-l}c^{j_{l},s_{l}}_{j}D^{j_{l},s_{l}}_{j,s,j}+\frac{1}{D_{0}}\sum_{l\in \overline{I_{j}}}(-1)^{r+1-l}c^{j,s_{l}}_{j_{l}}D^{j_{l},s_{l}}_{j,s,j}\\
&=\frac{-1}{D_{0}}\sum_{l\in I_{j}}(-1)^{r+1-l}c^{j_{l},s_{l}}_{j}D^{j_{l},s_{l}}_{j,s,j}-\frac{1}{D_{0}^{2}}\sum_{l\in \overline{I_{j}}, l'=1,...,r}(-1)^{l+l'}c^{j_{l'},s_{l'}}_{j_{j_{l}}}D^{j_{l},s_{l}}_{j,s,j}D^{j_{l'},s_{l'}}_{j,s_{l},j_{l}}\\
&=\frac{-1}{D_{0}}\sum_{l\in I_{j},j_{l}\neq j}(-1)^{r+1-l}c^{j_{l},s_{l}}_{j}D^{j_{l},s_{l}}_{j,s,j}-\frac{1}{D_{0}^{2}}\sum_{l\in \overline{I_{j}},j_{l'}\neq j_{l}}(-1)^{l+l'}c^{j_{l'},s_{l'}}_{j_{j_{l}}}D^{j_{l},s_{l}}_{j,s,j}D^{j_{l'},s_{l'}}_{j,s_{l},j_{l}}\\
&\qquad -\frac{1}{D_{0}}\sum_{l\in I_{j},j_{l}= j}(-1)^{r+1-l}\frac{\delta_{0}^{s_{l}}}{n}D^{j,s_{l}}_{j,s,j}-\frac{1}{D_{0}^{2}}\sum_{l\in \overline{I_{j}},j_{l'}= j_{l}}(-1)^{l+l'}\frac{\delta_{0}^{s_{l'}}}{D_{0}}D^{j_{l},s_{l}}_{j,s,j}D^{j_{l},s_{l'}}_{j,s_{l},j_{l}}
\end{split}
\end{equation}
et donc
\begin{equation}\label{cond16}
\begin{split}
&\sum_{l\in I_{j},j_{l}\neq j}(-1)^{r+1-l}c^{j_{l},s_{l}}_{j}D^{j_{l},s_{l}}_{j,s,j}+\frac{1}{D_{0}}\sum_{l\in \overline{I_{j}},j_{l'}\neq j_{l}}(-1)^{l+l'}c^{j_{l'},s_{l'}}_{j_{j_{l}}}D^{j_{l},s_{l}}_{j,s,j}D^{j_{l'},s_{l'}}_{j,s_{l},j_{l}}\\
&\qquad =-\frac{\delta_{0}^{s}}{n}D_{0}-\frac{1}{n}\sum_{l\in I_{j},j_{l}= j}(-1)^{r+1-l}D^{j,0}_{j,s,j}-\frac{1}{D_{0}}\sum_{l\in \overline{I_{j}},j_{l'}= j_{l}}(-1)^{l+l'}D^{j_{l},s_{l}}_{j,s,j}D^{j_{l},0}_{j,s_{l},j_{l}}
\end{split}
\end{equation}

Comme déjà dit, dans tous les cas, les coefficients $c^{j,s}_{i}$ s'écrivent en fonction des $c^{j_{k},s_{k}}_{j_{l}}\;\; avec\;\; k\;\; et\;\; l\leqslant r.$ Dans le premier terme de \ref{Rel18} (et, de façon analogue, dans le deuxième terme), on pourrait penser que l'indice $j_{k_{2}}$ ne peut s'écrire sous la forme $j_{l}$ avec $l\leqq r.$ Mais en fait, puisque dans ce premier terme $k\in I_{j_{k_{2}}},$ c'est que $(j_{k_{2}},s_{k})\in I,$ et donc il existe $l\leqq r$ tel que $(j_{k_{2}},s_{k})=(j_{l},s_{l}),$ et, par conséquent, $j_{k_{2}}=j_{l}$ avec $l\leqslant r.$ \\
Puisque tous les $c^{j,s}_{i}$ s'écrivent en fonction des $c^{j_{k},s_{k}}_{j_{l}}$ avec $k,l\leqslant r,$ nous devons choisir les $c^{j_{k},s_{k}}_{j_{l}}, \;\; k,l\leqslant r\;\; et\;\; j_{k}\neq j_{l,}$ tels que les $c^{j,s}_{i}$ ainsi définis vérifient \ref{cond4} et \ref{am2b}. \\
Pour ce faire, si $c^{j,s}_{i}$ relève du quatrième cas, il nous faut imposer \ref{Rel18} et \ref{cond16}. \\
Si $c^{j,s}_{i}$ relève du troisième ou du deuxième cas (qui sont les mêmes), nous avons forcément $i\neq j$ car, parmi $(i,s)$ et $(j,s),$ l'un est dans $I,$ l'autre dans $\overline{I},$ et alors $c^{j,s}_{i}=-c^{i,s}_{j}.$\\

Si $ k,l\leqslant r,\;\; tels \;\; que\;\; (j_{l},s_{k})\in \overline{I},$
\begin{equation}
	c^{j_{k},s_{k}}_{j_{l}}=-c^{j_{l},s_{k}}_{j_{k}}=\frac{1}{D_{0}}\sum_{l'=1}^{r}(-1)^{r+1-l'}c^{j_{l'},s_{l'}}_{j_{k}}D^{j_{l'},s_{l'}}_{j_{l},s_{k},j_{k}},
\end{equation}

ou encore
\begin{equation}\label{Rel20}
	D_{0}c^{j_{k},s_{k}}_{j_{l}}-\sum_{l'\;:\; j_{l'}\neq j_{k}}(-1)^{r+1-l'}c^{j_{l'},s_{l'}}_{j_{k}}D^{j_{l'},s_{l'}}_{j_{l},s_{k},j_{k}}=\sum_{l'\;:\;j_{l'}=j_{k},\;s_{l'=0}}^{r}\frac{(-1)^{r+1-l'}}{n}D^{j_{k},0}_{j_{l},s_{k},j_{k}}.
\end{equation}
Enfin, si $c^{j,s}_{i}$ relève du le premier cas, alors $c^{j,s}_{i}=c^{j_{k},s_{k}}_{j_{l}},\;\;k,l\leqslant r$  et, si $ j_{k}\neq j_{l},\;\; alors\;\; c^{j_{k},s_{k}}_{j_{l}}=-c^{j_{l},s_{k}}_{j_{l}}.$\\
Il nous faut donc rechercher les $c^{j_{k},s_{k}}_{j_{l}}\;\; avec\;\; k\;\; et\;\; l\leqslant r,\;\; j_{k}\neq j_{l},$ vérifiant les conditions \ref{Rel18}, \ref{cond16} et \ref{Rel20}, c'est à dire 

\begin{equation}\label{S}
\begin{cases}
\forall (i,s)\;et\; (j,s)\in \overline{I},\;\ i\neq j,\\
\sum_{l\in I_{j},\;j_{l}\neq j=j_{k}, k\leqslant r}(-1)^{r+1-l}c^{j_{l},s_{l}}_{j_{k}}D^{j_{l},s_{l}}_{i,s,j}+\sum_{l\in I_{j},\;j_{l}\neq i=j_{k},\; k\leqslant r}(-1)^{r+1-l}c^{j_{l},s_{l}}_{j_{k}}D^{j_{l},s_{l}}_{j,s,i}\\
\qquad +\frac{1}{D_{0}}\Big[\sum_{l\in \overline{I}_{j},\;\;l':j_{l'}\neq j_{l}}(-1)^{l+l'}c^{j_{l'},s_{l'}}_{j_{l}}D^{j_{l},s_{l}}_{i,s,j}D^{j_{l'},s_{l'}}_{j,s_{l},j_{l}}+\sum_{l\in \overline{I}_{i},\;\;l' :j_{l'}\neq j_{l}}(-1)^{l+l'}c^{j_{l'},s_{l'}}_{j_{l}}D^{j_{l},s_{l}}_{j,s,i}D^{j_{l'},s_{l'}}_{i,s_{l},j_{l}}\Big]\\
\quad =\frac{-1}{n}\Big\{\sum_{l\in I_{j},\;j_{l}= j,\;s_{l}=0}(-1)^{r+1-l}D^{j,0}_{i,s,j}+\sum_{l\in I_{i},\;j_{l}= i,\;s_{l}=0}(-1)^{r+1-l}D^{i,0}_{j,s,i}\\
\qquad +\frac{1}{D_{0}}\Big[\sum_{l\in \overline{I}_{j},\;\;l':j_{l'}= j_{l},\;s_{l'}=0}(-1)^{l+l'}D^{j_{l},s_{l}}_{i,s,j}D^{j_{l},0}_{j,s_{l},j_{l}}+\sum_{l\in \overline{I}_{i},\;\;l' :j_{l'}= j_{l},\;s_{l'}=0}(-1)^{l+l'}D^{j_{l},s_{l}}_{j,s,i}D^{j_{l},0}_{i,s_{l},j_{l}}\Big]\Big\}\\
\forall (j,s)\in \overline{I},\\
\sum_{l\in I_{j},\;j_{l}\neq j=j_{k},\; k\leqslant r}(-1)^{r+1-l}c^{j_{l},s_{l}}_{j_{k }}D^{j_{l},s_{l}}_{i,s,j}
 +\frac{1}{D_{0}}\sum_{l\in \overline{I}_{j},\;\;l':j_{l'}\neq j_{l}}(-1)^{l+l'}c^{j_{l'},s_{l'}}_{j_{l}}D^{j_{l},s_{l}}_{i,s,j}D^{j_{l'},s_{l'}}_{j,s_{l},j_{l}}\\
\quad =\frac{-1}{nD_{0}}\Big[D_{0}\delta_{0}^{s}+\sum_{l\in I_{j},\;j_{l}= j,\;s_{l}=0}(-1)^{r+1-l}D^{j,0}_{i,s,j} +\frac{1}{D_{0}}\sum_{l\in \overline{I}_{j},\;\;l':j_{l'}= j_{l},\;s_{l'}=0}(-1)^{l+l'}D^{j_{l},s_{l}}_{i,s,j}D^{j_{l},0}_{j,s_{l},j_{l}}\Big]\\
\forall k,l\leqslant r,\;\; et \; ,\; (j_{l},s_{k})\in \overline{I}\\
D_{0}c^{j_{k},s_{k}}_{j_{l}}-\sum_{l'=1,\;\,j_{l'}\neq j_{k}}^{r}(-1)^{r+1-l'}c^{j_{l'},s_{l'}}_{j_{k}}D^{j_{l'},s_{l'}}_{j_{l},s_{k},j_{k}}=\frac{1}{n}\sum_{l'=1,\;\;j_{l'}=j_{k},\;\; s_{l'}=0}^{r}(-1)^{r+1-l'}D^{j_{k},0}_{j_{l},s_{k},j_{k}}\\
\forall k,l\leqslant r,\;\; et\;\; (j_{l},s_{k})\in I,\;\; \\
c^{j_{k},s_{k}}_{j_{l}}=-c^{j_{l},s_{k}}_{j_{k}}.
\end{cases}
\end{equation}
C'est donc ce système dont les inconnues sont les $c^{j_{k},s_{k}}_{j_{l}},\;\;k,l\leqslant r,\;\; j_{k}\neq j_{l}$ que nous voulons résoudre.\\
Nous avons une formule de Cauchy pour la classe des fonctions vérifiant $d''f=0$ si et seulement si ce système admet des solutions.\\
Le nombre $N$ des équations de ce système et son rang $R$ dépendent des équations principales. Donc, nous n'irons pas plus loin dans l'étude du cas général. \\
Rappelons simplement que si $R=N,$ alors nous pouvons fixer arbitrairement les valeurs des inconnues non principales (sil y en a) et résoudre le système. Si $R<N,$ alors nous pouvons fixer les valeurs des inconnues non principales (s'il y en a), mais, pour que le système admette des solutions, nous avons $N-R$ conditions données par la nullité des déterminants caractéristiques. \\
Pour la mise en œuvre de cette méthode, on se reportera aux exemples ci-dessous, en particulier au cas de dimension 4.\\
\begin{rem}
	Les résultats précédents fournissent une formule de Cauchy pour la classe des fonctions vérifiant \ref{C2a} et continues jusqu'au bord. Toutefois, ils ne disent rien, a priori, de la taille, la richesse (le "nombre" de fonctions) de la classe considérée. \\
	Les fonctions de cette classe sont $\R-analytiques$ puisque $K(y,x)$ l'est. On peut donc utiliser le théorème de Cauchy-Kowalevskia (si l'un des coefficients de \ref{C2a} est inversible), ou les résultats de \cite{BCGGG} (torsion nulle et tableau involutif (voir p. 140)) pour s'assurer de la richesse de cette classe. \\
	Dans les exemples qui suivent, aussi bien dans le cas commutatif que non commutatif, le lecteur pourra facilement  s'assurer de la taille de la classe considérée. 
\end{rem}

\vspace{10mm}

\subsection{Le cas particulier d'une algèbre commutative.} 
Lorsque l'algèbre $A$ est commutative, nous supposerons que le système \ref{syst4} est de rang $q$, c'est-à-dire qu'il existe un déterminant $D_{0}$ (dont les entrées sont des éléments de l'algèbre commutative $A$) inversible et d'ordre maximal $q$ extrait de la matrice du système, et que ce déterminant est formé par les $q$ premières équations qui correspondent à $j=1,...,q.$ Comme le système comporte $n$ équations, nous obtiendrons $n-q$ conditions de compatibilité, pour que le système admette des solutions, en écrivant la nullité des déterminants caractéristiques obtenus en rajoutant à $D_{0}$ la ligne $q+k, \quad k=1,...,n-q,$ et les seconds membres correspondants. Si nous appelons $D^{k}_{m}$ le déterminant obtenu à partir de $D_{0}$ en supprimant la ligne $m$ et en rajoutant la ligne $q+k,$ les $n-q$ conditions de compatibilité, obtenues comme précisé ci-dessus en développant un déterminant caractéristique, s'écrivent :
\begin{equation} 
\sum^{q}_{m=1}(-1)^{q-1+m}D^{k}_{m}c^{m}_{i}\;+\;D_{0}c^{q+k}_{i}\;=\;0 \quad pour \quad tout \quad i=1,...,n, \quad et \quad tout \quad k=1,...,n-q. 
\end{equation} 
Ainsi 
\begin{equation}\label{cond5}
c^{q+k}_{i}\;=\;D^{-1}_{0}\sum^{q}_{m=1}(-1)^{q+m}D^{k}_{m}c^{m}_{i}
\end{equation} 
exprime $c^{q+k}_{i}$ en fonction des $c^{m}_{i},\quad m\leq q.$\\
Les coefficients $c^{j}_{i}$ étant choisis vérifiant \ref{cond4}, pour $1\leq i,j\leq q,$ nous pouvons, dans un premier temps, grâce à \ref{cond5}, calculer les coefficients $c^{q+k}_{i},\quad i=1,...,q, \quad k=1,...,n-q,$ puis, dans un second temps, exprimer les coefficients $c^{i}_{q+k}\;=\;-c^{q+k}_{i}\;=\;-D^{-1}_{0}\sum^{q}_{m=1}(-1)^{q+m}D^{k}_{m}c^{m}_{i},$ enfin, dans un troisième temps, obtenir, pour $1\leq k,l\leq n-q,$ les coefficients 
\begin{equation}
\begin{split}
c^{q+l}_{q+k}\;&=\;D^{-1}_{0}\sum^{q}_{m=1}(-1)^{q+m}D^{l}_{m}c^{m}_{q+k}\\
&=\;-D^{-2}_{0}\sum^{q}_{m=1}(-1)^{q+m}D^{l}_{m}\sum^{q}_{m'=1}(-1)^{q+m'}D^{k}_{m'}c^{m'}_{m}\\
&=\;-D^{-2}_{0}\sum^{q}_{m,m'=1}(-1)^{m+m'}D^{l}_{m}D^{k}_{m'}c^{m'}_{m}\\
&=\;-D^{-2}_{0}\sum^{q}_{m=1}D^{l}_{m}D^{k}_{m}\frac{e_{0}}{n}\;-\;D^{-2}_{0}\sum_{1\leq m<m'\leq q}(-1)^{m+m'}\big(D^{l}_{m}D^{k}_{m'}c^{m'}_{m}\,+\,D^{l}_{m'}D^{k}_{m}c^{m}_{m'}\big)\\
&=\;-D^{-2}_{0}\sum^{q}_{m=1}\frac{D^{l}_{m}D^{k}_{m}}{n}\;-\;D^{-2}_{0}\sum_{1\leq m<m'\leq q}(-1)^{m+m'}\big(D^{l}_{m}D^{k}_{m'}\,-\,D^{l}_{m'}D^{k}_{m'}\big)c^{m'}_{m}.
\end{split}
\end{equation}
Il reste maintenant à s'assurer que ces coefficients $c^{j}_{i}$ satisfont les conditions \ref{cond4}. Dans le cas où $1\leq i,j\leq q,$ cela relève de notre choix du départ. $c^{q+k}_{i}\;+\;c^{i}_{q+k}\;=\;0\quad i\leq q,\quad k\leq n-q$ est vérifié, de par la définition de $c^{i}_{q+k}.$ Il reste à vérifier que
\begin{equation}
\begin{split}
c^{q+k}_{q+k}\;=\;\frac{e_{0}}{n}\;&=\;D^{-1}_{0}\sum^{q}_{m=1}(-1)^{q+m}D^{k}_{m}c^{m}_{q+k}\\
&=\;-D^{-1}_{0}\sum^{q}_{m=1}(-1)^{q+m}D^{k}_{m}c^{q+k}_{m}\\
&=-D^{-2}_{0}\sum^{q}_{m=1}(-1)^{q+m}D^{k}_{m}\sum^{q}_{m'=1}(-1)^{q+m'}D^{k}_{m'}c^{m'}_{m}\\
&=\;-\frac{D^{-2}_{0}}{n}\sum^{q}_{m=1}(D^{k}_{m})^{2}\;-\;D^{-2}_{0}\sum_{m<m'}(-1)^{m+m'}D^{k}_{m}D^{k}_{m'}(c^{m'}_{m}\;+\;c^{m}_{m'})\\
&=\;-\frac{D^{-2}_{0}}{n}\sum^{q}_{m=1}(D^{k}_{m})^{2}
\end{split}
\end{equation}
et que, pour $1\leq k\neq l\leq n-q,$ 
\begin{equation}
\begin{split}
c^{q+l}_{q+k}\;+\;c^{q+k}_{q+l}\;&=\;-\frac{D^{-2}_{0}}{n}\sum_{m}D^{l}_{m}D^{k}_{m}\;-\;D^{-2}_{0}\sum_{m<m'}(-1)^{m+m'}(D^{l}_{m}D^{k}_{m'}\,-\,D^{l}_{m'}D^{k}_{m})C^{m'}_{m}\\
&+\;-\frac{D^{-2}_{0}}{n}\sum_{m}D^{k}_{m}D^{l}_{m}\;-\;D^{-2}_{0}\sum_{m<m'}(-1)^{m+m'}(D^{k}_{m}D^{l}_{m'}\,-\,D^{k}_{m'}D^{l}_{m})C^{m'}_{m}\\
&=\;-\frac{2D^{-2}_{0}}{n}\sum_{m}D^{l}_{m}D^{k}_{m}\\
&=\;0.
\end{split}
\end{equation}
L'avant dernière équation implique
$$D^{2}_{0}\;+\;\sum^{q}_{m=1}(D^{k}_{m})^{2}\;=\;0,$$
quand la dernière entraîne 
$$\sum^{q}_{m=1}D^{l}_{m}D^{k}_{m}\;=\;0.$$
Finalement, nous obtenons, comme conditions d'existence d'une représentation intégrale pour les fonctions $f$ satisfaisant nos conditions de Cauchy-Riemann \ref{C2a} : 
\begin{theorem}\label{cond20} Dans une algèbre commutative, les fonctions $f$ v\'erifiant les \' equations de Cauchy-Riemann admettent une repr\'esentation int\'egrale si les conditions ci-dessous sont v\'erifi\'ees
\begin{equation}
Conditions \quad(A):\qquad \sum^{q}_{m=1}D^{l}_{m}D^{k}_{m}\;+\;\delta^{l}_{k}D^{2}_{0}\;=\;0\quad pour \quad tout \quad 1\leq k,l\leq n-q.
\end{equation}
\end{theorem}

Pour terminer ce cas commutatif, il nous reste à remarquer que, les conditions ci-dessus étant satisfaites, les solutions obtenues pour le système \ref{syst1}, sont données par: 
$$\varphi_{m}\;=\;D^{-1}_{0}\widetilde{D^{m}}$$ 
où $\widetilde{D^{m}}$ est obtenu, à partir de $D_{0}$, en remplaçant la colonne $m$ par $(\Psi_{1},...,\Psi_{q}).$ Par conséquent, $\varphi_{m}(x,y)=\sum^{n}_{i=1}b^{i}_{m}(y_{i}-x_{i})$ avec
\begin{equation}\label{cond8}
b^{i}_{m}\;=\;\sum^{q}_{l=1}(-1)^{m+l}c^{l}_{i}D^{-1}_{0}\widehat{D^{l}_{m}}\;=\;D^{-1}_{0}\widetilde{D^{i}_{m}} 
\end{equation}
où $\widehat{D^{l}_{m}}$ est obtenu, à partir de $D_{0}$, en retirant la ligne $l$ et la colonne $m,$ et $\widetilde{D^{i}_{m}}$ en remplaçant la colonne $m$ de $D_{0}$ par la colonne $(c^{1}_{i},...,c^{q}_{i}).$\\
\begin{rem}
Le cas $q=n-1.$
\end{rem}
Les conditions (A) sont au nombre de $\frac{(n-q)(n-q+1)}{2}.$ Si $q=n-1,$ alors il n'y a qu'une seule condition qui s'écrit 
$$D^{2}_{0}\;+\;\sum^{n-1}_{m=1}(D^{1}_{m})^{2}\;=\;0.$$

\begin{ex}
Le cas d'une fonction A-différentiable.
\end{ex}
Supposons que $f$ est une fonction continue sur la ferméture $\overline{D}$ d'un ouvert $D$ de l'algèbre commutative $A,$ et $A-$différentiable, c'est-à-dire que pour tout $x\in D,$ la différentielle $df(x)$ est $A-$linéaire. Ceci est équivalent à dire que, pour tout $m=1,...,p,$ 
\begin{equation}
\begin{split}
\frac{\partial f(x)}{\partial x_{m}}\;&=\;df(x)(e_{m})\\
&=\;e_{m}df(x)(e_{0})\\
&=\;e_{m}\frac{\partial f(x)}{\partial x_{0}}.
\end{split}
\end{equation}
Ainsi, dire que $f$ est $A-$différentiable est équivalent à imposer $p$ conditions de Cauchy à $f:$ 
\begin{equation}\label{defdif}
\frac{\partial f}{\partial x_{m}}\;-\;e_{m}\frac{\partial f}{\partial x_{0}}\;=\;0,\quad m=1,...,p.
\end{equation}
Nous sommes dans le cas de la remarque précédente avec $a^{0}_{m}\,=\,-e_{m}, \quad a^{m}_{m}\,=\,e_{0},$ et tous les autres coefficients $a^{j}_{m}$ nuls. Adaptons la démonstration générale à ce cas particulier.\\
Le système \ref{syst1} s'écrit 
\begin{equation}
\begin{cases}
(y_{0}-x_{0})\frac{e_{0}}{n}+\sum_{i=1}^{p}(y_{i}-x_{i)}c^{0}_{i}=-Vol\sum_{m=1}^{p}\sum_{i=0}^{p}e_{m}b^{i}_{m}(y_{i}-x_{i})\\
(y_{j}-x_{j})\frac{e_{0}}{n}+\sum_{i\neq j}(y_{i}-x_{i)}c^{j}_{i}=Vol\sum_{i=0}^{p}b^{i}_{j}(y_{i}-x_{i}), 
\end{cases}
\end{equation}

qui se traduit par 
\begin{equation}
\begin{cases}
\frac{e_{0}}{n}=-Vol\sum_{m=1}^{p}e_{m}b^{0}_{m}\\
c^{0}_{i}=-Vol\sum_{m=1}^{p}e_{m}b^{i}_{m}\\
\frac{e_{0}}{n}=Vol\;b^{j}_{j},\;\;\forall j\geqslant 1\\
c^{j}_{i}=Vol\; b^{i}_{j},\;\;\forall j\geqslant 1,\;\; i\neq j.
\end{cases}
\end{equation}
Reportant $b^{j}_{j}$ et $b^{i}_{j}$ extraits des deux dernières équations dans les deux premières, nous obtenons :
\begin{equation}
\begin{cases}
\frac{e_{0}}{n}=-\sum_{m=1}^{p}e_{m}c^{m}_{0}\\
c^{0}_{i}-\sum_{m=1}^{p}e_{m}c^{m}_{i}.
\end{cases}
\end{equation}

Nous reportons $c^{0}_{i}=-c^{i}_{0},$ donné par la dernière équation, dans la première pour obtenir, en utilisant la commutativité de l'algèbre, 
\begin{equation}
\begin{split}
\frac{e_{0}}{n}&=-\sum_{m=1}^{p}e_{m}\sum_{m'=1}^{p}e_{m'}c^{m'}_{m}\\
&=-\sum_{m<m'}e_{m}e_{m'}c^{m'}_{m}-\sum_{m=1}^{p}e_{m}^{2}c^{m}_{m}-\sum_{m>m'}e_{m}e_{m'}c^{m'}_{m}\\
&=-\sum_{m<m'}e_{m}e_{m'}c^{m'}_{m}-\sum_{m=1}^{p}\frac{e_{m}^{2}}{n}-\sum_{m<m'}e_{m'}e_{m}c^{m}_{m'}\\
&=-\sum_{m<m'}e_{m}e_{m'}(c^{m'}_{m}+c^{m}_{m'})-\sum_{m=1}^{p}\frac{e_{m}^{2}}{n}\\
&=-\sum_{m=1}^{p}\frac{e_{m}^{2}}{n}.
\end{split}
\end{equation}

Ainsi, nous obtenons le corollaire, 
\begin{cor} Une fonction A-diff\'erentiable v\'erifie une formule de Cauchy si l'alg\` ebre A satisfait la condition
\begin{equation}\label{Acom}
\sum_{m=0}^{p}e_{m}^{2}=0.
\end{equation}
\end{cor}
Nous reconnaissons ici la condition $(A_{0})$ de \cite{BC1} et \cite{BC2}. En dimension 2, c'est-à-dire lorsque $A=\C,$ nous reconnaissons tout simplement la condition $1+i^{2}=0.$
\\
\\

\textbf{Le cas de la dimension 2.}\\
\\
Une algèbre $A$ unitaire de dimension 2 est toujours commutative. Si $(e_{0},e_{1})$ est une base, avec $e_{0}$ qui est l'élément unité, alors, il suffit, pour définir le produit dans l'algèbre $A,$ de donner 
$$(e_{1})^{2}=ae_{0}+be_{1}.$$
\begin{rem}\label{rq}
\end{rem}
Il est alors très facile de vérifier que les assertions suivantes sont equivalentes : 
\begin{equation}
\begin{split}
&\qquad \bullet A\;\;est\;\;integre,\\
&\qquad \bullet A\;\;est\;\;un\;\;corps,\\
&\qquad \bullet il\;\; existe\;\;e'_{1}\in A\;\;tel\;\;que\;\;(e'_{1})^{2}=-e_{0}\\
&\qquad \bullet b^{2}+4a<0.
\end{split}
\end{equation}
Ainsi, dans la base $(e_{0},e'_{1}),$ nous avons alors une version de $\C.$\\
\begin{theorem}
	Une algèbre $A$ définie par une seule condition de Cauchy \ref{C2a} (q=1) et admettant une formule de Cauchy est un corps, et donc est isomorphe à $\C.$
\end{theorem}
En effet, \ref{cond12} s'écrit 
\begin{equation}\begin{cases}
	a_{1}^{1}b_{1}^{1}=\frac{e_{0}}{2Vol}\\
	a_{1}^{2}b_{1}^{2}=\frac{e_{0}}{2Vol}\\
	a_{1}^{1}b_{1}^{2}+a_{1}^{2}b_{1}^{1}=0.\\
\end{cases}
\end{equation}
Les deux premières équations traduisent le fait que $b_{1}^{1}=\frac{(a_{1}^{1})^{-1}}{2Vol}$ et $b_{1}^{2}=\frac{((a_{1}^{2})^{-1}}{2Vol}.$ En reportant dans la troisième équation, $a_{1}^{1}(a_{1}^{2})^{-1}+a_{1}^{2}
(a_{1}^{1})^{-1}=0,$ c'est à dire $a_{1}^{1}(a_{1}^{2})^{-1}=-a_{1}^{2}
(a_{1}^{1})^{-1},$ ou encore 
\begin{equation}
	-e_{0}=\big(\frac{a_{1}^{2}}{a_{1}^{1}}\big)^{2}=\big(\frac{a_{1}^{1}}{a_{1}^{2}}\big)^{2},
\end{equation}
et, d'après la remarque ci-dessus, $A$ est une version de $\C.$ \\
\\

Nous allons maintenant revenir au cas général, c'est à dire au cas d'une algèbre non commutative de dimension quelconque.\\

Dans les exemples ci-dessous, les conditions $(A)$ nécessaires et suffisantes pour l'obtention d'une formule de Cauchy pour la classe des fonctions satisfaisant aux conditions de Cauchy  \ref{C2a} seront traduites par l'appartenance d'un vecteur à un sous-module du A-module $A^{C^{2}_{n-q+1}}.$ \\

\textbf{Un cas particulier} avec $q=p\;\;et\;\;n=p+1.$\\

Nous allons utiliser la condition matricielle \ref{condmat}.\\
$\widetilde{A}$ est alors une matrice à $p+1$ lignes et $p$ colonnes.\\
Nous supposons que la première ligne de $\widetilde{A}$ est:\\
$\begin{pmatrix}\alpha_{1}&\alpha_{2}&\dots&\alpha_{p}\end{pmatrix}$\\
et que les $p$ dernières lignes de $\widetilde{A}$ forment une matrice carrée d'ordre $p$ diagonale comportant, sur la diagonale principale $\beta_{1},\beta_{2},...,\beta_{p},$ éléments inversibles de $A.$ Nous cherchons une matrice à $p+1$ lignes et $p$ colonnes
$$B\;=\;\big(b^{j}_{m}\big)_{j=0,...,p;\; m=1,...,p}$$
vérifiant la condition \ref{condmat}.Nous calculons 
\begin{equation}
\widetilde{A}\,^{t}B\;=\;
\begin{pmatrix}
\sum^{p}_{m=1}\alpha_{m}b^{0}_{m}&\sum^{p}_{m=1}\alpha_{m}b^{1}_{m}&\dots&\sum^{p}_{m=1}\alpha_{m}b^{p}_{m}\\
\beta_{1}b^{0}_{1}&\beta_{1}b^{1}_{1}&\dots&\beta_{1}b^{p}_{1}\\
\beta_{2}b^{0}_{2}&\beta_{2}b^{1}_{2}&\dots&\beta_{2}b^{p}_{2}\\
\hdotsfor[2]{4}\\
\beta_{p}b^{0}_{p}&\beta_{p}b^{1}_{p}&\dots&\beta_{p}b^{p}_{p}
\end{pmatrix}
\end{equation}
et écrivons les conditions \ref{condmat}\\
$$\begin{cases}
\sum^{p}_{m=1}\alpha_{m}b^{0}_{m}\;=\;\beta_{k}b^{k}_{k}\;=\;\frac{e_{0}}{(p+1)Vol(B(0,1))}\quad pour\quad k=1,...,p,\\
\beta_{k}b^{0}_{k}\;=\;-\sum^{p}_{m=1}\alpha_{m}b^{k}_{m}\quad pour\quad k=1,...,p,\\
\beta_{k}b^{l}_{k}\;+\;\beta_{l}b^{k}_{l}\;=\;0,\quad pour k\neq l,
\end{cases}$$
qui se traduisent par\\
$$\begin{cases}
\sum^{p}_{m=1}\alpha_{m}b^{0}_{m}\;=\;\frac{e_{0}}{(p+1)Vol(B(0,1))},\\
b^{k}_{k}\;=\;\frac{\beta^{-1}_{k}}{(p+1)Vol(B(0,1))}\quad pour \quad k=1,...,p,\\
b^{0}_{k}\;=\;-\beta^{-1}_{k}\sum^{p}_{m=1}\alpha_{m}b^{k}_{m}\quad pour \quad  k=1,...,p,\\
b^{l}_{k}\;=\;-\beta^{-1}_{k}\beta_{l}b^{k}_{l}.
\end{cases}$$
Nous en déduisons la condition de possibilité:
\begin{equation}
\begin{split}
\frac{e_{0}}{(p+1)Vol(B(0,1))}\;&=\;\sum^{p}_{m=1}\alpha_{m}b^{0}_{m}\\
&=\;-\sum^{p}_{m,k=1}\alpha_{m}\beta^{-1}_{m}\alpha_{k}b^{m}_{k}\\
&=\;\frac{1}{(p+1)Vol(B(0,1))}\sum^{p}_{m=1}(\alpha_{m}\beta^{-1}_{m})^{2}\,-\;\sum_{k<m}(\alpha_{k}\beta^{-1}_{k}\alpha_{m}b^{k}_{m}-\alpha_{m}\beta^{-1}_{m}\alpha_{k}\beta^{-1}_{k}\beta_{m}b^{k}_{m})\\
&=\;\frac{-1}{(p+1)Vol(B(0,1))}\sum^{p}_{m=1}(\alpha_{m}\beta^{-1}_{m})^{2}\;-\;\sum_{k<m}(\alpha_{k}\beta^{-1}_{k}\alpha_{m}\beta^{-1}_{m}-\alpha_{m}\beta^{-1}_{m}\alpha_{k}\beta^{-1}_{k})\beta_{m}b^{k}_{m}.
\end{split}
\end{equation}
Comme les éléments $b^{k}_{m}\in A,\;k<m,$ peuvent être choisis arbitrairement, cette condition conduit au th\'eor\` eme,

\begin{theorem} Dans ce cas particulier, les fonctions de Cauchy-Riemann de A v\'erifient une formule de Cauchy si les conditions suivantes sont satisfaites :
\begin{equation}\label{Anc}
e_{0}\;+\;\sum^{p}_{m=1}(\alpha_{m}\beta^{-1}_{m})^{2}\;\in I(\alpha_{m}\beta^{-1}_{m}\alpha_{k}\beta^{-1}_{k}\;-\;\alpha_{k}\beta^{-1}_{k}\alpha_{m}\beta^{-1}_{m},\;\;1\leq k<m\leq p)
\end{equation}
où $I(\alpha_{m}\beta^{-1}_{m}\alpha_{k}\beta^{-1}_{k}\;-\;\alpha_{k}\beta^{-1}_{k}\alpha_{m}\beta^{-1}_{m},\;\;k<m)$ désigne l'idéal de $A$ engendré \\
par les éléments $\alpha_{m}\beta^{-1}_{m}\alpha_{k}\beta^{-1}_{k}\;-\;\alpha_{k}\beta^{-1}_{k}\alpha_{m}\beta^{-1}_{m}.$
\end{theorem}

\textbf{Un cas encore plus particulier :} le cas d'une fonction $A$-différentiable.\\

Une fonction $f$ définie sur un ouvert $D$ de l'algèbre $A,$ et à valeurs dans $A,$ est dite $A-$différentiable, à gauche, si sa différentielle vérifie les conditions (\ref{defdif}). Il s'agit donc d'un cas encore plus particulier du cas particulier précédent, cas plus particulier obtenu en posant $\alpha_{m}=-e_{m}$ et $\beta_{m}=e_{0},\;\;m=1,...,p.$ La condition $(A),$ qui doit être satisfaite pour que nous puissions obtenir une représentation intégrale pour $f,$ se déduit de (\ref{Anc}) et nous obtenons le corollaire
\begin{cor}Si $f$ est une fonction A-diff\'erentiable, elle v\'erifie une formule de Cauchy si les conditions ci-dessous sont v\'erifi\'ees
\begin{equation}
\sum^{p}_{m=0}e^{2}_{m}\in I\big(e_{m}e_{k}-e_{k}e_{m},\;\;1\leq k<m\leq p\big),
\end{equation}
où, encore une fois, $I\big(e_{m}e_{k}-e_{k}e_{m},\;\;1\leq k<m\leq p\big)$ désigne l'idéal à gauche de $A$ engendré par les éléments $e_{m}e_{k}-e_{k}e_{m},\;\;1\leq k<m\leq p$ de $A.$\end{cor}

Dans le cas commutatif, cet idéal est nul, et l'on retrouve la condition $(A)$ de (\ref{Acom}).\\
Bien sûr, dans le cas de la $A$-différentiabilité à droite, on obtient une condition analogue. \\

\textbf{Un autre cas particulier avec } $q=p-1$ et $n=p+1.$\\

Ici, nous prendrons pour $\widetilde{A}$ une matrice à $p+1$ lignes et $p-1$ colonnes (et donc $B$ aussi). Nous noterons 
\begin{equation}
\widetilde{A}\;=\;
\begin{pmatrix}
\alpha_{2}&\alpha_{3}&\alpha_{4}&\dots&\alpha_{p}\\
\beta_{2}&\beta_{3}&\beta_{4}&\dots&\beta_{p}\\
\gamma_{2}&0&0&\dots&0\\
0&\gamma_{3}&0&\dots&0\\
0&0&\gamma_{4}&\dots&0\\
\hdotsfor[2]{5}\\
0&0&0&\dots&\gamma_{p}
\end{pmatrix}
\end{equation}
et
\begin{equation}
B\;=\;
\begin{pmatrix}
b^{0}_{2}&b^{0}_{3}&\dots&b^{0}_{p}\\
b^{1}_{2}&b^{1}_{3}&\dots&b^{1}_{p}\\
\hdotsfor[2]{4}\\
b^{p}_{2}&b^{p}_{3}&\dots&b^{p}_{p}
\end{pmatrix}
\end{equation}
donc
\begin{equation}
\widetilde{A}\,^{t}B\;=\;
\begin{pmatrix}
\sum^{p}_{m=2}\alpha_{m}b^{0}_{m}&\sum^{p}_{m=2}\alpha_{m}b^{1}_{m}&\sum^{p}_{m=2}\alpha_{m}b^{2}_{m}&\dots&\sum^{p}_{m=2}\alpha_{m}b^{p}_{m}\\
\sum^{p}_{m=2}\beta_{m}b^{0}_{m}&\sum^{p}_{m=2}\beta_{m}b^{1}_{m}&\sum^{p}_{m=2}\beta_{m}b^{2}_{m}&\dots&\sum^{p}_{m=2}\beta_{m}b^{p}_{m}\\
\gamma_{2}b^{0}_{2}&\gamma_{2}b^{1}_{2}&\gamma_{2}b^{2}_{2}&\dots&\gamma_{2}b^{p}_{2}\\
\gamma_{3}b^{0}_{3}&\gamma_{3}b^{1}_{3}&\gamma_{3}b^{2}_{3}&\dots&\gamma_{3}b^{p}_{3}\\
\hdotsfor[2]{5}\\
\gamma_{p}b^{0}_{p}&\gamma_{p}b^{1}_{p}&\gamma_{p}b^{2}_{p}&\dots&\gamma_{p}b^{p}_{p}\\
\end{pmatrix}.
\end{equation}
Ainsi, les conditions \ref{condmat} s'écrivent: 
\begin{equation}
\begin{cases}
\frac{e_{0}}{(p+1)Vol(B(0,1))}\;=\;\sum^{p}_{m=2}\alpha_{m}b^{0}_{m}\;=\;\sum^{p}_{m=2}\beta_{m}b^{1}_{m}\;=\;\gamma_{k}b^{k}_{k}\quad pour \quad k=2,...,p;\\
\sum^{p}_{m=2}\beta_{m}b^{0}_{m}\;=\;-\sum^{p}_{m=2}\alpha_{m}b^{1}_{m};\\
\gamma_{k}b^{0}_{k}\;=\;-\sum^{p}_{m=2}\alpha_{m}b^{k}_{m}\quad pour \quad k=2,...,p;\\
\gamma_{k}b^{1}_{k}\;=\;-\sum^{p}_{m=2}\beta_{m}b^{k}_{m}\quad pour \quad k=2,...p;\\
\gamma_{k}b^{l}_{k}\;=\;-\gamma_{l}b^{k}_{l}\quad pour \quad 2\leq k<l\leq p.
\end{cases}
\end{equation}
La première ligne nous dit que les $\gamma_{k}$ sont inversibles et donne la valeur de $
b^{k}_{k}\;=\;\frac{\gamma^{-1}_{k}}{(p+1)Vol(B(0,1))},$ les troisième et quatrième permettent d'écrire $b^{0}_{k}\;=\;-\sum^{p}_{m=2}\gamma^{-1}_{k}\alpha_{m}b^{k}_{m}$ et $b^{1}_{k}\;=\;-\sum^{p}_{m=2}\gamma^{-1}_{k}\beta_{m}b^{k}_{m};$ la cinquième donne $b^{l}_{k}\;=\;-\gamma^{-1}_{k}\gamma_{l}b^{k}_{l};$\\
Il reste donc trois conditions données par les deux premières égalités de la première ligne, et par la deuxième ligne. Nous les écrivons: 
\begin{equation}
\begin{split}
\frac{e_{0}}{(p+1)vol(B(0,1))}\;&=\;-\sum^{p}_{m,m'=2}\alpha_{m}\gamma^{-1}_{m}\alpha_{m'}b^{m}_{m'}\\
&=\;-\sum^{p}_{m=2}\frac{(\alpha_{m}\gamma^{-1}_{m})^{2}}{(p+1)Vol(B(0,1))}\;-\;\sum_{m<m'}(\alpha_{m}\gamma^{-1}_{m}\alpha_{m'}\;-\;\alpha_{m'}\gamma^{-1}_{m'}\alpha_{m}\gamma^{-1}_{m}\gamma_{m'})b^{m}_{m'},\\
\end{split}
\end{equation}
et
\begin{equation}
\begin{split}
\frac{e_{0}}{(p+1)vol(B(0,1))}\;&=\;-\sum^{p}_{m,m'=2}\beta_{m}\gamma^{-1}_{m}\beta_{m'}b^{m}_{m'}\\
&=\;-\sum^{p}_{m=2}\frac{(\beta_{m}\gamma^{-1}_{m})^{2}}{(p+1)Vol(B(0,1))}\;-\;\sum_{m<m'}(\beta_{m}\gamma^{-1}_{m}\alpha_{m'}\;-\;\beta_{m'}\gamma^{-1}_{m'}\beta_{m}\gamma^{-1}_{m}\gamma_{m'})b^{m}_{m'},\\
\end{split}
\end{equation}
et aussi
\begin{equation}
\begin{split}
\sum^{p}_{m=2}\beta_{m}b^{0}_{m}\;&=\;-\sum^{p}_{m,m'=2}\beta_{m}\gamma^{-1}_{m}\alpha_{m'}b^{m}_{m'}\\
&=\;-\sum^{p}_{m=2}\frac{(\beta_{m}\gamma^{-1}_{m}\alpha_{m}\gamma^{-1}_{m})}{(p+1)Vol(B(0,1))}\;-\;\sum_{m<m'}(\beta_{m}\gamma^{-1}_{m}\alpha_{m'}\;-\;\beta_{m'}\gamma^{-1}_{m'}\alpha_{m}\gamma^{-1}_{m}\gamma_{m'})b^{m}_{m'}\\
&=\;-\sum^{p}_{m=2}\alpha_{m}b^{1}_{m}\\
&=\;\sum^{p}_{m,m'=2}\alpha_{m}\gamma^{-1}_{m}\beta_{m'}b^{m}_{m'}\\
&=\;\sum^{p}_{m=2}\frac{(\alpha_{m}\gamma^{-1}_{m}\beta_{m}\gamma^{-1}_{m})}{(p+1)Vol(B(0,1)}\;+\;\sum_{m<m'}(\alpha_{m}\gamma^{-1}_{m}\beta_{m'}\;-\;\alpha_{m'}\gamma^{-1}_{m'}\beta_{m}\gamma^{-1}_{m}\gamma_{m'})b^{m}_{m'}.
\end{split}
\end{equation}
Comme les éléments $b^{m}_{m'}$ de $A$ sont arbitraires, nous pouvons écrire ces trois conditions: 
\begin{equation}\label{C2b}
\begin{split}
\Big(e_{0}&\;+\;\sum^{p}_{m=2}(\alpha_{m}\gamma^{-1}_{m})^{2}\;,\;e_{0}\;+\;(\beta_{m}\gamma^{-1}_{m})^{2}\;,\;\sum^{p}_{m=2}(\alpha_{m}\gamma^{-1}_{m}\beta_{m}\gamma^{-1}_{m}\;+\;\beta_{m}\gamma^{-1}_{m}\alpha_{m}\gamma^{-1}_{m})\Big)\in\\ &I\Big[\Big((\alpha_{m'}\gamma^{-1}_{m'})(\alpha_{m}\gamma^{-1}_{m})\gamma_{m'}\;-\;(\alpha_{m}\gamma^{-1}_{m})\alpha_{m'}
\;,\;(\beta_{m'}\gamma^{-1}_{m'})(\beta_{m}\gamma^{-1}_{m})\gamma_{m'}\;-\;(\beta_{m}\gamma^{-1}_{m})\beta_{m'}
\;,\;\\
&(\beta_{m'}\gamma^{-1}_{m'})(\alpha_{m}\gamma^{-1}_{m})\gamma_{m'}\;+\;(\alpha_{m'}\gamma^{-1}_{m'})(\beta_{m}\gamma^{-1}_{m})\gamma_{m'}\;-\;(\beta_{m}\gamma^{-1}_{m})\alpha_{m'}\;-\;(\alpha_{m}\gamma^{-1}_{m})\beta_{m'}\Big),\;m<m'\Big],
\end{split}
\end{equation}
où $I\Big[\Big(a(m,m'),b(m,m'),c(m,m')\Big)\;:\; m<m'\Big]$ désigne le sous-module de $A^{3}$ engendré par les éléments $\Big(a(m,m'),b(m,m'),c(m,m')\Big),\;\;m<m',$ de $A^{3}.$\\
Dans ce cas, nous pouvons résumer les résultats  dans le théorème suivant. \\

\begin{theorem}
Dans ce cas particulier, la classe des fonctions à valeurs dans $A$ vérifiant les $q$ conditions de Cauchy (\ref{C2a}) admet une formule de Cauchy 
\begin{equation}
f(x)\;=\;\int_{\partial D}f(y)\sum^{q}_{m=1}\Big(\sum^{n}_{j=1}a^{j}_{m}\frac{\partial}{\partial y_{j}}\Big)\lrcorner \Big(\frac{\varphi_{m}(x,y)}{\left\|y-x\right\|^{n}}dV\Big) 
\end{equation}
avec 
\begin{equation}
\varphi_{m}(x,y)\;=\;\sum^{n}_{i=1}\sum^{p}_{l=0}b^{i,l}_{m}(y_{i}-x_{i})e_{l}
\end{equation}
si et seulement si la condition $(A)$ donnée ci-dessus par \ref{C2b} est vérifiée. 
\end{theorem}

Nous avons, précédemment, étudié le cas des algèbres de dimension 2.\\
Nous allons maintenant rechercher des exemples dans des algèbres de dimension 3 et 4.\\
Comme dans l'exemple de la dimension 2 précédemment étudié, nous allons considérer ici un domaine $D$ borné et à bord lisse dans l'algèbre $A$ et allons chercher des classes de fonctions \\
$f\;:\;\overline{D}\longrightarrow A$ admettant une formule de Cauchy \ref{CondCau} . Ces classes de fonctions seront définies par les $q$ conditions de Cauchy \ref{C2a}. \\
Nous rencontrerons ainsi des exemples d'algèbres possédant de telles classes de fonctions vérifiant une formule de Cauchy, et d'autres qui n'en possèdent point. 
\\

\textbf{Le cas de la dimension 3.}\\
\\

Si $A$ est une algèbre unitaire de dimension 3, nous noterons $(e_{0},e_{1},e_{2})$ une base de $A$ avec $e_{0}$ l'élément unité. Pour connaître le produit, il suffit de donner dans la base précédente 
$$e_{i}e_{j}=\sum_{k=0}^{2}\Gamma_{i,j}^{k}e_{k},\quad pour\;\;1\leq i,j\leq 2.$$
Nous allons considérer la classe de fonctions $f$ de $A$ dans $A$ vérifiant une condition de Cauchy 
\begin{equation}\label{Cau1}
a_{0}\frac{\partial f}{\partial x_{0}}+a_{1}\frac{\partial f}{\partial x_{1}}+a_{2}\frac{\partial f}{\partial x_{2}}=0\;\; avec\;\;a_{0},\; a_{1},\; a_{2}\in A,
\end{equation}
c'est-à-dire nous avons $q=1,$ et nous demander si cette classe vérifie une formule de Cauchy. En fait, nous allons démontrer le résultat suivant.\\

\begin{theorem}
Il n'y a pas d'algèbre associative unitaire $A$ de dimension 3, munie d'une classe de fonctions définie par la condition \ref{Cau1} admettant une formule de Cauchy \ref{CondCau}. Dans le cas commutatif, le résultat reste vrai, même si l'on n'a pas l'associativité. 
\end{theorem}

Nous allons d'abord considérer le cas  d'une algèbre commutative qui est très simple. 
Les égalités \ref{cond12} s'écrivent 
\begin{equation*}\label{cond14}
	\begin{cases}
		a_{0}b_{0}=\frac{e_{0}}{3Vol}\\
		a_{1}b_{1}=\frac{e_{0}}{3Vol}\\
		a_{2}b_{2}=\frac{e_{0}}{3Vol}\\
		a_{0}b_{1}+a_{1}b_{0}=0\\
		a_{0}b_{2}+a_{2}b_{0}=0\\
		a_{1}b_{2}+a_{2}b_{1}=0.
	\end{cases}
\end{equation*}
Les trois premières équations s'écrivent $b_{0}=\frac{a_{0}^{-1}}{3Vol},\;\; b_{1}=\frac{a_{1}^{-1}}{3Vol},\;\;b_{2}=\frac{a_{2}^{-1}}{3Vol}$ et on remporte dans les trois dernières

\begin{equation*}
	\begin{cases}
a_{0}a_{1}^{-1}+a_{1}a_{0}^{-1}=0\\
a_{0}a_{2}^{-1}+a_{2}a_{0}^{-1}=0\\
a_{1}a_{2}^{-1}+a_{2}a_{1}^{-1}=0,
	\end{cases}
\end{equation*}
et on multiplie chaque ligne respectivement par $a_{0}a_{1},\;\; a_{0}a_{2}\;\; et\;\; a_{1}a_{2},$ pour obtenir
\begin{equation*}
	\begin{cases}
a_{0}^{2}+a_{1}^{2}=0\\
a_{0}^{2}+a_{2}^{2}=0\\
a_{1}^{2}+a_{2}^{2}=0,
	\end{cases}
\end{equation*}
d'où l'on tire $a_{0}^{2}=-a_{1}^{2}=-a_{2}^{2}=a_{1}^{2}=a_{2}^{2},$ donc $a_{0}^{2}=a_{2}^{2}=0,$ et, de même, $a_{0}^{2}=0.$ Mais ceci est impossible, car si $a_{0}^{2}=0,$ par exemple, alors la première ligne du premier de ces trois  systèmes donne, en la multipliant par $a_{0}$ : $\frac{a_{0}}{3Vol}=a_{0}^{2}b_{0}=0.$ \\
\begin{rem}
	Nous n'avons pas utilisé l'associativité dans cette démonstration. 
\end{rem}

Il reste maintenant à envisager le cas d'une algèbre non commutative. \\

Pour une algèbre $A,$ la propriété la plus délicate à satisfaire dans le choix des coefficients de structure $\Gamma^{k}_{i,j}$ est l'associativité. Nous allons l'exprimer sous forme matricielle. Soit $\Gamma_{.k}=\Big(\Gamma^{i}_{j,k}\Big)$ la matrice $(3,3)$ qui comporte l'élément $\Gamma^{i}_{j,k}$ à la ligne d'indice $i$ et la colonne d'indice $j,$ pour $0\leq i,j\leq2.$ De même $\Gamma_{k.}=\Big(\Gamma^{i}_{k,j}\Big)$ comporte $\Gamma^{i}_{k,j}$ à la ligne $i$ et la colonne $j.$ Ainsi, l'indice de ligne est toujours en haut, et l'indice de colonne toujours en bas. On vérifie facilement que l'associativité s'écrit, pour tout $1\leq k,i\leq 2,$ 
$$\Gamma_{.k}\Gamma_{i.}=\Gamma_{i.}\Gamma_{.k}.$$
Ici, nous avons tenu compte, bien sûr, du fait que, $si\; i\; ou\; j\; ou\; k=0\; alors\; (e_{i}e_{j})e_{k}=e_{i}(e_{j}
e_{k})$ est toujours vérifié.\\
Avec les 4 couples $(k,i),$ nous obtenons ainsi 4 équations qui traduisent l'égalité de deux matrices $(3,3).$\\ Les équations qui correspondent aux trois termes sur la diagonale se simplifient un peu et s'écrivent 
\begin{equation}
\begin{cases}
\Gamma^{0}_{1,k}\Gamma^{1}_{i,0}+\Gamma^{0}_{2,k}\Gamma^{2}_{i,0}=\Gamma^{0}_{i,1}\Gamma^{1}_{0,k}
+\Gamma^{0}_{i,2}\Gamma^{2}_{0,k}\\
\Gamma^{1}_{0,k}\Gamma^{0}_{i,1}+\Gamma^{1}_{2,k}\Gamma^{2}_{i,1}=\Gamma^{1}_{i,0}\Gamma^{0}_{1,k}
+\Gamma^{1}_{i,2}\Gamma^{2}_{1,k}\\
\Gamma^{2}_{0,k}\Gamma^{0}_{i,2}+\Gamma^{2}_{1,k}\Gamma^{1}_{i,2}=\Gamma^{2}_{i,0}\Gamma^{0}_{2,k}
+\Gamma^{2}_{i,1}\Gamma^{1}_{2,k},
\end{cases}
\end{equation}
soit aussi
\begin{equation}
\begin{cases}
\Gamma^{0}_{1,k}\Gamma^{1}_{i,0}-\Gamma^{1}_{0,k}\Gamma^{0}_{i,1}=\Gamma^{2}_{0,k}\Gamma^{0}_{i,2}
-\Gamma^{0}_{2,k}\Gamma^{2}_{i,0}\\
\Gamma^{0}_{1,k}\Gamma^{1}_{i,0}-\Gamma^{1}_{0,k}\Gamma^{0}_{i,1}=\Gamma^{1}_{2,k}\Gamma^{2}_{i,1}
-\Gamma^{2}_{1,k}\Gamma^{1}_{i,2}\\
\Gamma^{2}_{0,k}\Gamma^{0}_{i,2}
-\Gamma^{0}_{2,k}\Gamma^{2}_{i,0}=\Gamma^{1}_{2,k}\Gamma^{2}_{i,1}
-\Gamma^{2}_{1,k}\Gamma^{1}_{i,2}
\end{cases}
\end{equation}
ou encore 
\begin{equation}
\Gamma^{0}_{1,k}\Gamma^{1}_{i,0}-\Gamma^{1}_{0,k}\Gamma^{0}_{i,1}=\Gamma^{1}_{2,k}\Gamma^{2}_{i,1}
-\Gamma^{2}_{1,k}\Gamma^{1}_{i,2}=\Gamma^{2}_{0,k}\Gamma^{0}_{i,2}
-\Gamma^{0}_{2,k}\Gamma^{2}_{i,0},
\end{equation}
qui s'écrit aussi, en tenant compte du fait que $\Gamma^{i}_{0,j}=\Gamma^{i}_{j,0}=\delta^{i}_{j}$ 
\begin{equation}
\Gamma^{0}_{1,k}\delta_{1}^{i}-\delta_{1}^{k}\Gamma^{0}_{i,1}=\Gamma^{1}_{2,k}\Gamma^{2}_{i,1}
-\Gamma^{2}_{1,k}\Gamma^{1}_{i,2}=\delta_{2}^{k}\Gamma^{0}_{i,2}
-\Gamma^{0}_{2,k}\delta_{2}^{i}.
\end{equation}
Nous introduisons les notations 
\begin{equation}
\Gamma^{i}=\Gamma^{i}_{1,2}-\Gamma^{i}_{2,1}\quad \forall i=0,1,2;
\end{equation}
qui sont nulles dans  le cas commutatif, et aussi 
\begin{equation}
S^{i}=\Gamma^{i}_{1,2}+\Gamma^{i}_{2,1},
\end{equation}
donnons alors au couple $(i,k)$ les quatre valeurs qu'il est susceptible de prendre, et obtenons 
\begin{equation}\label{conddim3}
\begin{cases}
\Gamma^{2}_{1,1}\Gamma^{1}=0\\
\Gamma^{1}_{2,2}\Gamma^{2}=0\\
\Gamma^{0}_{1,2}=\Gamma^{1}_{2,2}\Gamma^{2}_{1,1}-\Gamma^{2}_{1,2}\Gamma^{1}_{1,2}\\
\Gamma^{0}_{2,1}=\Gamma^{2}_{1,1}\Gamma^{1}_{2,2}-\Gamma^{1}_{2,1}\Gamma^{2}_{2,1}.\\
\end{cases}
\end{equation}
Nous remarquons, ce qui est attendu, que le système de ces quatre équations reste inchangé si nous intervertissons les indices $1$ et $2.$ Les deux premières équations sont échangées, et les deeux dernières aussi.\\
Mais il reste six équations sur les neuf qui traduisent l'égalité des matrices $(3,3).$ Nous les explicitons. Pour tout $i,k=1,2,$ 
\begin{equation}
\begin{cases}
\Gamma^{0}_{0,k}\Gamma^{0}_{i,1}+\Gamma^{0}_{1,k}\Gamma^{1}_{i,1}+\Gamma^{0}_{2,k}\Gamma^{2}_{i,1}
=\Gamma^{0}_{i,0}\Gamma^{0}_{1,k}+\Gamma^{0}_{i,1}\Gamma^{1}_{1,k}+\Gamma^{0}_{i,2}\Gamma^{2}_{1,k}\\
\Gamma^{0}_{0,k}\Gamma^{0}_{i,2}+\Gamma^{0}_{1,k}\Gamma^{1}_{i,2}+\Gamma^{0}_{2,k}\Gamma^{2}_{i,2}=
\Gamma^{0}_{i,0}\Gamma^{0}_{2,k}+\Gamma^{0}_{i,1}\Gamma^{1}_{2,k}+\Gamma^{0}_{i,2}\Gamma^{2}_{2,k}\\
\Gamma^{1}_{0,k}\Gamma^{0}_{i,0}+\Gamma^{1}_{1,k}\Gamma^{1}_{i,0}+\Gamma^{1}_{2,k}\Gamma^{2}_{i,0}=
\Gamma^{1}_{i,0}\Gamma^{0}_{0,k}+\Gamma^{1}_{i,1}\Gamma^{1}_{0,k}+\Gamma^{1}_{i,2}\Gamma^{2}_{0,k}\\
\Gamma^{1}_{0,k}\Gamma^{0}_{i,2}+\Gamma^{1}_{1,k}\Gamma^{1}_{i,2}+\Gamma^{1}_{2,k}\Gamma^{2}_{i,2}=
\Gamma^{1}_{i,0}\Gamma^{0}_{2,k}+\Gamma^{1}_{i,1}\Gamma^{1}_{2,k}+\Gamma^{1}_{i,2}\Gamma^{2}_{2,k}\\
\Gamma^{2}_{0,k}\Gamma^{0}_{i,0}+\Gamma^{2}_{1,k}\Gamma^{1}_{i,0}+\Gamma^{2}_{2,k}\Gamma^{2}_{i,0}
=\Gamma^{2}_{i,0}\Gamma^{0}_{0,k}+\Gamma^{2}_{i,1}\Gamma^{1}_{0,k}+\Gamma^{2}_{i,2}\Gamma^{2}_{0,k}\\
\Gamma^{2}_{0,k}\Gamma^{0}_{i,1}+\Gamma^{2}_{1,k}\Gamma^{1}_{i,1}+\Gamma^{2}_{2,k}\Gamma^{2}_{i,1}
=\Gamma^{2}_{i,0}\Gamma^{0}_{1,k}+\Gamma^{2}_{i,1}\Gamma^{1}_{1,k}+\Gamma^{2}_{i,2}\Gamma^{2}_{1,k}.
\end{cases}
\end{equation}
Les troisième et cinquième équations du système ci-dessus sont trivialement vérifiées pour toutes les valeurs du couple $(i,k).$ Il reste alors quatre équations pour chacune des quatre valeurs du couple $(i,k).$ Nous explicitons ces 16 équations et nous obtenons  
\begin{equation}\label{conddim3bis}
\begin{cases}
\Gamma^{0}_{2,1}\Gamma^{2}_{1,1}=\Gamma^{0}_{1,2}\Gamma^{2}_{1,1}\\
\Gamma^{0}_{1,2}\Gamma^{1}_{2,1}+\Gamma^{0}_{2,2}\Gamma^{2}_{2,1}=\Gamma^{0}_{2,1}\Gamma^{1}_{1,2}+\Gamma^{0}_{2,2}\Gamma^{2}_{1,2}\\
\Gamma^{0}_{1,2}\Gamma^{1}_{1,1}+\Gamma^{0}_{2,2}\Gamma^{2}_{1,1}=\Gamma^{0}_{1,1}\Gamma^{1}_{1,2}+\Gamma^{0}_{1,2}\Gamma^{2}_{1,2}\\
\Gamma^{0}_{1,1}\Gamma^{1}_{2,1}+\Gamma^{0}_{2;1}\Gamma^{2}_{2,1}=\Gamma^{0}_{2,1}\Gamma^{1}_{1,1}+\Gamma^{0}_{2,2}\Gamma^{2}_{1,1}\\
\Gamma^{0}_{1,1}\Gamma^{1}_{1,2}+\Gamma^{0}_{2,1}\Gamma^{2}_{1,2}=\Gamma^{0}_{1,1}\Gamma^{1}_{2,1}+
\Gamma^{0}_{1,2}\Gamma^{2}_{2,1}\\
\Gamma^{0}_{1,2}\Gamma^{1}_{2,2}=\Gamma^{0}_{2,1}\Gamma^{1}_{2,2}\\
\Gamma^{0}_{1,2}\Gamma^{1}_{1,2}+\Gamma^{0}_{2,2}\Gamma^{2}_{1,2}=\Gamma^{0}_{1,1}\Gamma^{1}_{2,2}+
\Gamma^{0}_{1,2}\Gamma^{2}_{2,2}\\
\Gamma^{0}_{1,1}\Gamma^{1}_{2,2}+\Gamma^{0}_{2,1}\Gamma^{2}_{2,2}=\Gamma^{0}_{2,1}\Gamma^{1}_{2,1}+
\Gamma^{0}_{2,2}\Gamma^{2}_{2,1}\\
\Gamma^{0}_{1,2}+\Gamma^{1}_{1,1}\Gamma^{1}_{1,2}+\Gamma^{1}_{2,1}\Gamma^{2}_{1,2}=\Gamma^{0}_{2,1}+\Gamma^{1}_{1,1}\Gamma^{1}_{2,1}+\Gamma^{1}_{1,2}\Gamma^{2}_{2,1}\\
\Gamma^{1}_{1,2}\Gamma^{1}_{2,2}=\Gamma^{1}_{2,1}\Gamma^{1}_{2,2}\\
\Gamma^{1}_{1,2}\Gamma^{1}_{1,2}+\Gamma^{1}_{2,2}\Gamma^{2}_{1,2}=\Gamma^{0}_{2,2}+
\Gamma^{1}_{1,1}\Gamma^{1}_{2,2}+\Gamma^{1}_{1,2}\Gamma^{2}_{2,2}\\
\Gamma^{0}_{2,2}+\Gamma^{1}_{1,1}\Gamma^{1}_{2,2}+\Gamma^{1}_{2,1}\Gamma^{2}_{2,2}=\Gamma^{1}_{2,1}\Gamma^{1}_{2,1}+\Gamma^{1}_{2,2}\Gamma^{2}_{2,1}\\
\Gamma^{2}_{2,1}\Gamma^{2}_{1,1}=\Gamma^{2}_{1,2}\Gamma^{2}_{1,1}\\
\Gamma^{0}_{2,1}+\Gamma^{2}_{1,2}\Gamma^{1}_{2,1}+\Gamma^{2}_{2,2}\Gamma^{2}_{2,1}=\Gamma^{0}_{1,2}+\Gamma^{1}_{1,2}\Gamma^{2}_{2,1}+\Gamma^{2}_{2,2}\Gamma^{2}_{1,2}\\
\Gamma^{0}_{1,1}+\Gamma^{2}_{1,2}\Gamma^{1}_{1,1}+\Gamma^{2}_{2,2}\Gamma^{2}_{1,1}=\Gamma^{2}_{1,1}\Gamma^{1}_{1,2}+\Gamma^{2}_{1,2}\Gamma^{2}_{1,2}\\
\Gamma^{2}_{1,1}\Gamma^{1}_{2,1}+\Gamma^{2}_{2,1}\Gamma^{2}_{2,1}=\Gamma^{0}_{1,1}+\Gamma^{1}_{1,1}\Gamma^{2}_{2,1}+\Gamma^{2}_{2,2}\Gamma^{2}_{1,1}.
\end{cases}
\end{equation}
Les équations 15, 16, 11 et 12 du système ci-dessus s'écrivent aussi 
\begin{equation}
\begin{cases}
\Gamma^{0}_{1,1}&=\Gamma^{2}_{1,2}(\Gamma^{2}_{1,2}-\Gamma^{1}_{1,1})+\Gamma^{2}_{1,1}(\Gamma^{1}_{1,2}-\Gamma^{2}_{2,2})\\
&=\Gamma^{2}_{2,1}(\Gamma^{2}_{2,1}-\Gamma^{1}_{1,1})+\Gamma^{2}_{1,1}(\Gamma^{1}_{2,1}-\Gamma^{2}_{2,2})\\
\Gamma^{0}_{2,2}&=\Gamma^{1}_{1,2}(\Gamma^{1}_{1,2}-\Gamma^{2}_{2,2})+\Gamma^{1}_{2,2}(\Gamma^{2}_{1,2}-\Gamma^{1}_{1,1})\\
&=\Gamma^{1}_{2,1}(\Gamma^{1}_{2,1}-\Gamma^{2}_{2,2})+\Gamma^{1}_{2,2}(\Gamma^{2}_{2,1}-\Gamma^{1}_{1,1}),
\end{cases}
\end{equation}
ou encore 
\begin{equation}\label{conddim3ter}
\begin{cases}
\Gamma^{0}_{1,1}=\Gamma^{2}_{1,2}(\Gamma^{2}_{1,2}-\Gamma^{1}_{1,1})+\Gamma^{2}_{1,1}(\Gamma^{1}_{1,2}-\Gamma^{2}_{2,2})\\
\Gamma^{0}_{2,2}=\Gamma^{1}_{1,2}(\Gamma^{1}_{1,2}-\Gamma^{2}_{2,2})+\Gamma^{1}_{2,2}(\Gamma^{2}_{1,2}-\Gamma^{1}_{1,1})\\
\Gamma^{2}(S^{2}-\Gamma^{1}_{1,1})+\Gamma^{1}\Gamma^{2}_{1,1}=0\\
\Gamma^{1}(S^{1}-\Gamma^{2}_{2,2})+\Gamma^{2}\Gamma^{1}_{2,2}=0.
\end{cases}
\end{equation}

Par ailleurs, les équations 9 et 14 donnent
\begin{equation}\label{conddim3qua}
\begin{cases}
\Gamma^{0}&=\Gamma^{2}_{1,2}\Gamma^{1}-\Gamma^{1}_{1,2}\Gamma^{2}-\Gamma^{1}_{1,1}\Gamma^{1}\\
&=\Gamma^{1}_{21}\Gamma^{2}-\Gamma^{2}_{2,1}\Gamma^{1}-\Gamma^{2}_{2,2}\Gamma^{2}.
\end{cases}
\end{equation}

Ainsi, une algèbre $A$ de dimension 3 doit satisfaire, pour être associative, à la fois les conditions \ref{conddim3} et \ref{conddim3bis}. \\
Nous remarquons, comme prévu, que le système des ces 20 équations reste inchangé si nous intervertissons les indices 1 et 2. Pour mettre ce fait en évidence, réécrivons ce système en plaçant une équation dans une ligne impaire, et l'équation qui s'en déduit par interversion des indices $1$ et $2$ dans la ligne paire qui suit immédiatement. De façon précise, nous écrivons les équations 1 et 2 de \ref{conddim3}, puis les équations 1, 6, 13 et 10 de \ref{conddim3bis}, suivies des deux dernières équations de \ref{conddim3}, des équations de \ref{conddim3ter}, puis \ref{conddim3qua}, et enfin les équations 3, 8, 7, 4, 5 et 2 de \ref{conddim3bis}. Nous obtenons
\begin{equation}\label{systsym}
\begin{cases}
\Gamma^{2}_{1,1}\Gamma^{1}=0\\
\Gamma^{1}_{2,2}\Gamma^{2}=0\\
\Gamma^{2}_{1,1}\Gamma^{0}=0\\
\Gamma^{1}_{2,2}\Gamma^{0}=0\\
\Gamma^{2}_{1,1}\Gamma^{2}=0\\
\Gamma^{1}_{2,2}\Gamma^{1}=0\\
\Gamma^{0}_{1,2}=\Gamma^{1}_{2,2}\Gamma^{2}_{1,1}-\Gamma^{2}_{1,2}\Gamma^{1}_{1,2}\\
\Gamma^{0}_{2,1}=\Gamma^{2}_{1,1}\Gamma^{1}_{2,2}-\Gamma^{1}_{2,1}\Gamma^{2}_{2,1}\\
\Gamma^{0}_{1,1}=-\Gamma^{2}_{1,2}(\Gamma^{1}_{1,1}-\Gamma^{2}_{1,2})-\Gamma^{2}_{1,1}(\Gamma^{2}_{2,2}-\Gamma^{1}_{1,2})\\
\Gamma^{0}_{2,2}=-\Gamma^{1}_{2,1}(\Gamma^{2}_{2,2}-\Gamma^{1}_{2,1})-\Gamma^{1}_{2,2}(\Gamma^{1}_{1,1}-\Gamma^{2}_{2,1})\\
\Gamma^{2}(S^{2}-\Gamma^{1}_{1,1})+\Gamma^{1}\Gamma^{2}_{1,1}=0\\
\Gamma^{1}(S^{1}-\Gamma^{2}_{2,2})+\Gamma^{2}\Gamma^{1}_{2,2}=0\\
\Gamma^{0}=\Gamma^{2}_{1,2}\Gamma^{1}-\Gamma^{1}_{1,2}\Gamma^{2}-\Gamma^{1}_{1,1}\Gamma^{1}\\
\Gamma^{0}=\Gamma^{1}_{21}\Gamma^{2}-\Gamma^{2}_{2,1}\Gamma^{1}-\Gamma^{2}_{2,2}\Gamma^{2}\\
\Gamma^{1}_{2,2}\Gamma^{2}_{1,1}\Gamma^{2}+\Gamma^{2}_{1,1}S^{1}\Gamma^{1}-\Gamma^{2}_{1,1}\Gamma^{2}_{2,2}\Gamma^{1}=0\\
\Gamma^{2}_{1,1}\Gamma^{1}_{2,2}\Gamma^{1}+\Gamma^{1}_{2,2}S^{2}\Gamma^{2}-\Gamma^{1}_{2,2}\Gamma^{1}_{1;1}\Gamma^{2}=0\\
\Gamma^{2}_{1,2}[\Gamma^{1}(S^{1}-\Gamma^{2}_{2,2})+\Gamma^{1}_{2,2}\Gamma^{2}]=0\\
\Gamma^{1}_{2,1}[\Gamma^{2}(S^{2}-\Gamma^{1}_{1,1})+\Gamma^{2}_{1,1}\Gamma^{1}]=0\\
\Gamma^{0}_{1,1}\Gamma^{1}=-\Gamma^{1}_{2,2}\Gamma^{2}_{1,1}\Gamma^{2}-\Gamma^{2}_{1,2}\Gamma^{2}_{2,1}\Gamma^{1}\\
\Gamma^{0}_{2,2}\Gamma^{2}=-\Gamma^{2}_{1,1}\Gamma^{1}_{2,2}\Gamma^{1}-\Gamma^{1}_{2,1}\Gamma^{1}_{1,2}\Gamma^{2}.
\end{cases}
\end{equation}
\\

Dans ce système, les équations 7, 8, 9 et 10 peuvent être  considérées comme des définitions de $\Gamma^{0}_{1,2},\;\Gamma^{0}_{2,1},\;\Gamma^{0}_{1,1},\;\Gamma^{0}_{2,2}.$ 
\\
\\
Puisque le cas d'une algèbre commutative a déjà été étudié, nous considérons désormais, le cas non commutatif
\\
Alors, $\Gamma^{1}\;\; ou\;\; \Gamma^{2}$ est non nul. Donc, d'après les six premières équations, $\Gamma^{1}_{2,2}=\Gamma^{2}_{1,1}=0.$ Si, par exemple, $\Gamma^{2}\neq 0$ (la cas $\Gamma^{1}\neq 0$ se traite de façon tout à fait similaire), alors  \ref{systsym} s'écrit

\begin{equation}
\begin{cases}
\Gamma^{2}_{1,1}=\Gamma^{1}_{2,2}=0\\
\Gamma^{0}_{1,2}=-\Gamma^{1}_{1,2}\Gamma^{2}_{1,2}\\
\Gamma^{0}_{2,1}=-\Gamma^{1}_{2,1}\Gamma^{2}_{2,1}\\
\Gamma^{0}_{1,1}=\Gamma^{2}_{1,2}(\Gamma^{2}_{1,2}-\Gamma^{1}_{1,1})\\
\Gamma^{0}_{2,2}=\Gamma^{1}_{2,1}(\Gamma^{1}_{2,1}-\Gamma^{2}_{2,2})\\
S^{2}=\Gamma^{1}_{1,1}\\
\Gamma^{1}(S^{1}-\Gamma^{2}_{2,2})=0\\
\Gamma^{0}=\Gamma^{2}_{1,2}\Gamma^{1}-\Gamma^{1}_{1,2}\Gamma^{2}-\Gamma^{1}_{1,1}\Gamma^{1}=(\Gamma^{1}_{21}\Gamma^{2}-\Gamma^{2}_{2,1}\Gamma^{1}-\Gamma^{2}_{2,2}\Gamma^{2})\\
\Gamma^{1}S^{2}-\Gamma^{2}S^{1}=\Gamma^{1}_{1,1}\Gamma^{1}-\Gamma^{2}_{2,2}\Gamma^{2}\\
\Gamma^{2}_{1,2}\Gamma^{1}(S^{1}-\Gamma^{2}_{2,2})=0\\
\Gamma^{0}_{1,1}\Gamma^{1}=-\Gamma^{2}_{1,2}\Gamma^{2}_{2,1}\Gamma^{1}\\
\Gamma^{0}_{2,2}\Gamma^{2}=-\Gamma^{1}_{2,1}\Gamma^{1}_{1,2}\Gamma^{2}.
\end{cases}
\end{equation}
Ici, la ligne 9 traduit l'égalité dans la parenthèse de la ligne 8.\\
La ligne 9 s'écrit $\Gamma^{1}(S^{2}-\Gamma^{1}_{1,1})=\Gamma^{2}(S^{1}-\Gamma^{2}_{2,2}),$ ce qui, puisque le membre de gauche est nul, implique $S^{1}=\Gamma^{2}_{2,2}.$ La dixième ligne est satisfaite.\\
La quatrième ligne s'écrit 
$\Gamma^{0}_{1,1}=(\Gamma^{2}_{1,2})^{2}-\Gamma^{2}_{1,2}(\Gamma^{2}_{1,2}+\Gamma^{2}_{2,1})=-\Gamma^{2}_{1,2}\Gamma^{2}_{2,1},$ et donc l'avant dernière ligne est vérifiée. De façon analogue, la cinquième ligne implique $\Gamma^{0}_{2,2}=-\Gamma^{1}_{1,2}\Gamma^{1}_{2,1},$ et la dernière ligne est satisfaite. La huitième aussi, comme on le vérifie aisément. \\
Au total, dans le cas non commutatif, l'associativité se traduit par 
\begin{equation}\label{condnoncom}
\begin{cases}
\Gamma^{2}_{1,1}=\Gamma^{1}_{2,2}=0\\
\Gamma^{0}_{1,2}=-\Gamma^{1}_{1,2}\Gamma^{2}_{1,2}\\
\Gamma^{0}_{2,1}=-\Gamma^{1}_{2,1}\Gamma^{2}_{2,1}\\
\Gamma^{0}_{1,1}=-\Gamma^{2}_{1,2}\Gamma^{2}_{2,1}\\
\Gamma^{0}_{2,2}=-\Gamma^{1}_{1,2}\Gamma^{1}_{2,1}\\
S^{2}=\Gamma^{1}_{1,1}\\
S^{1}=\Gamma^{2}_{2,2}.
\end{cases}
\end{equation}
\\
\\
Maintenant que nous connaissons les conditions d'obtention d'une algèbre de dimension 3, nous allons regarder si la classe des fonctions vérifiant \ref{Cau1} possède une formule de Cauchy, c'est à dire satisfait les conditions \ref{cond12} explicitées en \ref{cond14}.
\begin{equation}
\begin{cases}
b_{0}=\frac{a_{0}^{-1}}{3Vol}\\
b_{1}=\frac{a_{1}^{-1}}{3Vol}\\
b_{2}=\frac{a_{2}^{-1}}{3Vol},
\end{cases}
\end{equation}
et, en reportant dans les trois dernières lignes,
\begin{equation}
\begin{cases}
a_{0}a_{1}^{-1}+a_{1}a_{0}^{-1}=0\\
a_{0}a_{2}^{-1}+a_{2}a_{0}^{-1}=0\\
a_{1}a_{2}^{-1}+a_{2}a_{1}^{-1}=0.
\end{cases}
\end{equation}
Ces équations sont équivalentes à $\forall (i,j)=(0,1),\;(0,2),\;(1,2),\quad e_{0}+(a_{j}a_{i}^{-1})^{2}=0,$ que nous écrirons

\begin{equation}\label{C6}
e_{0}+(c^{i}_{j})^{2}=0,
\end{equation}
en posant $c^{i}_{j}=a_{j}a_{i}^{-1}.$ Nous allons donc chercher $c^{i}_{j}$ vérifiant $e_{0}+(c^{i}_{j})^{2}=0.$ Il nous faudra ensuite trouver $a_{0},\; a_{1},\; a_{2}$ vérifiant 
\begin{equation}
\begin{cases}
a_{1}a_{0}^{-1}=c^{0}_{1}\\
a_{2}a_{0}^{-1}=c^{0}_{2}\\
a_{2}a_{1}^{-1}=c^{1}_{2}.\\
\end{cases}
\end{equation}
Les deux premières lignes donnent 
\begin{equation}
\begin{cases}
a_{1}=c^{0}_{1}a_{0}\\
a_{2}=c^{0}_{2}a_{0},
\end{cases}
\end{equation}
et la troisième ligne devient
\begin{equation}
c^{1}_{2}=c^{0}_{2}a_{0}a_{0}^{-1}(c^{0}_{1})^{-1}=c^{0}_{2}(c^{0}_{1})^{-1},
\end{equation}
soit

\begin{equation}\label{C8}
c^{1}_{2}c^{0}_{1}=c^{0}_{2}.
\end{equation}
Par conséquent, nous cherchons d'abord $c=ae_{0}+be_{1}+ce_{2}$ vérifiant $e_{0}+c^{2}=0.$ Sur les composantes, cela donne 
\begin{equation}\label{3g}
\begin{cases}
1+a^{2}+b^{2}\Gamma^{0}_{1,1}+bcS^{0}+c^{2}\Gamma^{0}_{2,2}=0\\
2ab+b^{2}\Gamma^{1}_{1,1}+bcS^{1}+c^{2}\Gamma^{1}_{2,2}=0\\
2ac++b^{2}\Gamma^{2}_{1,1}+bcS^{2}+c^{2}\Gamma^{2}_{2,2}=0.
\end{cases}
\end{equation}

Nous utilisons \ref{condnoncom} et \ref{3g} devient
\begin{equation}
\begin{cases}
1+a^{2}-b^{2}\Gamma^{2}_{1,2}\Gamma^{2}_{2,1}-bc[\Gamma^{1}_{1,2}\Gamma^{2}_{1,2}+\Gamma^{1}_{2,1}\Gamma^{2}_{2,1}]-c^{2}\Gamma^{1}_{1,2}\Gamma^{1}_{2,1}=0\\
b[2a+bS^{2}+cS^{1}]=0\\
c[2a+bS^{2}+cS^{1}]=0.
\end{cases}
\end{equation}
D'après la première équation, $b\; ou\; c$ est non nul. Donc, les deux dernières équations disent simplement que 

\begin{equation}
a=-\frac{bS^{2}+cS^{1}}{2},
\end{equation}
ce qui donne, en reportant dans la première équation
\begin{equation}
1+\frac{b^{2}(S^{2})^{2}+2bcS^{1}S^{2}+c^{2}(S^{1})^{2}}{4}-b^{2}\Gamma^{2}_{1,2}\Gamma^{2}_{2,1}-bc[\Gamma^{1}_{1,2}\Gamma^{2}_{1,2}+\Gamma^{1}_{2,1}\Gamma^{2}_{2,1}]-c^{2}\Gamma^{1}_{1,2}\Gamma^{1}_{2,1}=0,
\end{equation}
ou encore, par un calcul facile, 
\begin{equation}
b^{2}\frac{(\Gamma^{2})^{2}}{4}-bc\frac{\Gamma^{1}\Gamma^{2}}{2}+c^{2}\frac{(\Gamma^{1})^{2}}{4}=-1,
\end{equation}
soit 
\begin{equation}
(b\Gamma^{2}-c\Gamma^{1})^{2}=-1,
\end{equation}
ce qui est manifestement impossible. Ainsi, il n'y a pas de formule de Cauchy. \\
\\
\\

\textbf{Le cas de la dimension 4.}\\
\\
Nous allons considérer deux exemples d'algèbres de dimension 4, et, d'abord, l'exemple de l'espace 
$\mathfrak{M}_{2}(\R)$ des matrices carrées d'ordre 2.
\begin{ex}
$A=\mathfrak{M}_{2}(\R)$
\end{ex}
Nous choisirons, par exemple, comme base de $A=\mathfrak{M}_{2}(\R),$ la base $(e_{0},\;e_{1},\;e_{2},\;e_{3})$ suivante 
\begin{equation}
e_{0}=\begin{pmatrix}
1&0\\
0&1
\end{pmatrix},
\quad e_{1}=\begin{pmatrix}
0&1\\0&0
\end{pmatrix},
\quad e_{2}=\begin{pmatrix}
0&0\\1&0
\end{pmatrix},
\quad e_{3}=\begin{pmatrix}
1&0\\0&0
\end{pmatrix},
\end{equation}
à laquelle est associée la table suivante pour le produit 
\begin{center}
\begin{tabular}{|c||c|c|c|c|}

 \hline

 $\times$ & $\, e_0 \,$ & $e_1$   & $\, e_2 \,$ & $e_{3} $\\
 \hline
 \hline
 $e_0$    & $e_0$ & $e_1$   & $e_2$ & $e_{3}$ \\
 \hline
 $e_1$ &    $e_1$ & $ 0 $   & $e_{3}$ & $0$ \\
 \hline
 $e_2$ & $e_2$ & $ e_{0}-e_{3} $ & $0$ &$ e_{2} $\\
 \hline
 $e_3$ & $e_3$ & $ e_{1} $ & $0$ &$ e_{3} $\\
 \hline
 \end{tabular}
\end{center}
Nous calculons les éléments $A_{m,l}^{j}=a_{m}^{j}e_{l}$ pour $m=1,...,q$ où $q$ reste à fixer :
\begin{equation}\label{E2}
\begin{split}
A_{m,0}^{j}&=a_{m}^{j,0}e_{0}+a_{m}^{j,1}e_{1}+a_{m}^{j,2}e_{2}+a_{m}^{j,3}e_{3}\\
A_{m,1}^{j}&=a_{m}^{j,0}e_{1}+a_{m}^{j,2}(e_{0}-e_{3})+a_{m}^{j,3}e_{1}\\
A_{m,2}^{j}&=a_{m}^{j,0}e_{2}+a_{m}^{j,1}e_{3}\\
A_{m,3}^{j}&=a_{m}^{j,0}e_{3}+a_{m}^{j,2}e_{2}+a_{m}^{j,3}e_{3}.
\end{split}
\end{equation}
Nous allons considérer la classe des fonctions de $\mathfrak{M}_{2}(\R)$ dans $\mathfrak{M}_{2}(\R)$ qui satisfont les conditions de Cauchy \ref{C2a} et cherchons une formule de Cauchy pour cette classe.\\

Nous voulons donc traduire les relations \ref{Rel4}. Pour tout $i=0,1,2,3,$ fixé, le système $S_{i}$ s'écrit
\begin{equation}\label{E}
\sum_{m=0}^{q}\sum_{l=0}^{3}b^{i,l}_{m}A^{j,s}_{m,l}=c^{j,s}_{i}\;\; avec\;\; c^{i,s}_{i}=\delta^{s}_{0}\;\; et \;\; c^{i,s}_{j}=-c^{j,s}_{i}\;\; si\;\; i\neq j.
\end{equation}
Nous avons remarqué que les équations principales jouent un rôle décisif dans les formules obtenues. Aussi, nous allons faire plusieurs hypothèses sur ces équations principales.\\
\\

\fbox{Première hypothèse}\\
\\
Nous supposons que $q=1.$ Alors, la seule valeur permise pour le paramètre $m$ est $m=1.$ Par conséquent, nous supprimons, dans cette partie, l'indice $m.$  Notre hypothèse est que les équations principales sont au nombre de 4, et sont données par les indices 

\begin{equation}
I=\{(j_{0},s_{0})=(0,0),\;\;(j_{1},s_{1})=(0,1),\;\;(j_{2},s_{2})=(0,2),\;\;(j_{3},s_{3})=(0,3)\}
\end{equation}

et donc

\begin{equation}
\overline{I}=\{(j,s)\; :\; j=1,2,3,\;\; s=0,1,2,3\}.
\end{equation}
Ainsi, \ref{am2b} et \ref{Rel18} s'écrivent
\begin{equation}
\begin{split}
\forall  j\geqslant 1,\;\;\;c^{j,s}_{i}&=-\frac{1}{D_{0}}\sum_{l=0}^{3}(-1)^{l}c^{0,l}_{i}D^{0,l}_{j,s,i}\\
&=\frac{-1}{4D_{0}}D^{0,0}_{j,s,i}-\frac{1}{D_{0}}\sum_{l=1}^{3}(-1)^{l}c^{i,l}_{0}D^{0,l}_{j,s,i}\\
&=\frac{-1}{4D_{0}}D^{0,0}_{j,s,i}-\frac{1}{D_{0}^{2}}\sum_{l=1,2,3,\;\;l'=0,1,2,3}(-1)^{l+l'}D^{0,l}_{j,s,i}c^{0,l'}_{0}D^{0,l'}_{i,l,0}\\
&=\frac{-1}{4D_{0}}D^{0,0}_{j,s,i}-\frac{1}{4D_{0}^{2}}\sum_{l=1}^{3}(-1)^{l}D^{0,l}_{j,s,i}D^{0,0}_{i,l,0}
\end{split}
\end{equation}
et
\begin{equation}\label{E4}
\forall i,j\geqslant 1,\quad D^{0,0}_{j,s,i}+D^{0,0}_{i,s,j}+\frac{1}{D_{0}}\sum_{l=1}^{3}(-1)^{l}[D^{0,l}_{j,s,i}D^{0,0}_{i,l,0}+D^{0,l}_{i,s,j}D^{0,0}_{j,l,0}]=0
\end{equation}
tandis que \ref{cond16} devient

\begin{equation}
\begin{split}
\forall j\geqslant 1,\quad c^{j,s}_{j}=\frac{\delta_{0}^{s}}{4}&=\frac{-1}{D_{0}}\sum_{l=0}^{3}(-1)^{l}c^{0,l}_{j}D^{0,l}_{j,s,j}=\frac{1}{D_{0}}\sum_{l=0}^{3}(-1)^{l}c^{j,l}_{0}D^{0,l}_{j,s,j}\\
&=\frac{-1}{D_{0}^{2}}\sum_{l,l'=0}^{3}(-1)^{l+l'}D^{0,l}_{j,s,j}c^{0,l'}_{0}D^{0,l'}_{j,l,0}\\
&=\frac{-1}{4D_{0}^{2}}\sum_{l=0}^{3}(-1)^{l}D^{0,l}_{j,s,j}D^{0,0}_{j,l,0},
\end{split}
\end{equation}

ou
\begin{equation}\label{E6}
\sum_{l=0}^{3}(-1)^{l}D^{0,l}_{j,s,j}D^{0,0}_{j,l,0}=-D_{0}^{2}\delta_{0}^{s}.
\end{equation}
Nous avons ici les deux premiers type d'équations qui constituent \ref{S}. Quant au troisième type, les équations conditionnées à $j_{k}\neq j_{l},$ il n'y en a pas ici, car $k,l\leqslant r$ implique $j_{k}=j_{l}=0,$ et donc il est impossible d'avoir $((i,s_{k})=(0,s_{k})\in \overline{I}.$\\
Ici, le système \ref{S} qui donne les conditions d'existence d'une formule de Cauchy pour la classe des fonctions vérifiant les conditions de Caucchy \ref{C2a} s'écrit : $\forall i,j\geqslant 1,$
\begin{equation}
\begin{cases}
D^{0,0}_{j,s,i}+D^{0,0}_{i,s,j}+\frac{1}{D_{0}}\sum_{l=1}^{3}(-1)^{l}[D^{0,l}_{j,s,i}D^{0,0}_{i,l,0}+D^{0,l}_{i,s,j}D^{0,0}_{j,l,0}]=0\\
\sum_{l=0}^{3}(-1)^{l}D^{0,l}_{j,s,j}D^{0,0}_{j,l,0}=-D_{0}^{2}\delta_{0}^{s}.
\end{cases}
\end{equation}
Ce sont les conditions $(A)$ pour qu'existe une formule de Cauchy. \\
\ref{E} donne les équations principales :

\begin{equation}
\begin{cases}
\sum_{l=0}^{3}b^{i,l}A^{0,0}_{l}=c^{0,0}_{i}\\
\sum_{l=0}^{3}b^{i,l}A^{0,1}_{l}=c^{0,1}_{i}\\
\sum_{l=0}^{3}b^{i,l}A^{0,2}_{l}=c^{0,2}_{i}\\
\sum_{l=0}^{3}b^{i,l}A^{0,3}_{l}=c^{0,3}_{i},
\end{cases}
\end{equation}
c'est à dire $\forall s=0,...,3,$ d'après \ref{E2}, 
\begin{equation}
\begin{split}
&b^{i,o}[a^{0,0}e_{0}+a^{0,1}e_{1}+a^{0,2}e_{2}+a^{0,3}e_{3}]+b^{i,1}[a^{0,0}e_{1}+a^{0,2}(e_{0}-e_{3})+a^{0,3}e_{1}]\\
&\qquad +b^{i,2}[a^{0,0}e_{2}+a^{0,1}e_{3}]+b^{i,3}[a^{0,0}e_{3}+a^{0,2}e_{2}+a^{0,3}e_{3}]=cc^{0}_{i},
\end{split}
\end{equation}
soit, en séparant les quatre composantes

\begin{equation}
\begin{cases}
b^{i,0}a^{0,0}+b^{i,1}a^{0,2}&=c^{0,0}_{i}\\
b^{i,0}a^{0,1}+b^{i,1}(a^{0,0}+a^{0,3})&=c^{0,1}_{i}\\
b^{i,0}a^{0,2}  \qquad   \qquad  \qquad      +b^{i,2}a^{0,0}+b^{i,3}a^{0,2}&=c^{0,2}_{i}\\
b^{i,0}a^{0,3}       -b^{i,1}a^{0,2}\qquad +b^{i,2}a^{0,1}+b^{i,3}(a^{0,0}+a^{0,3})&=c^{0,3}_{i}.
\end{cases}
\end{equation}
Nous allons maintenant supposer (ce qui n'est pas vraiment restrictif) que la matrice de ce système de 4 équations de rang 4 est la matrice Identité, c'est à dire
\begin{equation}
\begin{cases}
a^{0,0}=1\\
a^{0,1}=a^{0,2}=a^{0,3}=0.
\end{cases}
\end{equation}
Le déterminant de ce système est donc $D_{0}=0.$ \\
\ref{E} donne les équations non principales, et nous obenons les déterminants caractéristiques qui fournissent les conditions pour que les systèmes $S_{i}$ aient des solutions : $\forall i=0,...,3,\;\; \forall j=1,2,3,$

\begin{equation}
c^{j,s}_{i}=A^{j,s}_{0}c^{0,0}_{i}+A^{j,s}_{1}c^{0,1}_{i}+A^{j,s}_{2}c^{0,2}_{i}+A^{j,s}_{3}c^{0,3}_{i},
\end{equation}

c'est à dire que l'on a $\forall j\geqslant 1,\;\;\forall k=0,...,3,\;\; D^{0,k}_{j,s,i}=(-1)^{k+1}A^{j,s}_{k}.$\\
Ainsi, les conditions d'obtention d'une formule de Cauchy \ref{E4} et \ref{E6} s'écrivent \\
$\forall i,j\geqslant 1,$
\begin{equation}\label{E8}
\begin{cases}
-A^{j,s}_{0}-A^{i,s}_{0}+\sum_{l=0}^{3}\big(A^{j,s}_{l}A^{i,l}_{0}+A^{i,s}_{l}A^{j,l}_{0}\Big)=0\\
\sum_{l=0}^{3}A^{j,s}_{l}A^{j,s}_{0}=-\delta_{0}^{s}.
\end{cases}
\end{equation}
\ref{E2} donne les composantes $A^{j,s}_{i}$ de $A^{j}_{i},$ et donc nous pouvons expliciter ces conditions. Nous explicitons d'abord la deuxième équation de \ref{E8} pour $s=0,...,3,$ et $j=1, 2.$

\begin{equation}\label{E12}
\begin{cases}
1+(a^{1,0})^{2}+a^{1,2}a^{1,1}=0\\
1+(a^{2,0})^{2}+a^{2,2}a^{2,1}=0\\
1+(a^{3,0})^{2}+a^{3,2}a^{3,1}=0\\
a^{1,1}[2a^{1,0}+a^{1,3}]=0\\
a^{2,1}[2a^{2,0}+a^{2,3}]=0\\
a^{3,1}[2a^{3,0}+a^{3,3}]=0.
\end{cases}
\end{equation}
D'après la première équation, $a^{1,1},\;\; a^{2,1}\;\; et\;\; a^{3,1}$ sont non nuls. Et donc

\begin{equation}
\begin{cases}\label{E10}
a^{1,3}=-2a^{1,0}\\
a^{2,3}=-2a^{2,0}\\
a^{3,3}=-2a^{3,0}.
\end{cases}
\end{equation}
Explicitons, maintenant, la première équation de \ref{E8} pour $(i,j,s)=(1,2,1),\;\;(1,3,1)\;\; et\;\; (2,3,1).$
\begin{equation}
\begin{cases}
-a^{2,1}-a^{1,1}+(a^{2,1}a^{1,0}+a^{1,1}a^{2,0})+((a^{2,0}+a^{2,3})a^{1,1}+(a^{1,0}+a^{1,3})a^{2,1})=0\\
-a^{3,1}-a^{1,1}+(a^{3,1}a^{1,0}+a^{1,1}a^{3,0})+((a^{3,0}+a^{3,3})a^{1,1}+(a^{1,0}+a^{1,3})a^{3,1})=0\\
-a^{3,1}-a^{2,1}+(a^{3,1}a^{2,0}+a^{2,1}a^{3,0})+((a^{3,0}+a^{3,3})a^{2,1}+(a^{2,0}+a^{2,3})a^{3,1})=0.
\end{cases}
\end{equation}
Nous reportons les valeurs de \ref{E10} dans cette équation et obtenons 
\begin{equation}
\begin{cases}
-a^{2,1}-a^{1,1}=0\\
-a^{3,1}-a^{1,1}=0\\
-a^{3,1}-a^{2,1}=0,
\end{cases}
\end{equation}
donc, nous avons
\begin{equation}
a^{1,1}=-a^{2,1}=-a^{3,1}=a^{2,1},
\end{equation}
et, en comparant le deuxième et le quatrième terme, nous déduisons que tous ces termes sont nuls. La nullité de $a^{1,1}$ est rigoureusement incompatible avec la première équation de \ref{E12}.\\
Ainsi, dans le cas étudié, nous n'avons pas de formule de Cauchy.\\
Nous avons déjà remarqué que, moins nous avons de conditions de Cauchy, plus il est difficile d'avoir une formule de Cauchy, mais plus est grande la classe des fonctions considérées. Nous avons donc un certain équilibre à respecter, et, en choisissant $q=1,$ nous avons trop peu de conditions de Cauchy. \\
\\

\fbox{Deuxième hypothèse}\\

Nous supposons $q=3,$ et voulons traduire les relations \ref{S4}. Dans le cas présent, elles s'écrivent, pour tout $i,j=0,1,2,3,$ 
\begin{equation}
\begin{split}
&\sum_{m=1}^{3}b_{m}^{i,0}\big(a_{m}^{j,0}e_{0}+a_{m}^{j,1}e_{1}+a_{m}^{j,2}e_{2}
+a_{m}^{j,3}e_{3}\big)+b_{m}^{i,1}\big(a_{m}^{j,0}e_{1}+a_{m}^{j,2}(e_{0}-e_{3})+a_{m}^{j,3}e_{1}\big)\\
&\qquad\qquad+b_{m}^{i,2}\big(a_{m}^{j,0}e_{2}+a_{m}^{j,1}e_{3}\big)
+b_{m}^{i,3}\big(a_{m}^{j,0}e_{3}+a_{m}^{j,2}e_{2}+a_{m}^{j,3}e_{3}\big)=c_{i}^{j},
\end{split}
\end{equation}
soit, en séparant les quatre composantes, 
\begin{equation}
\begin{cases}
b_{1}^{i,0}a_{1}^{j,0}+b_{1}^{i,1}a_{1}^{j,2}+b_{2}^{i,0}a_{2}^{j,0}+b_{2}^{i,1}a_{2}^{j,2}
+b_{3}^{i,0}a_{3}^{j,0}+b_{3}^{i,1}a_{3}^{j,2}=c_{i}^{j,0}\\
b_{1}^{i,0}a_{1}^{j,1}+b_{1}^{i,1}(a_{1}^{j,0}+a_{1}^{j,3})+b_{2}^{i,0}a_{2}^{j,1}+b_{2}^{i,1}(a_{2}^{j,0}+a_{2}^{j,3})
+b_{3}^{i,0}a_{3}^{j,1}+b_{3}^{i,1}(a_{3}^{j,0}+a_{3}^{j,3})=c_{i}^{j,1}\\
b_{1}^{i,0}a_{1}^{j,2}+b_{1}^{i,2}a_{1}^{j,0}+b_{1}^{i,3}a_{1}^{j,2}+b_{2}^{i,0}a_{2}^{j,2}
+b_{2}^{i,2}a_{2}^{j,0}+b_{2}^{i,3}a_{2}^{j,2}+b_{3}^{i,0}a_{3}^{j,2}+b_{3}^{i,2}a_{3}^{j,0}
+b_{3}^{i,3}a_{3}^{j,2}=c_{i}^{j,2}\\
b_{1}^{i,0}a_{1}^{j,3}-b_{1}^{i,1}a_{1}^{j,2}+b_{1}^{i,2}a_{1}^{j,1}+b_{1}^{i,3}(a_{1}^{j,0}
+a_{1}^{j,3})
+b_{2}^{i,0}a_{2}^{j,3}-b_{2}^{i,1}a_{2}^{j,2}+b_{2}^{i,2}a_{2}^{j,1}+b_{2}^{i,3}(a_{2}^{j,0}
+a_{2}^{j,3})
+b_{3}^{i,0}a_{3}^{j,3}\\
\qquad\qquad-b_{3}^{i,1}a_{3}^{j,2}+b_{3}^{i,2}a_{3}^{j,1}+b_{3}^{i,3}(a_{3}^{j,0}+a_{3}^{j,3})=c_{i}^{j,3}.
\end{cases}
\end{equation}
Pour chaque valeur de $i=0,1,2,3,$ nous avons un système $S_{i}$ de 16 équations (4 équations pour chaque valeur de $j=0,1,2,3$) à 12 inconnues $b_{m}^{i,l},\;\;m=1,2,3;\;\;l=0,1,2,3.$ Dans chaque système $S_{i},$ les inconnues et les seconds membres dépendent de $i,$ mais pas les coefficients dans les premiers membres. Les seize équations sont rangées dans l'ordre suivant : d'abord, les quatre équations relatives à $j=0,$ (dans l'ordre ci-dessus), puis les quatre correspondant à $j=1,$ suivies des quatre de $j=2,$ et enfin les quatre dernières avec $j=3.$ \\
Nous supposons, c'est là notre deuxième hypothèse, que les équations principales de chaque système $S_{i}$ sont les 12 équations qui correspondent à $j=0,1,2.$ Nous notons toujours $D_{0}$ le déterminant de ce système principal. Les quatre équations non principales correspondent à $j=3.$ Pour énoncer les quatre conditions de compatibilité, nous désignerons, par $D_{l}^{13+k},\;\;l=1,...,12;\;\;k=0,1,2,3,$ le déterminant obtenu, à partir de $D_{0},$ en supprimant la ligne $l$ 
et en rajoutant, en bas, la ligne $13+k.$ Nous obtenons les quatre conditions de compatibilité en développant les déterminants caractéristiques et obtenons : pour tout $k=0,1,2,3,$ 
\begin{equation}\label{condcomp}
\begin{split}
&c_{i}^{3,k}D_{0}-c_{i}^{2,3}D_{12}^{13+k}+c_{i}^{2,2}D_{11}^{13+k}-c_{i}^{2,1}D_{10}^{13+k}
+c_{i}^{2,0}D_{9}^{13+k}-c_{i}^{1,3}D_{8}^{13+k}+c_{i}^{1,2}D_{7}^{13+k}-c_{i}^{1,1}D_{6}^{13+k}\\
&\qquad\qquad+c_{i}^{1,0}D_{5}^{13+k}-c_{i}^{0,3}D_{4}^{13+k}+c_{i}^{0,2}D_{3}^{13+k}-c_{i}^{0,1}D_{2}^{13+k}
+c_{i}^{0,0}D_{1}^{13+k}=0,
\end{split}
\end{equation}
donc 
\begin{equation}
\begin{split}
c_{i}^{3,k}&=\frac{1}{D_{0}}\big[c_{i}^{2,3}D_{12}^{13+k}-c_{i}^{2,2}D_{11}^{13+k}+c_{i}^{2,1}D_{10}^{13+k}
-c_{i}^{2,0}D_{9}^{13+k}+c_{i}^{1,3}D_{8}^{13+k}-c_{i}^{1,2}D_{7}^{13+k}+c_{i}^{1,1}D_{6}^{13+k}\\
&\qquad\qquad-c_{i}^{1,0}D_{5}^{13+k}+c_{i}^{0,3}D_{4}^{13+k}-c_{i}^{0,2}D_{3}^{13+k}+c_{i}^{0,1}D_{2}^{13+k}
-c_{i}^{0,0}D_{1}^{13+k}\big]
\end{split}
\end{equation}
ou, sous forme matricielle, 
\begin{equation}\label{dim4a}
c_{i}^{3}=M_{1}c_{i}^{2}+M_{2}c_{i}^{1}+M_{3}c_{i}^{0},
\end{equation}
en notant
\begin{equation}
\begin{split}
&c_{i}^{j}=\begin{pmatrix}
c_{i}^{j,0}\\c_{i}^{j,1}\\c_{i}^{j,2}\\c_{i}^{j,3}
\end{pmatrix},
\quad M_{1}=\frac{1}{D_{0}}\begin{pmatrix}
-D_{9}^{13}&+D_{10}^{13}&-D_{11}^{13}&+D_{12}^{13}\\
-D_{9}^{14}&+D_{10}^{14}&-D_{11}^{14}&+D_{12}^{14}\\
-D_{9}^{15}&+D_{10}^{15}&-D_{11}^{15}&+D_{12}^{15}\\
-D_{9}^{16}&+D_{10}^{16}&-D_{11}^{16}&+D_{12}^{16}
\end{pmatrix},\\
& M_{2}=\frac{1}{D_{0}}\begin{pmatrix}
-D_{5}^{13}&+D_{6}^{13}&-D_{7}^{13}&+D_{8}^{13}\\
-D_{5}^{14}&+D_{6}^{14}&-D_{7}^{14}&+D_{8}^{14}\\
-D_{5}^{15}&+D_{6}^{15}&-D_{7}^{15}&+D_{8}^{15}\\
-D_{5}^{16}&+D_{6}^{16}&-D_{7}^{16}&+D_{8}^{16}
\end{pmatrix},
\quad M_{3}=\frac{1}{D_{0}}\begin{pmatrix}
-D_{1}^{13}&+D_{2}^{13}&-D_{3}^{13}&+D_{4}^{13}\\
-D_{1}^{14}&+D_{2}^{14}&-D_{3}^{14}&+D_{4}^{14}\\
-D_{1}^{15}&+D_{2}^{15}&-D_{3}^{15}&+D_{4}^{15}\\
-D_{1}^{16}&+D_{2}^{16}&-D_{3}^{16}&+D_{4}^{16}
\end{pmatrix}. 
\end{split}
\end{equation}
Conformément à \ref{cond4}, nous posons alors, pour $i\leq 2,\;\;c_{i}^{i}=\frac{1}{4}\begin{pmatrix}
1\\0\\0\\0
\end{pmatrix}:=\frac{1}{4}\delta_{0},$ et choisissons arbitrairement $c_{i}^{j}$ pour $0\leq j<i\leq 2$ (c'est-à-dire $c_{1}^{0},\;c_{2}^{0},\;c_{2}^{1}).$ \\
on en déduit, de façon à satisfaire \ref{cond4}, $c_{j}^{i}=-c_{i}^{j}$ pour $0\leq j<i\leq 2$ (c'est-à-dire $c_{0}^{1},\;c_{0}^{2},\;c_{1}^{2}).$ On calcule alors, à l'aide de \ref{dim4a},
\begin{equation}
\begin{split}
c_{0}^{3}&=M_{1}c_{0}^{2}+M_{2}c_{0}^{1}+M_{3}c_{0}^{0}=-M_{1}c_{2}^{0}-M_{2}c_{1}^{0}+\frac{1}{4}M_{3}\delta_{0}\\
c_{1}^{3}&=M_{1}c_{1}^{2}+M_{2}c_{1}^{1}+M_{3}c_{1}^{0}=-M_{1}c_{2}^{1}+\frac{1}{4}M_{2}\delta_{0}+M_{3}c_{1}^{0}\\
c_{2}^{3}&=M_{1}c_{2}^{2}+M_{2}c_{2}^{1}+M_{3}c_{2}^{0}=\frac{1}{4}M_{1}\delta_{0}+M_{2}c_{2}^{1}+M_{3}c_{2}^{0}.
\end{split}
\end{equation}
On en déduit, d'après \ref{cond4}, 
\begin{equation}
\begin{split}
c_{3}^{0}&=-c_{0}^{3}=M_{1}c_{2}^{0}+M_{2}c_{1}^{0}-\frac{1}{4}M_{3}\delta_{0}\\
c_{3}^{1}&=-c_{1}^{3}=M_{1}c_{2}^{1}-\frac{1}{4}M_{2}\delta_{0}-M_{3}c_{1}^{0}\\
c_{3}^{2}&=-c_{2}^{3}=-\frac{1}{4}M_{1}\delta_{0}-M_{2}c_{2}^{1}-M_{3}c_{2}^{0},
\end{split}
\end{equation}
et l'on calcule, d'après \ref{dim4a},
\begin{equation}
\begin{split}
c_{3}^{3}&=M_{1}c_{3}^{2}+M_{2}c_{3}^{1}+M_{3}c_{3}^{0}\\
&=-\frac{1}{4}M_{1}^{2}\delta_{0}-M_{1}M_{2}c_{2}^{1}-M_{1}M_{3}c_{2}^{0}+M_{2}M_{1}c_{2}^{1}-\frac{1}{4}M_{2}^{2}\delta_{0}-M_{2}M_{3}c_{1}^{0}+M_{3}M_{1}c_{2}^{0}+M_{3}M_{2}c_{1}^{0}-\frac{1}{4}M_{3}^{2}\delta_{0}\\
&=-\frac{1}{4}\big(M_{1}^{2}+M_{2}^{2}+M_{3}^{2}\big)\delta_{0}+\big(M_{2}M_{1}-M_{1}M_{2}\big)c_{2}^{1}+\big(M_{3}M_{1}-M_{1}M_{3}\big)c_{2}^{0}+\big(M_{3}M_{2}-M_{2}M_{3}\big)c_{1}^{0}, 
\end{split}
\end{equation}
dont on veut, pour que \ref{cond4} soit toujours vérifié, que la valeur soit $\frac{\delta_{0}}{4},$ c'est-à-dire 
\begin{equation}\label{dim4e}
\frac{\delta_{0}}{4}+\frac{1}{4}\big(M_{1}^{2}+M_{2}^{2}+M_{3}^{2}\big)\delta_{0}=\big(M_{2}M_{1}-M_{1}M_{2}\big)c_{2}^{1}+\big(M_{3}M_{1}-M_{1}M_{3}\big)c_{2}^{0}+\big(M_{3}M_{2}-M_{2}M_{3}\big)c_{1}^{0}. 
\end{equation}
Puisque $c_{2}^{1},\;c_{2}^{0},\;c_{1}^{0}$ peuvent être choisis arbitrairement, cela signifie simplement que 
\begin{equation}\label{dim4b}
\delta_{0}+\big(M_{1}^{2}+M_{2}^{2}+M_{3}^{2}\big)\delta_{0}\in Im\big(M_{2}M_{1}-M_{1}M_{2}\big)+Im\big(M_{3}M_{1}-M_{1}M_{3}\big)+Im\big(M_{3}M_{2}-M_{2}M_{3}\big),
\end{equation}
où, si $M\in \mathfrak{M}_{4}(\R),$ nous notons $Im\big(M\big)=\big\{MX\;:\; X\in \R^{4}\big\}.$\\
\\
Nous supposons maintenant que, dans les membres de gauche des 12 premières équations du système $S_{i}$ ci-dessus (c'est-à-dire dans les membres de gauche des équations principales du système $S_{i}$), tous les coefficients $a_{i}^{j,k}$ sont nuls, sauf $a_{1}^{0,0},\;\;a_{2}^{1,0}\;\;et\;\;a_{3}^{2,0}$ qui sont non nuls. La matrice du système principal est alors diagonale, et son déterminant est $D_{0}=(a_{1}^{0,0})^{4}(a_{2}^{1,0})^{4}(a_{3}^{2,0})^{4}\neq 0.$ \\
Nous calculons alors les $D_{l}^{13+k}:$ 
\begin{equation}\label{dim4d}
\begin{split}
si\;\;9\leq l\leq 12,\quad D_{l}^{13+k}&=\frac{(-1)^{l}D_{0}}{a_{3}^{2,0}}A_{3,3-(12-l)}^{3,k}\\
si\;\;5\leq l\leq 8,\quad D_{l}^{13+k}&=\frac{(-1)^{l}D_{0}}{a_{2}^{1,0}}A_{2,3-(8-l)}^{3,k}\\
si\;\;1\leq l\leq 4,\quad D_{l}^{13+k}&=\frac{(-1)^{l}D_{0}}{a_{1}^{0,0}}A_{1,3-(4-l)}^{3,k},
\end{split}
\end{equation}
et, par conséquent, 
\begin{equation}
\begin{split}
M_{3}=\frac{1}{a_{1}^{0,0}}\begin{pmatrix}
a_{1}^{3,0}&a_{1}^{3,1}&a_{1}^{3,2}&a_{1}^{3,3}\\
a_{1}^{3,2}&a_{1}^{3,0}+a_{1}^{3,3}&0&-a_{1}^{3,2}\\
0&0&a_{1}^{3,0}&a_{1}^{3,1}\\
0&0&a_{1}^{3,2}&a_{1}^{3,0}+a_{1}^{3,3}
\end{pmatrix},
\quad M_{2}=\frac{1}{a_{2}^{1,0}}\begin{pmatrix}
a_{2}^{3,0}&a_{2}^{3,1}&a_{2}^{3,2}&a_{2}^{3,3}\\
a_{2}^{3,2}&a_{2}^{3,0}+a_{2}^{3,3}&0&-a_{2}^{3,2}\\
0&0&a_{2}^{3,0}&a_{2}^{3,1}\\
0&0&a_{2}^{3,2}&a_{2}^{3,0}+a_{2}^{3,3}
\end{pmatrix},\\
M_{1}=\frac{1}{a_{3}^{2,0}}\begin{pmatrix}
a_{3}^{3,0}&a_{3}^{3,1}&a_{3}^{3,2}&a_{3}^{3,3}\\
a_{3}^{3,2}&a_{3}^{3,0}+a_{3}^{3,3}&0&-a_{3}^{3,2}\\
0&0&a_{3}^{3,0}&a_{3}^{3,1}\\
0&0&a_{3}^{3,2}&a_{3}^{3,0}+a_{3}^{3,3}
\end{pmatrix}.
\end{split}
\end{equation}
On calcule alors $M_{3}M_{2},$ puis $M_{3}M_{2}-M_{2}M_{3},$ et l'on trouve 
\begin{equation}
M_{3}M_{2}-M_{2}M_{3}=\frac{1}{a_{3}^{2,0}a_{2}^{1,0}}\begin{pmatrix}
\widetilde{A}&C&-B&-2\widetilde{A}\\
-B&-\widetilde{A}&0&B\\
0&0&\widetilde{A}&C\\
0&0&-B&-\widetilde{A}
\end{pmatrix},
\end{equation}
où 
\begin{equation}\label{rellin}
\begin{cases}
\widetilde{A}=a_{1}^{3,1}a_{2}^{3,2}-a_{2}^{3,1}a_{1}^{3,2}\\
B=a_{1}^{3,2}a_{2}^{3,3}-a_{2}^{3,2}a_{1}^{3,3}\\
C=a_{1}^{3,1}a_{2}^{3,3}-a_{2}^{3,1}a_{1}^{3,3}.
\end{cases}
\end{equation}
\begin{rem}\label{rq2}
\end{rem}
Remarquons, d'une part, que $\widetilde{A},\;B,\;C$ sont alors liés : en effet, si, par exemple, $a_{1}^{3,1}\neq 0,$ alors, d'après les première et dernière équations, 
\begin{equation}
a_{2}^{3,2}=\frac{\widetilde{A}+a_{2}^{3,1}a_{1}^{3,2}}{a_{1}^{3,1}}\quad et\quad a_{2}^{3,3}=\frac{C+a_{2}^{3,1}a_{1}^{3,3}}{a_{1}^{3,1}},
\end{equation}
et, en reportant dans la deuxième équation,
\begin{equation}
B=a_{1}^{3,2}\frac{C+a_{2}^{3,1}a_{1}^{3,3}}{a_{1}^{3,1}}-a_{1}^{3,3}\frac{\widetilde{A}+a_{2}^{3,1}a_{1}^{3,2}}{a_{1}^{3,1}}=\frac{a_{1}^{3,2}C-a_{1}^{3,3}\widetilde{A}}{a_{1}^{3,1}}, 
\end{equation}
c'est-à-dire
\begin{equation}\label{rellin2}
\widetilde{A}a_{1}^{3,3}+Ba_{1}^{3,1}-Ca_{1}^{3,2}=0;
\end{equation}
et, d'autre part, que la réciproque est vraie : si $\widetilde{A},\;B,\;C$ sont trois nombres vérifiant \ref{rellin2}, alors on peut trouver trois nombres $a_{2}^{3,1},\;a_{2}^{3,2},\;a_{2}^{3,3}$ tels que \ref{rellin} soit vérifiée. En effet, on choisit arbitrairement $a_{2}^{3,1},$ et l'on pose $a_{2}^{3,2}=\frac{\widetilde{A}+a_{2}^{3,1}a_{1}^{3,2}}{a_{1}^{3,1}}$ et $a_{2}^{3,3}=\frac{C+a_{2}^{3,1}a_{1}^{3,3}}{a_{1}^{3,1}}.$ \\
\\

Nous calculons le déterminant de $M_{3}M_{2}-M_{2}M_{3},$ 
\begin{equation}
Det(M_{3}M_{2}-M_{2}M_{3})=\frac{1}{a_{2}^{1,0}a_{3}^{2,0}}\big(\widetilde{A}^{2}-BC\big)^{2},
\end{equation}
et utilisons la relation ci-dessus entre $\widetilde{A},\;B,\;C,$ pour obtenir 
\begin{equation}
Det(M_{3}M_{2}-M_{2}M_{3})=\frac{1}{a_{2}^{1,0}a_{3}^{2,0}(a_{1}^{3,2})^{2}}\big(a_{1}^{3,2}\widetilde{A}^{2}-a_{1}^{3,3}\widetilde{A}B-a_{1}^{3,1}B^{2}\big)^{2}
\end{equation}
où $\widetilde{A}$ et $B$ peuvent être choisis arbitrairement. Nous pouvons les choisir de façon à ce que le déterminant soit non nul, ou bien choisir $\widetilde{A}$ et $B$ quelconques (mais pas tous deux nuls) si $(a_{1}^{3,3})^{2}-4a_{1}^{3,2}a_{1}^{3,1}<0,$ ce qui assure aussi la non nullité du déterminant. Alors $Im(M_{3}M_{2}-M_{2}M_{3})=A=\mathfrak{M}_{2}(\R).$ \\
Par conséquent, la condition \ref{dim4b} est satisfaite, et nous pouvons calculer $c_{1}^{0},\;\;c_{2}^{0},\;\;c_{2}^{1}$ vérifiant \ref{dim4e}, puis les autres $c_{i}^{j}$ conformément aux formules ci-dessus. Nous sommes alors en mesure de résoudre les systèmes $(S_{i}),\;\;i=0,1,2,3$ ci-dessus, dont les solutions $(b_{m}^{i,l})$ permettent d'expliciter les $\varphi_{m}.$ Alors, la classe des fonctions vérifiant \ref{C2a} vérifie la formule de Cauchy \ref{CondCau}. \\
Nous pouvons faire des calculs tout à fait analogues avec $M_{1}M_{3}-M_{3}M_{1}$ et $M_{1}M_{2}-M_{2}M_{1};$ il suffit, pour obtenir une formule de Cauchy que l'un des trois déterminants soit non nul, mais ce n'est bien sûr pas nécessaire puisque la condition recherchée est \ref{dim4b} qui est la vraie condition $(A)$ pour que la classe des fonctions vérifiant les conditions de Cauchy \ref{C2a} admette une formule de Cauchy. \\
\\
\\

\begin{ex}
Algèbres de Clifford de dimension 4. 
\end{ex}
Pour une définition des algèbres de Clifford, nous renvoyons, par exemple à \cite{LS}. Soit $E$ un espace vectoriel muni d'une forme bilinéaire $g,$ et d'une base orthogonale (pour $g$) $(v_{1},v_{2}).$ Nous posons $g(v_{i},v_{i})=a_{i}.$ \\
Nous considérons l'espace vectoriel $A$ engendré par une base indicée par les sous ensembles de $\{1,2\}.$ Donc, nous notons $e_{0}=e_{\varnothing},\;\;e_{1}=e_{\{1\}},\;\;e_{2}=e_{\{2\}},\;\;e_{3}=e_{\{1,2\}}$ et appelons algèbre de Clifford associée à $g,$ l'algèbre de base $(e_{0},\;e_{1},\;e_{2},\;e_{3})$ dans laquelle la table du produit est donnée par 
\begin{center}
\begin{tabular}{|c||c|c|c|c|}

 \hline

 $\times$ & $\, e_0 \,$ & $e_1$   & $\, e_2 \,$ & $e_{3} $\\
 \hline
 \hline
 $e_0$    & $e_0$ & $e_1$   & $e_2$ & $e_{3}$ \\
 \hline
 $e_1$ &    $e_1$ & $ a_{1}e_{0} $   & $e_{3}$ & $a_{1}e_{2}$ \\
 \hline
 $e_2$ & $e_2$ & $ -e_{3} $ & $a_{2}e_{0}$ &$ -a_{2}e_{1} $\\
 \hline
 $e_3$ & $e_3$ & $ -a_{1}e_{2} $ & $a_{2}e_{1}$ &$ -a_{1}a_{2}e_{0} $\\
 \hline
 \end{tabular}
\end{center}
En vue de résoudre \ref{syst4}, nous calculons d'abord $A_{m,l}^{j}=a_{m}^{j}e_{l}$ pour $m=1,...,q,$
\begin{equation}
\begin{cases}
A_{m,0}^{j}=a_{m}^{j,0}e_{0}+a_{m}^{j,1}e_{1}+a_{m}^{j,2}e_{2}+a_{m}^{j,3}e_{3}\\
A_{m,1}^{j}=a_{m}^{j,0}e_{1}+a_{1}a_{m}^{j,1}e_{0}-a_{m}^{j,2}e_{3}-a_{1}a_{m}^{j,3}e_{2}\\
A_{m,2}^{j}=a_{m}^{j,0}e_{2}+a_{m}^{j,1}e_{3}+a_{2}a_{m}^{j,2}e_{0}+a_{2}a_{m}^{j,3}e_{1}\\
A_{m,3}^{j}=a_{m}^{j,0}e_{3}+a_{1}a_{m}^{j,1}e_{2}-a_{2}a_{m}^{j,2}e_{1}+a_{1}a_{2}a_{m}^{j,3}e_{0},
\end{cases}
\end{equation}
et, dès lors, résoudre le système \ref{syst4} revient à résoudre, pour chaque $i=0,1,2,3,$ 
\begin{equation}
\sum_{m=0}^{2}\sum_{l=0}^{3}b_{m}^{i,l}A_{m,l}^{j}=c_{i}^{j},\quad 0\leq j\leq 3,
\end{equation}
ce qui, donne le système $S_{i}$ de 16 équations à 12 inconnues $\forall j=0,...,3,$\\
\begin{equation}
\begin{cases}
\sum _{m=0}^{3}b^{i,0}_{m}a^{j,0}_{m}+b^{i,1}_{m}a^{j,1}_{m}+b^{i,2}_{m}a^{j,2}_{m}-b^{i,3}_{m}a_{1}a_{2}a^{j,3}_{m}=c^{j,0}_{i}\\
\sum _{m=0}^{3}b^{i,0}_{m}a^{j,1}_{m}+b^{i,1}_{m}a^{j,0}_{m}+b^{i,2}_{m}a_{2}a^{j,3}_{m}-b^{i,3}_{m}a_{1}a_{2}a^{j,2}_{m}=c^{j,1}_{i}\\
\sum _{m=0}^{3}b^{i,0}_{m}a^{j,2}_{m}-b^{i,1}_{m}a_{1}a^{j,3}_{m}+b^{i,2}_{m}a^{j,0}_{m}+b^{i,3}_{m}a^{j,1}_{m}=c^{j,2}_{i}\\
\sum _{m=0}^{3}b^{i,0}_{m}a^{j,3}_{m}-b^{i,1}_{m}a^{j,2}_{m}+b^{i,2}_{m}a^{j,1}_{m}+b^{i,3}_{m}a^{j,0}_{m}=c^{j,3}_{i}.
\end{cases}
\end{equation}
Nous supposons que les équations principales de chaque $S_{i}$ sont les équations qui correspondent à $j=1,2,3.$  \\
Pour que le système $S_{i}$ ait des solutions, nous écrivons, comme dans l'exemple précédent, deuxième hypothèse, la nullité des déterminants caractéristiques, et traduisons, comme dans cet exemple, le résultat sous forme matricielle : $\forall i=0,...,3,$
\begin{equation}
c^{0}_{i}=M_{1}c^{1}_{i}+M_{2}c^{2}_{i}+M_{3}c^{3}_{i}
\end{equation}
où $M_{1},\;M_{2},\;M_{3}$ sont définis de façon analogue à l'exemple cité. Nous en déduisons

\begin{equation}
\begin{split}
c^{0}_{0}&=M_{1}c^{1}_{0}+M_{2}c^{2}_{0}+M_{3}c^{0}_{0}\\
&=-M_{1}(M_{1}c^{1}_{1}+M_{2}c^{2}_{1}+M_{3}c^{3}_{1})-M_{2}(M_{1}c^{1}_{2}+M_{2}c^{2}_{2}+M_{3}c^{3}_{2})-M_{3}(M_{1}c^{1}_{3}+M_{2}c^{2}_{3}+M_{3}c^{3}_{3})\\
&=-(M_{1}^{2}+M_{2}^{2}+M_{3}^{2})\frac{\delta_{0}}{4}-(M_{2}M_{1}-M_{1}M_{2})c^{1}_{2}-(M_{3}M_{1}-M_{1}M_{3})c^{1}_{3}-(M_{3}M_{2}-M_{2}M_{3})c^{2}_{3},
\end{split}
\end{equation}
soit

\begin{equation}
\frac{\delta_{0}}{4}+(M_{1}^{2}+M_{2}^{2}+M_{3}^{2})\frac{\delta_{0}}{4}=-(M_{2}M_{1}-M_{1}M_{2})c^{1}_{2}-(M_{3}M_{1}-M_{1}M_{3})c^{1}_{3}-(M_{3}M_{2}-M_{2}M_{3})c^{2}_{3},
\end{equation}
ou encore

\begin{equation}
(I+M_{1}^{2}+M_{2}^{2}+M_{3}^{2})\frac{\delta_{0}}{4}\in Im(M_{2}M_{1}-M_{1}M_{2})+Im(M_{3}M_{1}-M_{1}M_{3})+Im(M_{3}M_{2}-M_{2}M_{3}).
\end{equation}
Nous avons ici la condition $(A)$ nécessaire et suffisante pour obtenir notre formule de Cauchy. \\
Comme précédemment, nous allons expliciter, dans un cas particulier, une condition suffisante. En effet, si l'une des matrices $M_{2}M_{1}-M_{1}M_{2},$ ou $M_{3}M_{1}-M_{1}M_{3}$ ou $M_{2}M_{3}-M_{3}M_{2}$ est inversible, alors la condition $(A)$ est vérifiée.\\
Nous supposons que la matrice donnée par les équations principales et les variables principales de $S_{i}$ est la matrice identité $I$ (dans l'exemple précédent, nous supposions quelle était diagonale). Nous calculons alors

\begin{equation}
\begin{split}
M_{1}=&\begin{pmatrix}
a^{0,0}_{1}&a^{0,1}_{1}&a^{0,2}_{1}&-a_{1}a^{0,3}_{1}\\
a^{0,1}_{1}&a^{0,0}_{1}&a_{2}a^{0,3}_{1}&-a_{2}a^{0,2}_{1}\\
a^{0,2}_{1}&-a_{1}a^{0,3}_{1}&a^{0,0}_{1}&a_{1}a^{0,1}_{1}\\
a^{0,3}_{1}&-a^{0,2}_{1}&a^{0,1}_{1}&-a^{0,0}_{1}
\end{pmatrix},\qquad
M_{2}=\begin{pmatrix}
a^{0,0}_{2}&a^{0,2}_{2}&a^{0,2}_{2}&-a_{1}a^{0,3}_{2}\\
a^{0,1}_{2}&a^{0,0}_{2}&a_{2}a^{0,3}_{2}&-a_{2}a^{0,2}_{2}\\
a^{0,2}_{2}&-a_{1}a^{0,3}_{2}&a^{0,0}_{2}&a_{1}a^{0,1}_{2}\\
a^{0,3}_{2}&-a^{0,2}_{2}&a^{0,1}_{2}&-a^{0,0}_{2}
\end{pmatrix}\\
&M_{3}=\begin{pmatrix}
a^{0,0}_{3}&a^{0,1}_{3}&a^{0,2}_{3}&-a_{1}a^{0,3}_{3}\\
a^{0,1}_{3}&a^{0,0}_{3}&a_{2}a^{0,3}_{3}&-a_{2}a^{0,2}_{3}\\
a^{0,2}_{3}&-a_{1}a^{0,3}_{3}&a^{0,0}_{3}&a_{1}a^{0,1}_{3}\\
a^{0,3}_{3}&-a^{0,2}_{3}&a^{0,1}_{3}&-a^{0,0}_{3}
\end{pmatrix}.
\end{split}
\end{equation}

Un calcul facile mais long  donne 

\begin{equation}
M_{3}M_{2}-M_{2}M_{3}=\begin{pmatrix}
0&a_{1}(1+a_{2})A_{2,3}&-(a_{1}+a_{2})A_{1,3}&a_{2}(1+a_{1})A_{1,2}\\
2a_{2}A_{2,3}&0&-2a_{2}A_{1,2}&a_{1}(1+a_{2})A_{1,3}\\
-2a_{1}A_{1,3}&(1+a_{1})A_{1,2}&0&a_{1}(1+a_{2})A_{2,3}\\
-2A_{1,2}&A_{1,3}-A_{0,2}&2a_{2}A_{2,3}&0\\
\end{pmatrix}
\end{equation}
où $A_{i,j}=a^{0,i}_{2}a^{0,j}_{3}-a^{0,j}_{2}a^{0,i}_{3}.$ On calcule le déterminant de cette matrice et l'on obtient
\begin{equation}
\begin{split}
&2\Big\{2a_{2}^{2}(1+a_{1})^{2}A_{1,2}^{4}+a_{1}^{2}(1+a_{1})(1+a_{2})(a_{1}+a_{2})A_{1,3}^{4}+2a_{1}^{2}a_{2}^{2}(1+a_{1})^{2}A_{2,3}^{4}\\
&\qquad -a_{1}(1+a_{1})[2a_{2}^{2}(1+a_{1})+(1+a_{2})(a_{1}+a_{2})]A_{1,2}^{2}A_{1,3}^{2}\\
&\qquad -2a_{1}^{2}[a_{1}(1+a_{1})^{2}+a_{2}(1+a_{2})^{2}]A_{1,2}^{2}A_{2,3}^{2}\\
&\qquad +a_{1}a_{2}(1+a_{2})[(1+a_{1})(a_{1}+a_{2})+2a_{1}^{2}(1+a_{2})]A_{1,3}^{2}A_{2,3}^{2}\Big\}.
\end{split}
\end{equation}
Lorsque ce déterminant est non nul, la condition $(A)$ est vérifiée. Et, il est facile de voir que cela peut être. Par exemple, si $a_{1}=0,$ alors ce déterminant est $4a_{2}^{2}A_{1,2}^{2}(A_{1,2}^{2}-a_{2}A_{2,3}^{2}).$\\
\\
Nous pouvons également donner des exemples plus explicites avec $q=3$ : choisissons $a_{1}$ et $a_{2}$ tels que $1+a_{1}+a_{2}-a_{1}a_{2}=0,$ alors nous avons $e_{0}+(e_{1})^{2}+(e_{2})^{2}+(e_{3})^{2}=0$ et, d'après ce qui précède , nous avons une formule de Cauchy pour les fonctions $A$-différentiables.\\

\section{Un autre exemple dans une alg\`ebre de dimension 4 et quelques applications}
Consid\'erons $H$ le corps des quaternions, $q$ un \'el\'ement de $H$ s'\'ecrit de mani\`ere habituelle $x_0+ix_1+jX_2+kx_3$ ce qui permet d'identifier $H$ \`a $\R^4$. Avec cette identification, on peut introduire l'op\'erateur de Cauchy-Fueter $\bar{\partial}_q:= \partial_{x_0}+i\partial_{x_1}+j\partial_{x_2}+k\partial_{x_3}$ et les fonctions de Cauchy-Fueter qui sont les fonctions $$f\ :\ \R^4\rightarrow H, \ \ (x_0,x_1,x_2,x_3)\rightarrow f_0+if_1+jf_2+kf_3$$ v\'erifiant $\bar{\partial}_q f=0$. Pour ces fonctions, les conditions de Cauchy correspondantes, avec les notations $(2.1)$ est le suivant: $a_{1}^{0}=1,\  a_{1}^1=i,\  a_{1}^2=j,\ a_{1}^3=k.$ Pour avoir une formule de Cauchy pour ces fonctions, il faut pouvoir r\'esoudre le syst\`eme $(2.11)$ soit trouver des $(b_{1}^i)$ tels que  :
$$a_{1}^ib_{1}^i={1\over 4Vol(B_{\R^4}(0,1))}\ (1),\ \ a_{1}^jb_{1}^i+a_{1}^ib_{1}^j=0\ (2),$$
ceci pour tous $i,j\in \{0,1,2,3\}$ et $i\not =j$. Les \'equations $(1)$ fixent les $b_{1}^i$, posons $\alpha={1\over 4Vol(B_{\R^4}(0,1))}$ alors $b_{1}^0=\alpha$, $b_{1}^1=-i\alpha$, $b_{1}^2=-j\alpha$, $b_{1}^3=-k\alpha$. Avec ces $(b_1^{i})$, il est ais\'e de v\'erifier que les \'equations $(2)$ sont satisfaites et nous avons donc un noyau de Cauchy.

Avant de d\'evelopper, notons $\alpha_n:={1\over 4nVol_{R^{4n}}}$. On introduit tr\`es souvent les fonctions de Cauchy-Fueter de plusieurs variables quaternioniques similairement aux fonctions holomorphes de plusieurs variables complexes, ce sont les fonctions $f$ de $\R^{4n}\rightarrow H$ v\'erifiant $\bar\partial_{q_l}=0$ avec $l\in\{0,\cdots, n-1\}$ o\`u $\bar\partial_{q_l}=\partial_{x_{4l}}+i\partial_{x_{4l+1}}+j\partial_{x_{4l+2}}+k\partial_{x_{4l+3}}$. Le syst\`eme $(2.11)$ correspondant admet alors les solutions suivantes $b_{l}^k=0$ si $k<4l$ ou $k>4l+3$ et $b_{l}^{4l}=1$, $b_{l}^{4l+1}=-i\alpha_n$, $b_{l}^{4l+2}=-j\alpha_n$, $b_{l}^{4l+3}=-k\alpha_n$. On obtient donc le noyau reproduisant correspondant habituel pour ces fonctions.

En fait, le passage des fonctions Cauchy-Fueter d'une variable \`a celles de plusieurs variables se fait simplement car on it\`ere les conditions de Cauchy d'une variable quaternionique  sur les variables r\'eelles $\{4l,\cdots, 4l+3\}$ avec $l\in\{0,\cdots, n-1\}$. Il est bien entendu que ceci est un fait g\'en\'eral, si $A$ est une alg\`ebre de dimension $k$ et si $(e_0,e_1,\cdots e_{k-1})$ est une base avec $e_0$ unitaire, on peut identifier $A$ \`a $\R^k$ en tant qu'espace vectoriel en consid\'erant les coordonn\'ees dans cette base; maintenant les fonctions de $A$ dans $A$ sont identifi\'ees aux fonctions de $R^k$ dans $A$. Consid\'erons des conditions de Cauchy qui admettent un noyau reproduisant pour les fonctions de $A$ dans $A$ v\'erifiant ces conditions, alors nous pouvons consid\'erer les fonctions de $A^n$ dans $A$ v\'erifiant les conditions de Cauchy d'une variable dans $A$ it\'er\'ees sur les variables r\'eelles prises par paquet de $k$ s\'epar\'ement. De la m\^eme mani\`ere que ci-dessus, le syst\`eme $(2.11)$ aura des solutions et donc les fonctions de $A^n$ dans $A$ v\'erifiant ces conditions de Cauchy admettront un noyau reproduisant. Ce r\'esultat n'est pas nouveau pour les fonctions de Cauchy-Fueter de plusieurs variables, il avait \'et\'e obtenu dans \cite{WW1} et \cite{WW2} par des m\'ethodes diff\'erentes.

\subsection{Th\'eor\`eme de Hartogs}

Regardons les fonctions Cauchy-Fueter d'une variable quaternionique, clairement au vu du noyau de Cauchy nous avons l'estimation suivante sur les d\'eriv\'ees de $f$:
$$\partial_{x_i}f(x)\leq MSup \vert f\vert_{\partial B(x,r)}/R.$$ En effet, la formule $(2.2)$ donne imm\'ediatement l'estimation suivante sur les d\'eriv\'ees du noyau:
$$\vert \partial _{x_i}K(y,x)\vert \leq  {1\over \vert \vert y-x\vert\vert ^{4}}.$$
Maintenant toute fonction de Cauchy-Fueter satisfait la formule :
$$\partial _{x_i}f(x)=\int_{\partial B_4(x,R)}f(y)\partial_{x_i}K(y,x),$$
et donc, en notant $ds$ l'\'el\'ement de mesure d'aire de la sph\'ere de dimension $3$ et de rayon R:
$$\vert \partial_{x_i} f(x)\vert\leq \int_{\partial B_4(x,R)}{\vert f(y)\vert\over \vert\vert y-x\vert\vert^{4}}ds(y) ,$$
ce qui donne l'estimation de Cauchy voulue.

Si $f$ est une fonction enti\`ere et born\'ee (c'est \`a dire, Cauchy-Fueter au voisinage de tout point de $H$) alors 
$$\vert \partial_{x_i}f(q)\vert \leq M \vert\vert f \vert\vert_{\infty}/R,$$ sur toute boule de rayon $R$, on en d\'eduit donc le th\'eor\`eme de Liouville pour ces fonctions.
Il en serait de m\^eme pour toutes les fonctions de $A$ dans $A$, avec $A$ une alg\`ebre unitaire et associative, v\'erifiant des conditions de Cauchy admettant un noyau reproduisant, en it\'erant le raisonnement au vu de la forme des noyaux (voir $(2.2)$).

Nous aurons d\'esormais besoin des deux d\'efinitions suivantes: 

\begin{defn} Soit $A$ une alg\`ebre unitaire et associative, si des conditions de Cauchy du type $(2.1)$ admettent un noyau de Cauchy pour les fonctions de $A^n$ dans $A$ v\'erifiant ces conditions, on dira que ces conditions de Cauchy sont admissibles dans $A^n$.
	
\end{defn}
\begin{defn} Soit $A$ une alg\`ebre unitaire, associative, et des conditions de Cauchy admissibles pour $A$. Les conditions de Cauchy sur $A^n$ issues de celles de $A$ par le proc\'ed\'e d\'ecrit ci-dessus seront appel\'ees conditions de Cauchy sur $A^n$ induites par $A$.
	
\end{defn}
Pour simplifier les \'ecritures, nous dirons qu'une fonction v\'erifiant des conditions de Cauchy admissibles est de Cauchy.

Il est clair d'apr\`es ce que l'on a vu que des conditions de Cauchy sur $A^n$ induites par $A$ sont admissibles pour $A^n$. Nous avons donc d\'emontrer le principe de Liouville pour les fonctions de Cauchy:

\begin{prop} Soit $A$ une alg\`ebre unitaire, associative, et des conditions de Cauchy admissibles sur $A$. Alors les fonctions de $A$ dans $A$ enti\`eres et born\'ees sur $A$ sont constantes.
\end{prop} 
De ce r\'esultat d\'ecoule le principe de Hartogs pour ces fonctions:

\begin{theorem} Soit $A$ une alg\`ebre unitaire, associative, et des conditions de Cauchy sur $A^n$ induites par $A$. Soit $f$ une fonction  de $A^n$ dans $A$ qui est de Cauchy sur une couronne contenant $\partial D$ le bord d'un domaine born\'e lisse, $D$, alors $f$ s'\'etend en une fonction de Cauchy sur $D$ tout entier.
\end{theorem}
Pour $A=H$ le corps des quaternions et les conditions de Cauchy-Fueter, le r\'esultat a \'et\'e  obtenu dans \cite{WW1} en r\'esolvant les \'equations de Cauchy-Fueter non homog\`enes \`a support compact: $\bar{\partial}_{q_l}u=f$ pour tout $l\in\{0,\cdots,n-1\}$ avec $u$ lisse et  \`a support compact, si $f$ est  lisse \`a support compact et v\'erifie les conditions de compatibilit\'e ad hoc; De plus, le support de $u$ est inclu dans le support de $f$.

Preuve: nous allons reproduire la preuve classique du cas holomorphe voir par exemple \cite{R}. Supposons que $D:=\{\rho < 0\}$ avec $\rho$ lisse au voisinage de $D$ alors consid\'erons $D_2:=\{\rho <\varepsilon \}$ et $D_1:=\{\rho < -\varepsilon \}$ avec $\varepsilon$ suffisamment petit et appliquons la formule de Cauchy sur la couronne $D_2/D_1$, on obtient:
$$f(x)=\int_{\partial D_2}f(y)K(y,x)-\int_{\partial D_1}f(y)K(y,x),  \ \forall x\in D_2/D_1.$$
Si nous pouvons d\'emontrer que $\int_{\partial D_1}f(y)K(y,x)=0$ pour tout $x$ dans le compl\'ementaire de $\bar D_1$ alors pour tout $x$ dans $D_2/D_1$, $f(x)=\int_{\partial D_2}f(y)K(y,x)$ et la fonction \`a droite de l'\'egalit\'e \'etant  de Cauchy sur $D$, car $\hat{d^{''}}K(y,x)=0$ hors de la diagonale (voir $(2.3)$), nous avons le r\'esultat souhait\'e.

Notons $x=(x^{'},x_n)$ o\`u $x^{'}$ est une variable de $A^{n-1}$ et $x_n$ la derni\`ere composante dans $A$ du vecteur $x\in A^n$; si $x^{'}$ est fix\'e suffisamment grand alors la fonction en $x$, $\int_{\partial D_1}f(y)K(y,x)$, est enti\`ere en $x_n$ pour les conditions de Cauchy sur $A$ qui induisent celles sur $A^n$. De plus, cette fonction tend vers z\'ero quand $x_n$ tend vers l'infini et elle est donc born\'ee sur $A$. D'apr\`es la propri\'et\'e de Liouville, elle est constante et donc nulle pour tout $y^{'}$ suffisamment grand, ce qui entra\^\i ne $\int_{\partial D_1}f(y)K(y,x)=0$ sur la composante non born\'ee de $A^n/\partial D_1$.

\subsection{Solution de Cauchy \`a support compact}

Consid\'erons des conditions de Cauchy sur $A^n$ induites par les conditions de Cauchy admissibles sur $A$. Si l' alg\`ebre $A$ est de dimension $p$ alors elle peut être repr\'esent\'ee par $\R^p$ muni d'une forme bilin\' eaire donn\'ee par le produit dans $A$ (le produit n'est pas forc\' ement  le produit composante par composante). Ainsi toute fonction de $A^n$ vers $A$  devient une fonction de $\R^{pn}$ dans $\R^p$. Introduisons les notations suivantes pour les conditions de Cauchy induites:
$$\sum_{i} a^{ij}_{m}{\partial f_{j}\over  \partial x_{i}}\ \ m=1, \cdots , p,\ \ 1\leq j \leq p,\ \ 1\leq i\leq np.$$ Dans ce qui suit, nous voulons r\'esoudre le syt\`eme d'EDP suivant :
$$\sum_{i} a^{ij}_{m}{\partial f_j\over  \partial x_{i}}=g_m\ \ m=1, \cdots , p.$$

Supposons que le tableau induit par ces conditions de Cauchy soit en involution (voir \cite{BM1} pour les d\'etails) alors les conditions de compatibilit\'es pour pouvoir résoudre le syst\`eme pr\'ec\'edent est \'egalement un syst\`eme d'EDP d'ordre $1$ \`a coefficients constants. Ces conditions de compatibilit\'es sont donn\'ees par le calcul de la torsion et not\'ees $tor_1(g)$ si $g$ est le second membre du syt\`eme de d\'epart (voir \cite{BM1} pour la d\'efinition pr\'ecise de la torsion). 

Nous allons prouver le th\'eor\`eme suivant:
\begin{theorem} Soit $A$ une alg\`ebre unitaire, associative, et des conditions de Cauchy sur $A^n$ induites par $A$ avec un tableau associ\'e en involution au sens de Cartan. Alors si $g$ est une fonction lisse de $\R^{pn}$ dans $\R^p$, \`a support compact et v\'erifiant $tor_1(g)=0$ alors le syst\`eme :
	$$\sum_{i} a^{ij}_{m}{\partial f_j\over  \partial x_{i}}=g_m\ \ m=1, \cdots , p.$$
	a une solution $f$ lisse \`a support compact.
	
	De plus, $g$ est donn\'ee par la formule:
	
	$$f(y)=\int_{B(0,r)}g(x)K(x,y), \ \ y\in A^n,$$
	o\`u $K$ est le noyau de Cauchy associ\'e aux conditions de Cauchy sur $A^n$ induite par $A$ et $B(0,r)$ est une boule de rayon suffisamment grand.
\end{theorem}
\begin{rem} Il n'est pas n\'ecessaire d'avoir le tableau associ\'e au syst\`eme des \'equations de Cauchy en involution. La seule diff\'erence est qu'alors la torsion peut \^etre un  op\'erateur \`a coefficients constants d'ordre sup\'erieur à 1 mais  le th\'eor\`eme reste valable pourvu que la torsion s'annule \'evidemment. C'est le cas pour les \'equations de Cauchy-Fueter, la torsion est un op\'erateur \`a coefficients constants d'ordre 2 et le th\'eor\`eme est vrai pour les formes \`a support compact v\'erifiant $tor_1(g)=0$. Ce r\'esultat pour les conditions de Cauchy-Fueter \'etait connu mais obtenu par d'autres m\' ethodes (voir \cite{WW1} et \cite{WW2}). Le cas du tableau en involution est celui qui se rapproche le plus du complexe de De Rham ou Dolbeault, et nous avons pris le parti de traiter uniquement ce cas.
\end{rem}

Preuve : D'apr\`es \cite{BM1}, sous les hypoth\`eses du th\'eor\`eme, il existe une solution $f$ au syst\`eme sur n'importe quel domaine convexe, par exemple sur n'importe quelle boule, $B(0,r)$. En g\'en\'eral, la solution n'est pas \`a support compact mais nous allons la corriger. Prenons $B(0,r)$ une boule contenant le support de $g$ alors sur la couronne $C(r,2r)$, la fonction $f$ est de Cauchy et elle se prolonge en une fonction de Cauchy $\tilde{f}$ sur $B(0,2r)$ d'apr\`es le th\'eor\`eme de Hartogs. Maintenant $f-\tilde f$ v\'erifie le syst\`eme car $\tilde f$ est de Cauchy et elle est \`a support compact dans $A^n$.

Il reste \`a d\'emontrer la formule: il suffit d'appliquer la formule de Cauchy (13) \`a la solution $f$ \`a support compact sur une boule $B(0,r)$ contenant son support.

\section{Conditions de Cauchy à coefficients variables.} 
Dans ce paragraphe, nous allons changer au minimum les notations du précédent paragraphe.\\
Soit encore $f\,:\,\overline{D}\longrightarrow A$ une fonction de classe $C^{1}$, à valeurs dans l'algèbre unitaire
$A.$ Soient aussi $a^{j}_{m}(x);\;\;j=1,...,n;\;\;m=1,...,q;$ des fonctions de classe $C^{1}$ de $\overline{D}$ dans $A.$ \\
Nous définissons, aussi, les opérateurs de Cauchy-Riemann: 
$$d''\;=\;\sum^{q}_{m=1}dx_{i_{m}}\sum^{n}_{j=1}\;\frac{\partial}{\partial x_{j}}a^{j}_{m}(x), \quad 1\leq i_{m}<i_{m+1}\leq n.$$
et
$$\widehat{d''}\;=\;\sum^{q}_{m=1}dx_{i_{m}}\sum^{n}_{j=1}\;a^{j}_{m}(x)\frac{\partial}{\partial x_{j}}, \quad 1\leq i_{m}<i_{m+1}\leq n.$$
En fait, pour définir la singularité sur la diagonale, nous n'avons pas besoin d'utiliser la distance euclidienne, comme nous avons fait ci-dessus, ni même une norme. Aussi, nous noterons:
$$\left|\left\|.\right|\right\|\;:\;\R^{n}\rightarrow \R$$
une application de $\R^{n}$ dans $\R$ telle que $\left|\left\|x\right|\right\|\geq 0,\; \;\left|\left\|\lambda x\right|\right\|\;=\;\left|\lambda\right|\;\left|\left\|x\right\|\right|,\quad \lambda\in\R,$ et $\left|\left\|x\right\|\right|=0\quad implique\quad x=0$.\\
Nous chercherons $K(y,x),$ solution fondamentale pour $\widehat{d''},$ sous la forme: 
$$K(y,x)\;=\;\frac{1}{\left|\left\|y-x\right|\right\|^{n}}\sum^{q}_{m=1}\widehat{dy_{i_{m}}}(-1)^{i_{m}-1}\varphi_{m}(x,y),$$ avec $\varphi_{m}$ fonction $C^{1},$ à valeurs dans $A$ et vérifiant $\varphi_{m}(x,x)=0$ pour tout $x\in D.$ Notant $dV\;=\;dy_{1}\wedge ...\wedge dy_{n},$ nous avons, hors de la diagonale: 
\begin{equation}
\begin{split}
\widehat{d''_{y}}K(y,x)\;&=\;dV\sum^{q}_{m=1}\sum^{n}_{j=1}a^{j}_{m}(y)\frac{\partial}{\partial y_{j}}\Big(\frac{\varphi_{m}(x,y)}{\left|\left\|y-x\right\|\right|^{n}}\Big)\\
&=\;\frac{dV}{\left|\left\|y-x\right\|\right|^{n}}\sum_{m,j}a^{j}_{m}(y)\Big[\frac{\partial \varphi_{m}(x,y)}{\partial y_{j}}\;-\;\frac{n}{2\left|\left\|y-x\right\|\right|^{2}}\frac{\partial\left|\left\|y-x\right\|\right|^{2}}{\partial y_{j}}\varphi_{m}(x,y)\Big].
\end{split}
\end{equation}
Donc, $\widehat{d''_{y}}K(y,x)=0,$ hors de $\Delta,$ si et seulement si: 
\begin{equation}\label{C1}
\sum_{m,j}a^{j}_{m}(y)\Big[\frac{\partial \varphi_{m}(x,y)}{\partial y_{j}}\;-\;\frac{n}{2\left|\left\|y-x\right\|\right|^{2}}\frac{\partial\left|\left\|y-x\right\|\right|^{2}}{\partial y_{j}}\varphi_{m}(x,y)\Big]\;=\;0.
\end{equation}
Nous allons, maintenant, étudier la singularité sur la diagonale.  
\begin{equation}
\begin{split}
\int_{D}\sum_{m,j}&\frac{\partial}{\partial y_{j}}\Big(f(y)a^{j}_{m}(y)\Big)\frac{\varphi_{m}(x,y)}{\left|\left\|y-x\right\|\right|^{n}}dV\;=\;\lim_{\epsilon\rightarrow 0}\int_{D-B(x,\epsilon)}\sum_{m,j}\frac{\partial}{\partial y_{j}}\Big(f(y)a^{j}_{m}(y)\Big)\frac{\varphi_{m}(x,y)}{\left|\left\|y-x\right\|\right|^{n}}dV\\
&=\;\lim_{\epsilon\rightarrow 0}\int_{D-B(x,\epsilon)}d\Big[\sum_{m,j}f(y)a^{j}_{m}(y)\frac{\varphi_{m}(x,y)}{\left|\left\|y-x\right\|\right|^{n}}(-1)^{j-1}\widehat{dy_{j}}\Big]-\sum_{m,j}f(y)a^{j}_{m}(y)\frac{\partial}{\partial y_{j}}\Big(\frac{\varphi_{m}(x,y)}{\left|\left\|y-x\right\|\right|^{n}}\Big)dV\\
&=\;\int_{\partial D}\sum_{m,j}f(y)a^{j}_{m}(y)\frac{\varphi_{m}(x,y)}{\left|\left\|y-x\right\|\right|^{n}}(-1)^{j-1}\widehat{dy_{j}}-\lim_{\epsilon\rightarrow 0}\int_{\partial B(x,\epsilon)}\sum_{m,j}f(y)a^{j}_{m}(y)\frac{\varphi_{m}(x,y)}{\left|\left\|y-x\right\|\right|^{n}}(-1)^{j-1}\widehat{dy_{j}}\\
&\qquad\qquad -\int_{D-B(x,\epsilon)}f(y)\widehat{d''_{y}}K(y,x)\\
&=\;\int_{\partial D}f(y)\sum_{m}\Big[\Big(\sum_{j}a^{j}_{m}(y)\frac{\partial}{\partial y_{j}}\Big)\rfloor\Big(\frac{\varphi_{m}(x,y)}{\left|\left\|y-x\right\|\right|^{n}}dV\Big)\Big]-f(x)Vol(B(0,1)\sum_{m,j}a^{j}_{m}(x)\frac{\partial \varphi_{m}(x,x)}{\partial y_{j}}.
\end{split}
\end{equation}
Comme l'on veut que ce dernier membre de droite soit $f(x),$ on doit prendre $\sum_{m,j}a^{j}_{m}(x)\frac{\partial\varphi_{m}(x,x)}{\partial y_{j}}\;=\;\frac{e_{0}}{Vol(B(0,1))}.$ En tenant compte de \ref{C1}, nous obtenons les conditions suivantes pour avoir une représentation intégrale: 
\begin{equation}\label{C2}
\begin{cases}
\sum^{q}_{m=1}\sum^{n}_{j=1}a^{j}_{m}(x)\frac{\partial\varphi_{m}(x,x)}{\partial y_{j}}\;=\;\frac{e_{0}}{Vol(B(0,1))}\\
\sum_{m,j}a^{j}_{m}(y)\Big[\frac{\partial \varphi_{m}(x,y)}{\partial y_{j}}\;-\;\frac{\partial log\left|\left\|y-x\right\|\right|^{n}}{\partial y_{j}}\varphi_{m}(x,y)\Big]\;=\;0\quad hors\quad de\quad \Delta.
\end{cases}
\end{equation}
Ces conditions sur $K(y,x)$ étant satisfaites, nous avons la formule de représentation intégrale 
\begin{equation}
f(x)\;=\;\int_{\partial D}f(y)\sum_{m}\Big[\Big(\sum_{j}a^{j}_{m}(y)\frac{\partial}{\partial y_{j}}\Big)\rfloor\Big(\frac{\varphi_{m}(x,y)}{\left|\left\|y-x\right\|\right|^{n}}dV\Big)\Big]-\int_{D}\sum_{m,j}\frac{\partial}{\partial y_{j}}\Big(f(y)a^{j}_{m}(y)\Big)\frac{\varphi_{m}(x,y)}{\left|\left\|y-x\right\|\right|^{n}}.
\end{equation}
Pour toute fonction $f$ de la classe des fonctions vérifiant les conditions de Cauchy 
\begin{equation}\label{conCauvar}
\sum^{n}_{j=1}\frac{\partial}{\partial y_{j}}\Big(f(y)a^{j}_{m}(y)\Big)\;=\;0\quad pour\quad tout\quad m=1,...,q,
\end{equation}
nous avons la formule de Cauchy 
\begin{equation}\label{repintvar}
f(x)\;=\;\int_{\partial D}f(y)\sum_{m}\Big[\Big(\sum_{j}a^{j}_{m}(y)\frac{\partial}{\partial y_{j}}\Big)\rfloor\Big(\frac{\varphi_{m}(x,y)}{\left|\left\|y-x\right\|\right|^{n}}dV\Big)\Big].
\end{equation}

Nous rappelons que $\varphi_{m}(x,x)=0$ pour tout $x\in D.$ \\
Par conséquent, nous posons: 
$$\varphi_{m}(x,y)\;=\;\sum^{n}_{j=1}b^{j}_{m}(x)(y_{j}-x_{j})\;+\;\Phi_{m}(x,y)$$
avec $\Phi_{m}(x,y)\;=\;O(\left|\left\|y-x\right\|\right|^{2}).$
Dès lors, \ref{C2} s'écrit: 
\begin{equation}\label{C3}
\begin{cases}
\sum^{q}_{m=1}\sum^{n}_{j=1}a^{j}_{m}(x)b^{j}_{m}(x)\;=\;\frac{e_{0}}{Vol(B(0,1))}\\
\begin{split}
\sum_{m,j}&a_{m}^{j}(y)\big[\frac{\partial\Phi_{m}(x,y)}{\partial y_{j}}\;-\;\frac{n}{2\left|\left\|y-x\right\|\right|^{2}}\frac{\partial\left|\left\|y-x\right\|\right|^{2}}{\partial y_{j}}\Phi_{m}(x,y)\big]\;=\\
&=\;-\sum_{m,j}\big(a^{j}_{m}(y)-a^{j}_{m}(x)\big)\big[b^{j}_{m}(x)-\frac{n}{2\left|\left\|y-x\right\|\right|^{2}}\frac{\partial\left|\left\|y-x\right\|\right|^{2}}{\partial y_{j}}\sum^{n}_{i=1}b^{i}_{m}(x)(y_{i}-x_{i})\big]\\
&\;-\;\frac{e_{0}}{Vol(B(0,1))}\;+\;\frac{n}{2\left|\left\|y-x\right\|\right|^{2}}\sum_{m,i,j}a^{j}_{m}(x)\big[\frac{\partial\left|\left\|y-x\right\|\right|^{2}}{\partial y_{j}}b^{i}_{m}(x)(y_{i}-x_{i})\big].
\end{split}
\end{cases}
\end{equation}
$\frac{\partial\Phi_{m}(x,y)}{\partial y_{j}}$ est un $O(\left|\left\|y-x\right\|\right|),\quad (a^{j}_{m}(y)-a^{j}_{m}(x))$ aussi ; donc si, dans la dernière équation du système précédent, nous identifions les termes d'ordre $0$ en $\left|\left\|y-x\right\|\right|$ d'une part (ils constituent la troisième ligne de la deuxième équation de \ref{C3}), et les termes d'ordre supérieur ou égal à un (les deux premières lignes de la deuxième équation) d'autre part, nous obtenons:
\begin{equation}\label{C4}
\begin{cases}
\sum^{q}_{m=1}\sum^{n}_{j=1}a^{j}_{m}(x)b^{j}_{m}(x)\;=\;\frac{e_{0}}{Vol(B(0,1))},\\
\left|\left\|y-x\right\|\right|^{2}\frac{e_{0}}{Vol(B(0,1))}\;=\;\frac{n}{2}\sum^{n}_{i,j=1}(y_{i}-x_{i})\frac{\partial\left|\left\|y-x\right\|\right|^{2}}{\partial y_{j}}\sum^{q}_{m=1}a^{j}_{m}(x)b^{i}_{m}(x),\\
\begin{split}
\sum_{m,j}&\big[\frac{\partial\Phi_{m}(x,y)}{\partial y_{j}}\;-\;\frac{n}{2\left|\left\|y-x\right\|\right|^{2}}\frac{\partial\left|\left\|y-x\right\|\right|^{2}}{\partial y_{j}}\Phi_{m}(x,y)\big]\;=\\
&=\;-\sum_{m,j}\big(a^{j}_{m}(y)-a^{j}_{m}(x)\big)\big[b^{j}_{m}(x)-\frac{n}{2\left|\left\|y-x\right\|\right|^{2}}\frac{\partial\left|\left\|y-x\right\|\right|^{2}}{\partial y_{j}}\sum^{n}_{i=1}b^{i}_{m}(x)(y_{i}-x_{i})\big].
\end{split}
\end{cases}
\end{equation}
La deuxième équation de \ref{C4} devient spécialement intéressante si nous prenons pour $\left|\left\|.\right\|\right|$ la norme euclidienne. Alors elle s'écrit: 
\begin{equation}
\begin{split}
\left\|y-x\right\|^{2}\frac{e_{0}}{Vol(B(0,1))}\;&=\;n\sum^{n}_{i,j=1}(y_{i}-x_{i})(y_{j}-x_{j})\sum_{m}a^{j}_{m}(x)b^{i}_{m}(x)\\
&=\,n\sum^{n}_{j=1}(y_{j}-x_{j})^{2}\sum_{m}a^{j}_{m}(x)b^{j}_{m}(x)\;\\
&+\;n\sum_{i<j}(y_{i}-x_{i})(y_{j}-x_{j})\sum_{m}\big(a^{j}_{m}(x)b^{i}_{m}(x)+a^{i}_{m}(x)b^{j}_{m}(x)\big),
\end{split}
\end{equation}
et, en identifiant les deux membres, on obtient: 
\begin{equation}\label{C5}
\begin{cases}
\sum^{q}_{m=1}a^{j}_{m}(x)b^{j}_{m}(x)\;=\;\frac{e_{0}}{nVol(B(0,1))}\\
\sum^{q}_{m=1}a^{j}_{m}(x)b^{i}_{m}(x)\;+\;a^{i}_{m}(x)b^{j}_{m}(x)\;=0,\quad si\quad i\neq j.
\end{cases}
\end{equation}
Nous retrouvons ici les conditions \ref{cond12}. \\
Dans le paragraphe 2, ces conditions nous ont conduits, dans le cas commutatif, aux conditions (A) (voir \ref{cond20}). De même ici, si nous appelons $\widetilde{\varphi_{m}}(x,y)\;=\;\sum^{n}_{j=1}b^{j}_{m}(x)(y_{j}-x_{j})$ la partie d'ordre 1 de $\varphi_{m}(x,y),$ alors, dans le cas commutatif, les mêmes calculs que ceux du paragraphe 2 nous mènent à la condition nécessaire (A). De façon plus précise, nous suivons les calculs effectués au paragraphe 2 dans le cas commutatif. Nous définissons 
$$\Psi^{j}(x,y)\;=\;Vol(B(0,1))\sum^{q}_{m=1}a^{j}_{m}(x)\widetilde{\varphi_{m}}(x,y)\;=\;\sum^{n}_{i=1}\,c^{j}_{i}(x)(y_{i}-x_{i}),$$ 
sommes conduits à  
$$\Psi^{j}(x,y)\;=\;(y_{j}-x_{j})\frac{e_{0}}{n}\;+\;\sum_{i\neq j}\,(y_{i}-x_{i})c^{j}_{i}(x)\quad avec \quad  c^{j}_{i}\;=\;-c^{i}_{j},$$
et amenés à résoudre le système, analogue au système \ref{syst1}, à coefficients fonctions de x:
\begin{equation}\label{syst2}
\Psi^{j}(x,y)\;=\;\sum^{q}_{m=1}\,a^{j}_{m}(x)\widetilde{\varphi_{m}}(x,y).
\end{equation}
Comme au paragraphe 2, nous appelons $D_{0}(x)$ un déterminant, extrait du système, d'ordre $q$ et inversible. Des calculs analogues au paragraphe 2 nous imposent les conditions nécessaires (A):
\begin{equation}\label{Avar}
\quad \quad (A):\quad\quad\quad\quad\quad\sum^{q}_{m=1}D^{l}_{m}(x)D^{k}_{m}(x)\;+\;\delta^{l}_{k}(D_{0}(x))^{2}\;=\;0\quad pour \quad tout\quad 1\leq k,l\leq n-q,
\end{equation}
et, dans le cas où elles sont satisfaites, nous donnent les solutions (voir \ref{cond8}):
\begin{equation}\label{solvar}
b^{m}_{i}(x)=\sum^{q}_{l=1}(-1)^{m+l}c^{l}_{i}(x)(D_{0}(x))^{-1}\widehat{D^{l}_{m}(x)}=(D_{0})^{-1}(x)\widetilde{D^{i}_{m}}(x)
\end{equation}
avec des notations analogues à celles du paragraphe 2. \\
\\
Nous avons donc le: 
\begin{theorem}
Soient $D$ un domaine borné à bord lisse de $\R^{n},$ des fonctions 
$$a^{j}_{m}(x):\overline{D}\rightarrow A;\;\; m=1,...,q;\;\;j=1,...,n,$$
de classe $C^{1}$ et à valeurs dans $A,$ et $f:\overline{D}\rightarrow A$ une fonction de classe $C^{1}$ sur $\overline{D}$ à valeurs dans $A$ vérifiant les conditions de Cauchy (\ref{conCauvar}).\\
Alors, la formule de représentation intégrale (\ref{repintvar}) est vérifiée avec $\varphi_{m}(x,y)$ vérifiant les conditions (\ref{C2}), ou, dans le cas de la norme euclidienne, les conditions \ref{C5}.\\
Pour que ces conditions puissent être vérifiées, il faut, dans le cas commutatif, que les conditions $(A),$ données par (\ref{Avar}) soient vérifiées, et les solutions sont données par (\ref{solvar}).\\
\end{theorem}

Nous allons appliquer ceci à des exemples, en choisissant, d'abord, pour les coefficients $a^{j}_{m}$ des conditions de Cauchy des polynômes de degré 1 en $y-x.$ 
\begin{ex}
Conditions de Cauchy-Riemann à coefficients de degré 1, dans le cas A commutatif. 
\end{ex}
Ici encore, nous choisirons A commutatif, afin de mener les calculs à l'aide de l'outil que sont les déterminants. 
Nous choisissons pour $\left|\left\|.\right\|\right|$ la norme euclidienne. Nous allons imposer, à $f,$ $q$ conditions de Cauchy à coefficients affines $a^{j}_{m}(y)\;=\;a^{j}_{m}(x)\;+\;\sum^{n}_{i=1}d^{j,i}_{m}(y_{i}\;-\;x_{i}),$ avec $d^{j,i}_{m}\in A$ vérifiant les conditions: 
\begin{equation}\label{cex4}
\begin{split}
&d^{j,i}_{m}\;+\;d^{i,j}_{m}\;=\;0\quad si\quad m=1,...,q;\quad et \quad i\neq j,\\
&d^{j,j}_{m}\;=\;d_{m},
\end{split}
\end{equation}
(d'autres conditions seront encore imposées ci-dessous), et chercher $\varphi_{m}(x,y)$ sous la forme d'un polynôme de degré 1 : $\varphi_{m}(x,y)\;=\;\sum^{n}_{i=1}b^{i}_{m}(x)(y_{i}-x_{i})$  (rappelons que $\varphi_{m}$ doit s'annuler sur la diagonale). \\
En tenant compte de \ref{C5}, les conditions \ref{C4} se réduisent à:
\begin{equation}
\begin{cases}
\sum^{q}_{m=1}a^{j}_{m}(x)b^{j}_{m}(x)\;=\;\frac{e_{0}}{nVol(B(0,1))}\\
\sum^{q}_{m=1}a^{j}_{m}(x)b^{i}_{m}(x)\;+\;a^{i}_{m}(x)b^{j}_{m}(x)\;=0,\quad si\quad i\neq j,\\
\left\|y-x\right\|^{2}\sum_{m,j,i}d^{j,i}_{m}(y_{i}-x_{i})b^{j}_{m}(x)\;=\;n\sum^{q}_{m=1}\sum^{n}_{j,i,l=1}d^{j,l}_{m}(y_{l}-x_{l})(y_{j}-x_{j})(y_{i}-x_{i})b^{i}_{m}(x).
\end{cases}
\end{equation}
Les deux premières lignes, comme nous l'avons expliqué ci-dessus, nous conduisent aux conditions (A) et à la solution fournie par \ref{cond8}. En identifiant, dans la dernière ligne, les polynômes dans chaque membre, nous obtenons: 
\begin{equation}
\begin{cases}
\sum^{q}_{m=1}\sum^{n}_{l=1}d^{l,j}_{m}b^{l}_{m}\;=\;\sum_{m}d^{j,j}_{m}b^{j}_{m}, \quad pour \quad tout\quad j=1,...,n;\\
n\sum_{m}d^{j,j}_{m}b^{i}_{m}\;+\;d^{i,j}_{m}b^{j}_{m}\;+\;d^{j,i}_{m}b^{j}_{m}\;=\;\sum_{m,l}d^{l,i}_{m}b^{l}_{m},\quad pour\quad i\neq j;\\
\sum_{m}b^{l}_{m}(d^{i,j}_{m}+d^{j,i}_{m})\;+\;b^{j}_{m}(d^{i,l}_{m}+d^{l,i}_{m})\;+\;b^{i}_{m}(d^{j,l}_{m}+d^{l,j}_{m})\;=\;0,\quad si\quad i,j\;et\;l\;tous\;\;distincts.
\end{cases}
\end{equation}
Nous pouvons remplacer la deuxième équation par l'équation obtenue en lui soustrayant la première. Nous obtenons ainsi: 
$$\sum_{m}(d^{i,i}_{m}-d^{j,j}_{m})b^{i}_{m}\;+\;(d^{i,j}_{m}+d^{j,i}_{m})b^{j}_{m}\;=\;0,\quad i\neq j.$$
D'après \ref{cex4}, cette équation est satisfaite, ainsi que la troisième équation du système. Il ne reste donc plus que la première. D'après \ref{cond8}, celle-ci s'écrit: 
$$\sum_{m,l}d^{l,i}_{m}\widetilde{D^{l}_{m}}\;=\;n\sum_{m}d_{m}\widetilde{D^{i}_{m}},\quad pour\quad i=1,...,n,$$
ou encore:
$$\sum_{m,l\neq i}d^{l,i}_{m}\widetilde{D^{l}_{m}}\;=\;(n-1)\sum_{m}d_{m}\widetilde{D^{i}_{m}},\quad pour\quad i=1,...,n.$$
Supposons que l'un (au moins) des $a^{j}_{m}(x)$ est inversible, pour tout $x\in D,\quad a^{1}_{1},$ par exemple. Alors, puisque $D_{0}\;=\;\sum^{q}_{l=1}(-1)^{l+1}a^{1}_{l}\widehat{D^{l}_{1}},$ et puisque $\widetilde{D^{1}_{1}}\;=\;\sum^{q}_{l=1}(-1)^{l+1}c^{l}_{1}\widehat{D^{l}_{1}}\;=\;\frac{\widehat{D^{1}_{1}}}{n}\;+\;\sum^{q}_{l=2}(-1)^{l+1}c^{l}_{1}\widehat{D^{l}_{1}},$ il est possible de choisir les éléments $(c^{2}_{1},c^{3}_{1},...,c^{q}_{1})$ qui sont arbitraires, alors que $c^{1}_{1},$ lui, ne l'est pas (voir \ref{cond4}), de telle façon que $\widetilde{D^{1}_{1}}$ soit inversible. Alors, pour tout $i=2,...,n,$ la condition ci-dessus s'écrit: 
$$d^{1,i}_{1}\widetilde{D^{1}_{1}}+\sum_{2\leq l\neq i}d^{l,i}_{1}\widetilde{D^{l}_{1}}+\sum_{m\geq 2,l}d^{l,i}_{m}\widetilde{D^{l}_{m}}\;=\;(n-1)\sum_{m}d_{m}\widetilde{D^{i}_{m}},$$
et, puisque $\widetilde{D^{1}_{1}}$ est inversible, on peut choisir $d^{1,i}_{1}$ pour que cette équation soit satisfaite. Il reste encore l'équation pour $i=1.$ Elle s'écrit $(n-1)\sum_{m}d_{m}\widetilde{D^{1}_{m}}\;=\;\sum_{l\neq i,m}d^{l,i}_{m}\widetilde{D^{l}_{m}},$ et, pour la même raison, on peut choisir $d_{1}$ pour qu'elle soit vérifiée.

\begin{ex}
Un exemple dans le cas non commutatif.
\end{ex}
Ici aussi, il semble difficile, dans le cas non commutatif, de mener à terme les calculs dans le cas général. Nous donnerons donc un exemple.\\
Les conditions pour l'obtention d'une représentation intégrale sont données, rappelons le, en \ref{C2}. Nous choisirons encore la norme euclidienne, et supposerons, dans cet exemple, que $f$ est une fonction d'un ouvert connexe $D$ de $A$ dans $A$ à laquelle nous imposons $p$ conditions de Cauchy 
$$\sum^{p}_{j=0}a^{j}_{m}(x)\frac{\partial f}{\partial x_{j}}(x)\;=\;0,\quad pour\quad m=1,...,p.$$
avec tous les coefficients $a^{j}_{m}(x)$ nuls, sauf $a^{0}_{m}(x)$ et $a^{m}_{m}(x).$ Nous cherchons $\varphi_{m}(x,y)$ sous la forme: 
$$\varphi_{m}(x,y)\;=\;\sum^{p}_{j=0}b^{j}_{m}(x)(y_{j}-x_{j}).$$
Les conditions \ref{C2} s'écrivent: 
\begin{equation}
\begin{cases}
\sum^{p}_{m=1}a^{0}_{m}(x)b^{0}_{m}(x)\;+\;a^{m}_{m}(x)b^{m}_{m}(x)\;&=\;\frac{e_{0}}{Vol(B(0,1))},\\
\sum^{p}_{m=1}a^{0}_{m}(y)b^{0}_{m}(x)\;+\;a^{m}_{m}(y)b^{m}_{m}(x)\;&=\;\frac{p+1}{\left\|y-x\right\|^{2}}\sum_{m}[a^{0}_{m}(y)(y_{0}-x_{0})\;\\
&\qquad +\;a^{m}_{m}(y)(y_{m}-x_{m})]\sum^{p}_{j=0}b^{j}_{m}(x)(y_{j}-x_{j}).
\end{cases}
\end{equation}
En multipliant la deuxième équation par $\left\|y-x\right\|^{2},$ nous obtenons:
\begin{equation}
\begin{split}
\left\|y-x\right\|^{2}\sum_{m}(a^{0}_{m}(y)b^{0}_{m}(x)+&a^{m}_{m}(y)b^{m}_{m}(x))\;=\;(p+1)(y_{0}-x_{0})^{2}\sum_{m}a^{0}_{m}(y)b^{0}_{m}(x)\\
&+(p+1)\sum^{p}_{j=1}(y_{0}-x_{0})(y_{j}-x_{j})\sum_{m}a^{0}_{m}(y)b^{j}_{m}(x)\\
&+(p+1)\sum_{m}(y_{0}-x_{0})(y_{m}-x_{m})a^{m}_{m}(y)b^{0}_{m}(x)\\
&+(p+1)\sum_{m}(y_{m}-x_{m})^{2}a^{m}_{m}(y)b^{m}_{m}(x)\\
&+(p+1)\sum_{1\leq j<m}(y_{j}-x_{j})(y_{m}-x_{m})(a^{m}_{m}(y)b^{j}_{m}(x)+a^{j}_{j}(y)b^{m}_{j}(x)).
\end{split}
\end{equation}
Nous identifions les termes d'ordre 2 en $y-x$ dans les deux membres de cette égalité, et, tenant compte de la première équation du système précédent, nous obtenons: 
\begin{equation}
\begin{cases}
\begin{split}
\sum_{m}(a^{0}_{m}(x)b^{0}_{m}(x)+a^{m}_{m}(x)b^{m}_{m}(x))\;&=\;(p+1)\sum_{m}a^{0}_{m}(x)b^{0}_{m}(x)\\
&=\;(p+1)a^{j}_{j}(x)b^{j}_{j}(x)\;=\;\frac{e_{0}}{Vol(B(0,1))},\quad pour\quad j=1,...,p;
\end{split}\\
a^{m}_{m}(x)b^{j}_{m}(x)+a^{j}_{j}(x)b^{m}_{j}(x)\;=\;0,\quad si\quad 1\leq j\neq m\geq 1;\\
\sum_{m}a^{0}_{m}(x)b^{j}_{m}(x)\;+\;a^{j}_{j}(x)b^{0}_{j}(x)\;=\;0,\quad pour\quad j=1,...,p;
\end{cases}
\end{equation}
ce qui est équivalent à:
\begin{equation}
\begin{cases}
a^{j}_{j}(x)b^{j}_{j}(x)\;=\;\frac{e_{0}}{(p+1)Vol(B(0,1))}\;=\;\sum_{m}a^{0}_{m}(x)b^{0}_{m}(x)\quad pour\quad j=1,...,p;\\
a^{m}_{m}(x)b^{j}_{m}(x)+a^{j}_{j}(x)b^{m}_{j}(x)\;=\;0,\quad si\quad 1\leq j\neq m\geq 1;\\
a^{j}_{j}(x)b^{0}_{j}(x)\;=\;-\sum_{m}a^{0}_{m}(x)b^{j}_{m}(x),\quad pour\quad j=1,...,p.
\end{cases}
\end{equation}
De la première égalité dans la première équation du système, nous déduisons que $a^{j}_{j}(x)$ et $b^{j}_{j}(x)$ sont inversibles, et que $b^{j}_{j}(x)=\frac{(a^{j}_{j}(x))^{-1}}{(p+1)Vol(B(0,1))}.$ De la deuxième équation, nous obtenons $b^{j}_{m}(x)=-(a^{m}_{m}(x))^{-1}a^{j}_{j}(x)b^{m}_{j}(x),$ et de la troisième $b^{0}_{j}(x)=-(a^{j}_{j}(x))^{-1}\sum_{m}a^{0}_{m}(x)b^{j}_{m}(x).$\\
Il reste la deuxième égalité de la première équation, dont nous déduisons la condition suivante qui, dans le contexte de notre exemple, est l'équivalent de la condition $(A):$ 
\begin{equation}
\begin{split}
\frac{e_{0}}{(p+1)Vol(B(0,1))}\;&=\;\sum_{m}a^{0}_{m}(x)b^{0}_{m}(x)\\
&=\;-\sum_{m}a^{0}_{m}(x)(a^{m}_{m}(x))^{-1}\sum^{p}_{j=1}a^{0}_{j}(x)b^{m}_{j}(x)\\
&=\;-\sum_{m}a^{0}_{m}(x)(a^{m}_{m}(x))^{-1}a^{0}_{m}(x)b^{m}_{m}(x)\\
&\quad-\,\sum_{1\leq j<m}a^{0}_{m}(x)(a^{m}_{m}(x))^{-1}a^{0}_{j}(y)b^{m}_{j}(x)\;-\;\sum_{m<j}a^{0}_{m}(x)(a^{m}_{m}(x))^{-1}a^{0}_{j}(x)b^{m}_{j}(x)\\
&=\;-\sum_{m}\frac{a^{0}_{m}(x)(a^{m}_{m}(x))^{-1}a^{0}_{m}(x)(a^{m}_{m}(x))^{-1}}{(p+1)Vol(B(0,1))}\\
&\quad-\;\sum_{j<m}\big[a^{0}_{m}(x)(a^{m}_{m}(x))^{-1}a^{0}_{j}(x)\;\\
&\qquad \qquad -\;a^{0}_{j}(x)(a^{j}_{j}(x))^{-1}a^{0}_{m}(x)(a^{m}_{m}(x))^{-1}a^{j}_{j}(x)\big]b^{m}_{j}(x)\\
&=\;-\sum_{m}\frac{a^{0}_{m}(x)(a^{m}_{m}(x))^{-1}a^{0}_{m}(x)(a^{m}_{m}(x))^{-1}}{(p+1)Vol(B(0,1))}\\
&\quad-\;\sum_{j<m}\big[a^{0}_{m}(x)(a^{m}_{m}(x))^{-1}a^{0}_{j}(x)(a^{j}_{j}(x))^{-1}\\
&\quad-\;a^{0}_{j}(x)(a^{j}_{j}(x))^{-1}a^{0}_{m}(x)(a^{m}_{m}(x))^{-1}\big]a^{j}_{j}(x)b^{m}_{j}(x),
\end{split}
\end{equation}
ce que nous pouvons aussi écrire:
\begin{equation}
\begin{split}
e_{0}&+\sum^{p}_{m=1}\big[a^{0}_{m}(x)(a^{m}_{m}(x))^{-1}\big]^{2}\in\\
&\in I\Big[\Big(a^{0}_{m}(x)(a^{m}_{m}(x))^{-1}\Big)\Big(a^{0}_{j}(x)(a^{j}_{j}(x))^{-1}\Big)-\Big(a^{0}_{j}(x)(a^{j}_{j}(x))^{-1}\Big)\Big(a^{0}_{m}(x)(a^{m}_{m}(x))^{-1}\Big):1\leq j<m\Big]
\end{split}
\end{equation}
où $I\Big[a(j,m)\;:\;1\leq j<m\Big]$ désigne l'idéal de $A$ engendré par les éléments $a(j,m)\;:\;1\leq j<m$ de $A.$\\
Remarquons, encore une fois, que cette condition porte sur la structure d'algèbre de $A.$\\
\\
\section{Annexe : quelques exemples exotiques}
Ici, nous allons donner des exemples de formule de Cauchy pour des alg\`ebres $A$ plus exotiques.

Le premier concerne l'alg\`ebre des Tessarines qui peut \^etre obtenue par la donn\'ee d'une base $(1,i,j,k)$ v\'erifiant $i^2=-1$, $j^2=1$, $k=ij$, et le reste de la table de multiplication s'obtient en imposant l'associativit\'e et la commutativit\'e. En utilisant ceci, $k^2=-1$. Maintenant consid\'erons $f: A\rightarrow  A$ ; $A$ \'etant commutative, en it\'erant la preuve de la non existence en dimension $3$ et dans le cas commutatif d'une formule de Cauchy  satisfaisant une seule \'equation de Cauchy, nous obtenons la non existence d'une formule de Cauchy, pour $f$, dans le cas des Tessarines avec une seule condition dans l'alg\`ebre  (voir preuve dans le cas commutatif du th\'eor\`eme $2.12$). Par contre, pour les fonctions $A$-diff\'erentiables, nous avons une formule de Cauchy  d'apr\`es le corollaire 2.6 car $1+i^2+j^2+k^2=0$.

Pour l'obtention d'une formule de Cauchy pour les fonctions de $A$ dans $A$ dans le cas d'une seule condition de Cauchy dans l'alg\`ebre $A$, le caract\'ere associatif impose des contraintes assez strictes : nous avons des exemples, bien s\^ur dans $\C$, et dans  l'alg\`ebre des quaternions munie de la condition de Cauchy-Fueter (voir section $3$). Si l'on omet cette condition, nous avons plus d'exemples et nous allons montrer comment les g\'en\'erer. Si $A$ est une alg\`ebre unitaire ou unif\`ere de dimension $p$ ayant pour base $(1,e_2,\cdots ,e_p)$, alors les conditions \ref{cond12} avec $m=1$ sont v\'erifi\'ees si la base v\'erifie $e_{i}^2=-1$ et $e_ie_j=-e_je_i$ pour tous les $i,j>1$ et $i\not =j$, et nous avons donc une formule de Cauchy pour les fonctions $f$ satisfaisant la  condition de Cauchy :
$$\frac{\partial f} {\partial x_1}+e_2\frac{\partial f}{\partial x_2}+ \cdots e_p\frac{ \partial f}{ \partial x_p}=0.$$
Ces conditions sont satisfaites, par exemple, avec $A$ les octonions ou les s\'ed\'enions munis de leur base standard. 

Enfin, pour finir, dans le cas d'une condition de Cauchy ($m=1$) alors le noyau $K$ est une solution fondamentale et nous avons donc beaucoup de fonctions $f$ de $A$ dans $A$ v\'erifiant la condition de Cauchy. Si $m>1$, ce n'est plus le cas, et la classe de fonction v\'erifiant les conditions de Cauchy peut \^etre plus restreinte. C'est exactement le ph\'enom\`ene que nous rencontrons dans $\C$ et $\C^n$ pour les fonctions d'une et plusieurs variables complexes, en une variable le noyau de Cauchy est holomorphe, en plusieurs variables, le noyau de Bochner ne l'est plus.


\begin{thebibliography}{<00>}

\bibitem[BM1]{BM1}
P.BONNEAU et E.MAZZILLI,
\newblock Complex associated to some systems of PDE,
\newblock J. Geom. Phys. 169 (2021), Paper No. 104319, 28 pp.

\bibitem[WW1]{WW1}
W.WANG,
\newblock On the non-homogeneous Cauchy-Fueter equations and Hartog’s phenomenon in several quater- nionic variables, 
\newblock Journal of Geometry and Physics, 58 (2008), 1203-1210.

\bibitem[WW2]{WW2}
W.WANG
\newblock The k-Cauchy-Fueter complex, Penrose transformation and Hartogs phenomenon for quater-nionic k-regular functions, 
\newblock Journal of Geometry and Physics, 60 (2010), 513-530.

\bibitem[BCGGG]{BCGGG}
R. Bryant, S. Chern, R. Gardner, H. Goldschmidt, P. Griffiths,
\newblock Exterior differential systems, 
\newblock Springer Verlag, Berlin (1991).


\bibitem[BC1]{BC1}
P. BONNEAU et A. CUMENGE,
\newblock  Entre analyse complexe et superanalyse,
\newblock C. R. Acad Sci. Paris, Sér. 1 (2009).



\bibitem[BC2]{BC2}
P. BONNEAU et A. CUMENGE,
\newblock  A mi-chemin entre analyse complexe et superanalyse,
\newblock Publ. Mat. 56 (2012), 3-40.





\bibitem[BM]{BM}
R. BOTT and J. MILNOR,
\newblock On the parallelizability of the spheres,
\newblock     Bull. AMS, {\bf 64},  (1958), 87-89.



\bibitem[HP]{HP}
R. HARVEY and J. POLKING,
\newblock Fundamental solutions in complex analysis, Parts 1 and 2,

\newblock Duke Math. J., Vol 46, 253-340 (1978).

\bibitem[H]{H}
L. HÖRMANDER,
\newblock , Complex Analysis in Several Variables
\newblock North-Holland, Amsterdam (1973).


\bibitem[K]{K}
M. KERVAIRE,
\newblock Non-parallelizability of the n-sphere for n>7,
\newblock     Proc. Nat. Acad. Sci. USA, {\bf 44},  (1958), 280-283.



\bibitem[LS]{LS}
S. LANG,
\newblock , Algebra,
\newblock Addison-Wesley, Amsterdam (1965).






\bibitem[L]{L}
P. LELONG,
\newblock Fonctions plurisousharmoniques et fonctions analytiques de variables réelles,
\newblock Ann. Inst. Fourier, Vol. 11, 515-562 (1961).



\bibitem[N]{N}
NARASIMHAN,
\newblock Analysis on real and complex manifolds,
\newblock     Masson,  Paris.









\bibitem[R]{R}
R. M. RANGE, 
\newblock Holomorphic functions and integral representations in several complex variables, 
\newblock Springer-Verlag, Berlin (1986)



\end{thebibliography}
\end{document}